\documentclass[12pt]{mythesis}
\title{Jones Pairs}
\author{Ada Chan}
\discipline{Combinatorics and Optimization}
\degree{Doctor of Philosophy}

\usepackage{amsmath}
\usepackage{amssymb}
\usepackage{theorem}
\usepackage{enumerate}
\usepackage{epsfig}
\usepackage{makeidx}
\usepackage{float}

\newtheorem{theorem}{Theorem}[section]
\newtheorem{corollary}[theorem]{Corollary}

\newtheorem{lemma}[theorem]{Lemma}

{\theorembodyfont{\rmfamily} \newtheorem{example}[theorem]{Example}}

\newcommand{\eop}{\nobreak\ \hfill $\Box$ \medbreak} 
\newcommand{\SL}[1]{{\sl {#1}}{\index{#1}}}        

\def\N#1#2{{\cal N}_{#1,#2}}        
\def\dualN#1#2{{\cal N}'_{#1,#2}}   
\def\Nom#1{{\cal N}_{#1}}           
\def\dualNom#1{{\cal N}'_{#1}}      
\def\XDX#1#2#3{X_{#1} \Delta_{#2} X_{#3}}       
\def\DXD#1#2#3{\Delta_{#1} X_{#2} \Delta_{#3}}  
\def\A{{\cal A}}                    
\def\BM{{\cal B}}                   
\def\Ter#1{{\cal T}_{#1}}           
\def\Br#1{{\mathbf{B_{#1}}}}        
\def\Sinv#1{{#1}^{(-)}}             
\def\SinvT#1{{#1}^{(-)T}}           
\def\inv#1{{#1}^{-1}}               
\def\invT#1{{#1}^{-T}}              
\def\schur{\circ}
\def\com{\mathbb{C}}                
\def\D{\Delta} 
\def\T#1#2{\Theta_{#1,#2}}          
\def\Evector#1#2#3#4{#1 e_{#2} \schur #3 e_{#4}} 
\def\Eij{E_{ij}}                    
\def\range#1#2#3{#1=#2, \ldots, #3}  
\def\seq#1#2#3{#1_{#2},\ldots,#1_{#3}}
\def\mat#1#2{{\mathbf{M}_{#1}(#2)}} 
\def\nbyn{n \times n}               
\def\tr{\mathop{\hbox{\rm tr}}\nolimits}         
\def\spn{\mathop{\hbox{\rm span}}\nolimits}     
\def\elsm{\mathop{\hbox{\rm sum}}\nolimits}      
\def\ip#1#2{{\langle #1, \ #2 \rangle}}          

\makeindex

\begin{document}
\prepages
\maketitle
\sigpages
\begin{abstract}

Motivated by Jones' braid group representations constructed from
spin models, we define {\sl a Jones pair} to be a pair of 
$\nbyn$ matrices $(A,B)$ such that the endomorphisms $X_A$ and $\D_B$
form a representation of a braid group.  
When $A$ and $B$ are type-II matrices, we call $(A,B)$ 
{\sl an invertible Jones pair}.
We develop the theory of Jones pairs in this thesis.

Our aim is to study the connections among association schemes,
spin models and four-weight spin models using the viewpoint of Jones pairs.
We use Nomura's method to construct a pair of algebras from the 
matrices $(A,B)$, which we call the Nomura algebras of $(A,B)$.
These algebras become the central tool in this thesis.  
We explore their properties in Chapters~\ref{Nomura} and \ref{IINom}.

In Chapter~\ref{JP}, we introduce Jones pairs.  
We prove the equivalence of four-weight spin models and invertible
Jones pairs.
We extend some existing concepts for four-weight spin models to Jones pairs.
In Chapter~\ref{SpinModels}, 
we provide new proofs for some well-known results on the Bose-Mesner algebras 
associated with spin models.

We document the main results of the thesis in Chapter~\ref{InvJP}.
We prove that every four-weight spin model
comes from a symmetric spin model (up to odd-gauge equivalence).
We present four Bose-Mesner algebras associated to each
four-weight spin model.  We study the relations among these algebras.
In particular,
we provide a strategy to search for four-weight spin models.
This strategy is analogous to the method given by
Bannai, Bannai and Jaeger for finding spin models.

\end{abstract}

\begin{acknowledgements}

First of all, I wish to thank Chris Godsil for being a great mentor.  
He has introduced me to the fascinating area of algebraic combinatorics and 
given me excellent guidance on my research.
I am very grateful to the amount of time and energy he has spent on
helping me grow as a mathematician.
His advice on doing research, giving presentations, 
writing and teaching will always be helpful.

I feel very fortunate that I was a part of the Department of Combinatorics and 
Optimization.  
The professors and staffs have provided a comfortable learning environment.
In particular, I would like to thank Ian Goulden and Adrian Lewis for
being my lecturing requirements examiners.
I am also grateful to Ian Goulden, Sylvia Hobart, 
Peter Hoffman and David Wagner for being in my examining committee.

Without the love and the unwavering support from
my parents and my husband, I would never have accomplished
this work.  Lo-Dou, Mummy, Stewart, thank you for always believing
in me.

I appreciate the kind words of encouragement and valuable help from 
Miguel Anjos, Lise Arseneau, Mike Best, Eva Chan, Daniel Chan, Jason Chen, 
Kevin Cheung, Michael Cheung, Philip Fong, Jackie Ho, Peter Hoffman, 
Julia Huang, Sa$\check{\mathrm s}$a Juri$\check{\mathrm s}$i\'c, Florence Kong,
Kai Kwok, Serge Kruk, Charles Lam, Annie Lee, Patrick Lau, Bill Martin,
Roxana and Stanley Poon.  Thank you very much, my friends.

The most important of all, I admire God for creating the beauty in Mathematics.
I am grateful for having a chance to get a glimpse of this beauty.
\end{acknowledgements}

\newpage
\begin{center}
\vspace*{2.5in}
To Hing, Fung and Stewart
\end{center}
\tableofcontents
\mainbody
\chapter{Introduction}
\label{Intro}

We first give an overview of this thesis.
In Sections~\ref{Intro_Spin} and \ref{Intro_Braids},
we give the background materials on 
spin models and braids.  Then we present a historical overview of 
the research on the association schemes attached to spin models.

\section{Overview}
\label{Intro_Intro}

The purpose of this thesis is to introduce Jones pairs and
to extend the existing theory of association schemes attached
to four-weight spin models.

We now define Jones pairs.
Given two $\nbyn$ matrices $M$ and $N$, their Schur product 
$M \schur N$ is defined by
\begin{equation*}
(M \schur N)_{i,j} = M_{i,j} N_{i,j},
\end{equation*}
for all $\range{i,j}{1}{n}$.
If $C$ is an $\nbyn$ matrix, we define two endomorphisms of
$\mat{n}{\com}$, $X_C$ and $\D_C$, by
\begin{equation*}
X_C(M) = CM \quad {\text and} \quad \D_C(M)=C\schur M,
\end{equation*}
for any $\nbyn$ matrix $M$.
We say that the pair of $\nbyn$ matrices $(A,B)$ is a \SL{Jones pair}\  
if $X_A$ and $\D_B$ are invertible, and they satisfy
\begin{eqnarray*}
\XDX{A}{B}{A} &=& \DXD{B}{A}{B}, \quad {\text and}\\
\XDX{A}{B^T}{A} &=& \DXD{B^T}{A}{B^T}.
\end{eqnarray*}
We will show in this thesis that Jones pairs give representations of
braid groups.
Moreover, we will see that spin models and four-weight spin models belong
to a special class of Jones pairs, called the invertible Jones pairs.
Consequently, we obtain a representation of braid group from every
four-weight spin model.  This fact was not known previously.

The connections between spin models and association schemes
have been the focus of existing research.
We have always found association schemes intriguing,
mainly because of their connections to a vast number of combinatorial
objects such as distance regular graphs, codes and designs.
It follows naturally that we are interested in the results due to
Jaeger \cite{J_Class} and Nomura \cite{N_Alg} which say that 
every spin model belongs to the Bose-Mesner algebra of some 
association scheme.
In \cite{J_Class}, Jaeger asked for an intrinsic characterization
of the association schemes whose Bose-Mesner algebras contain a spin model
and this is the question which motivates the work in this thesis.

Nomura used the type-II condition of spin models to obtain the 
result mentioned above.  
We say an $\nbyn$ matrix $A$ satisfies the \SL{type-II condition}\  if
\begin{equation*}
\sum_{x=1}^n \frac{A_{i,x}}{A_{j,x}} = \delta_{ij}n, \quad
\text{ for all }\range{i,j}{1}{n}
\end{equation*}
and it is called a \SL{type-II matrix}.
Spin models and the matrices in four-weight spin models are examples of
type-II matrices.
In \cite{N_Alg}, Nomura constructed a Bose-Mesner algebra of some
association scheme from each type-II matrix.  
This algebra is called the Nomura algebra of the type-II matrix.  
He also showed that every spin model belongs to its Nomura algebra.
In Chapter~\ref{Nomura}, we generalize his construction to a pair of
algebras built from a pair of matrices $(A,B)$.
We also study the properties of this pair of algebras.
In this process, we find new and simpler proofs of
many known results on the Nomura algebras of type-II matrices and spin models.
We are now convinced that Jones pairs provide a natural setting for 
the study of the problems related to spin models and four-weight spin models.

Relying on the fact that every spin model belongs to its Nomura algebra,
Bannai, Bannai and Jaeger \cite{BBJ} designed a strategy to search 
for spin models.
However, the matrices in a four-weight spin model do not belong to their
Nomura algebras.
So Bannai, Bannai and Jaeger's method does not apply directly to 
four-weight spin models.
In Chapter~\ref{InvJP}, we provide a construction of 
an $4n \times 4n$ symmetric spin model from each $\nbyn$ four-weight spin model.
This construction generalizes a construction due to Nomura in \cite{N_Twist}.
As a result, four-weight spin models are not very different from symmetric
spin models.
Moreover, we design a strategy to find four-weight spin models using
this newly constructed $4n \times 4n$ symmetric spin model.

In addition to the Nomura algebras obtained from the matrices in 
an $\nbyn$ four-weight spin models and 
the $4n \times 4n$ symmetric spin model, we will construct 
two more Nomura algebras
from the four-weight spin model. 
In Chapter~\ref{InvJP}, we will see how
these Nomura algebras are intricately related to each other.

\section{Spin Models}
\label{Intro_Spin}

Using a statistical mechanical model called a spin model, Jones
\cite{VJ_Knot}
constructed invariants of unoriented links in the form of partition functions.
His definition of a spin model is essentially a symmetric matrix that
satisfies certain properties such that its partition function 
is invariant under the Reidemeister moves of link diagrams.
The Potts model, which is a linear combination of the identity
matrix and the matrix of all ones, is the simplest spin model.
From the Potts model,
we can obtain the infamous Jones polynomial.

In 1994, Kawagoe, Munemasa and Watatani \cite{KMW_Gen} generalized 
the definition of spin models by removing the symmetry condition.  
We adopt their definition in this thesis:
an $\nbyn$ matrix $W$ is a \SL{spin model}\ if there exists a non-zero
scalar $a$ and $d=\pm \sqrt{n}$ such that
\begin{enumerate}[(I)]
\item
For $\range{i}{1}{n}$, 
\begin{equation*}
W_{i,i} = a
\end{equation*} 
and
\begin{equation*}
\sum_{x=1}^n W_{i,x} = \sum_{x=1}^n W_{x,i} = da^{-1}.
\end{equation*}
\item
For all $\range{i,j}{1}{n}$, 
\begin{equation*}
\sum_{x=1}^{n} \frac{W_{i,x}}{W_{j,x}} = 
\delta_{ij} n,
\end{equation*}
where $\delta_{ij}$ equals one when $i=j$ and zero otherwise.
\item
For all $\range{i,j,k}{1}{n}$,
\begin{equation*}
\sum_{x=1}^n  \frac{W_{k,x}W_{x,i}}{W_{j,x}} =
d\frac{W_{k,i}}{W_{j,i}W_{j,k}}.
\end{equation*}
\end{enumerate}
They also showed that the partition functions give invariants of oriented links.  
So far, we only know three infinite families of spin models.
The first family is the Potts models.  The second one comes from finite
Abelian groups.  In \cite{BBJ}, Bannai, Bannai and Jaeger built 
a spin model from each finite Abelian group.
The third family consists of the symmetric and the non-symmetric Hadamard
spin models, constructed by Jaeger and Nomura \cite{N_Had, JN}.

It turns out that the type-II condition for spin models 
is more important than the other two conditions.  
In particular, the type-II condition of a spin model
plays a key role in Nomura's construction of the Bose-Mesner algebra
containing it.
Moreover, unitary type-II matrices are important objects in the study of
Von~Neumann algebras \cite{JS}.

In 1995, Bannai and Bannai \cite{BB_4wt} gave a further
generalization by defining the four-weight spin models.  
In \cite{J_4wt}, Jaeger normalized the partition function of a
four-weight spin model to give an invariant of oriented links.
A \SL{four-weight spin model}\ is a $5$-tuple 
$(W_1,W_2,W_3,W_4;d)$ with $d^2=n$ and a non-zero scalar $a$ satisfying
\begin{enumerate}[(I)]
\item
For all $\range{\alpha}{1}{n}$,
\begin{equation*}
(W_3)_{\alpha,\alpha} = a^{-1},\quad (W_1)_{\alpha,\alpha}=a,
\end{equation*}
and
\begin{eqnarray*}
\sum_{x=1}^n (W_2)_{\alpha,x} =& \sum_{x=1}^n (W_2)_{x,\alpha}& =
da^{-1},\\
\sum_{x=1}^n (W_4)_{\alpha,x} =& \sum_{x=1}^n (W_4)_{x,\alpha}& = da.
\end{eqnarray*}
\item
For all $\range{\alpha,\beta}{1}{n}$,
\begin{eqnarray*}
(W_1)_{\alpha,\beta} (W_3)_{\beta,\alpha} = 1, 
&\quad & (W_2)_{\alpha,\beta} (W_4)_{\beta,\alpha} = 1, \\
\sum_{x=1}^n (W_1)_{\alpha,x} (W_3)_{x,\beta} &=& 
\delta_{\alpha \beta} n,\\
\sum_{x=1}^n (W_2)_{\alpha,x} (W_4)_{x,\beta} &=& 
\delta_{\alpha \beta} n.
\end{eqnarray*}
\item

For all $\range{\alpha,\beta,\gamma}{1}{n}$,
\begin{eqnarray*}
\sum_{x=1}^n (W_1)_{\alpha,x} (W_1)_{x,\beta} (W_4)_{\gamma,x} &=& 
d (W_1)_{\alpha,\beta} (W_4)_{\gamma,\alpha} (W_4)_{\gamma,\beta},\\
\sum_{x=1}^n (W_1)_{x,\alpha} (W_1)_{\beta,x} (W_4)_{x,\gamma} &=& 
d (W_1)_{\beta,\alpha} (W_4)_{\alpha,\gamma} (W_4)_{\beta,\gamma}.
\end{eqnarray*}
\end{enumerate}
In the same paper, Bannai and Bannai listed sixteen equations that are
equivalent to the type-III conditions of four-weight spin models.
Having to decide which equations are more suitable for a given problem
complicates any analysis of four-weight spin models.
One advantage of Jones pairs is that they save us from 
having to deal with these sixteen equations.

\section{Braids}
\label{Intro_Braids}

Jones pairs are designed to give representations of braid groups,
so braids are naturally the topological object of interest here.
In this section, we provide some essential background of braids and 
braid groups. 

A \SL{braid}\ on $m$ strands is a set of $m$ disjoint arcs in 3-space joining
$m$ points on a horizontal plane to $m$ points on another horizontal
plane directly below the first $m$ points.  
Given a braid, we form the \SL{closure of the braid}\ by joining
the $m$ points on top to the $m$ points at the bottom as illustrated in
the following figure.
It is obvious that the closure of a braid is a link.

\begin{figure}[h]
\centering\epsfig{file=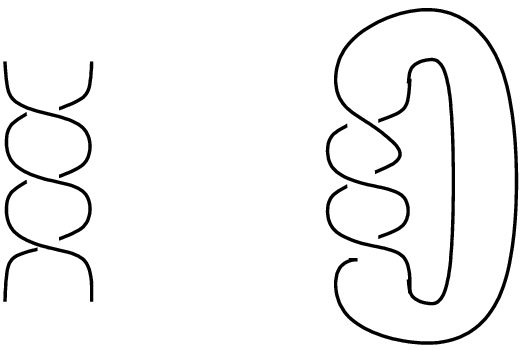}
\label{fig_Trefoil}
\end{figure}

The following theorem due to Alexander \cite{Murasugi} 
explains the connection between braids and links.
\begin{theorem}
\label{thm_Braid_Link}
Any link is isotopic to the closure of some braid.
\eop
\end{theorem}

When we stack two braids on $m$ strands, we get a new braid, which we call
the \SL{product of the two braids}.  This operation is associative.
If we define $\sigma_i$ and $\sigma_i^{-1}$ as shown in the figure below,
then it is easy to see that any braid on $m$ strands 
is a product of $\seq{\sigma}{1}{m-1}$.
Hence the inverse of a braid always exists.
As a result, we get a group structure on the set of braids on $m$ strands.

\begin{figure}[H]
\centering\epsfig{file=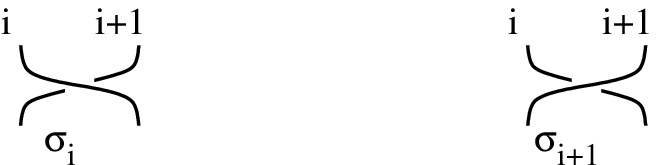}
\label{fig_GenBraid}
\end{figure}

The \SL{braid group on $m$ strands}, $\Br{m}$, 
is generated by $\seq{\sigma}{1}{m-1}$ subject to the following relations:
\begin{enumerate}[(a)]
\item
For all $|i-j| \geq 2 $,
\begin{equation*}
\sigma_i \sigma_j = \sigma_j \sigma_i.
\end{equation*}
\item
For all $\range{i}{1}{m-2}$,
\begin{equation*}
\sigma_i \sigma_{i+1} \sigma_i= \sigma_{i+1} \sigma_i \sigma_{i+1}.
\end{equation*}
\end{enumerate}


Given a spin model $W$, Jones \cite{VJ_Knot} built a representation of
the braid group $\Br{m}$.  His construction uses the endomorphisms,
$X_{W}$ and $\D_{\inv{W}}$, of $\mat{n}{\com}$.  
He showed that the trace of the element representing a 
braid equals the partition function invariant evaluated at the closure of the
braid, up to normalization.

In summary, given a spin model $W$, we can obtain the same link invariant
either using the partition function or the trace of a representation of the
braid group.  However Jones pointed out a puzzling distinction between
the two approaches, which is that the type-II condition on 
$W$ is needed for the first method, but not for the second.
This distinction motivates our definition of Jones pairs in
Chapter~\ref{JP}.
In particular, we do not assume the type-II condition on the matrices
in Jones pairs.
In the same chapter, 
we will generalize Jones' construction to four-weight spin models.

\section{Association Schemes}
\label{Intro_Schemes}

In 1992,
Jaeger \cite{J_SRG} made the first connection between association schemes 
and spin models.  
He showed that a spin model gives an evaluation of the
Kauffman polynomial if and only if it belongs to a triply-regular
formally self-dual two-class association scheme.  
Now two-class association schemes are equivalent to strongly regular graphs.
Moreover, a triply-regular two-class association scheme corresponds to
a strongly regular graph $G$ with the property that the neighbourhoods
of any vertex in $G$ and its complement induce strongly regular graphs.  
From this result, we see strong combinatorial properties
attached to spin models.  

Four years later, Jaeger \cite{J_Class} obtained the surprising
result that every spin model belongs to the Bose-Mesner algebra of
a formally self-dual association scheme.

Jaeger's profound discoveries caught the interests of 
researchers in Algebraic Combinatorics.  
They produced a series of results on spin models and association schemes.
Here, we document some other classical findings of this area.

Shortly after Jaeger announced that every spin model belongs to a
formally self-dual Bose-Mesner algebra, Nomura \cite{N_Alg}
came up with a substantially simpler algebraic construction of 
a Bose-Mesner algebra from a spin model. 
In fact, Nomura showed that each type-II matrix gives rise to a
Bose-Mesner algebra, now known as the Nomura algebra of the type-II
matrix.
Jaeger, Matsumoto and Nomura \cite{JMN} examined further the Nomura
algebras of type-II matrices,
and concluded that a type-II matrix $W$ belongs to its Nomura algebra  
if and only if it is a spin model up to scalar multiplication.

In 1997, Bannai, Bannai and Jaeger \cite{BBJ} found a necessary
condition for a formally self-dual Bose-Mesner algebra to be 
the Nomura algebra of some spin model.  They showed that the matrix of
eigenvalues of the Nomura algebra of a spin model must satisfy the
modular invariance property.
Unfortunately, we still do not have any way to tell whether a 
Bose-Mesner algebra is the Nomura algebra of some type-II matrix.

Much less is known about the connection between four-weight spin
models and association schemes.  
The only result we know is due to Bannai, Guo and Huang \cite{Ban1, Guo_4wt}.
They proved that if $(W_1,W_2,W_3,W_4;d)$ is a four-weight spin model, 
then the Nomura algebras of $W_1$, $W_2$, $W_3$ and $W_4$ coincide,
and this algebra is formally self-dual.   
\section{Directory}
\label{Dir}

In this section, we provide a layout of this thesis.

Chapter~\ref{Nomura} lays the groundwork for subsequent chapters.
In Sections~\ref{N_AB} and \ref{N_AB_XDX},
we introduce the Nomura algebras of two matrices $(A,B)$ and the duality map.
In Section~\ref{Exchange}, we present a useful tool called the Exchange Lemma.
These three sections are joint work by Godsil, Munemasa and the author
\cite{CGM}.
In Section~\ref{N_AB_Transform}, we examine effects on the Nomura algebras of
$(A,B)$ when $A$ and $B$ are multiplied by some monomial matrices.
In Section~\ref{N_AB_Tensor}, we show how the Nomura algebras of 
$(A_1\otimes A_2,\ B_1 \otimes B_2)$ are related to the Nomura algebras
of $(A_1,B_1)$ and $(A_2, B_2)$.  
The findings in these two sections 
generalize similar results from Jaeger, Matsumoto and Nomura \cite{JMN}.

Sections~\ref{Type_II}, \ref{N_A} and \ref{TypeII_Families}
survey existing results about the type-II  matrices
and their Nomura algebras.  
We present new proofs using the tools developed in Chapter~\ref{Nomura}
and they are due to Godsil, Munemasa and the author.
Sections~\ref{Scheme_BM_Alg} to \ref{Dual_Scheme} overview the standard
theory of association schemes.  
In Section~\ref{N_II}, we study the properties of the Nomura algebras of
two type-II matrices.  These properties are used substantially in subsequent chapters.  
The results in this section are joint work by Godsil, Munemasa and the
author.

We introduce one-sided Jones pairs and Jones pairs in Chapter~\ref{JP}.
Based on Jones' construction,
we build a braid group representation from a Jones pair.
More importantly, we show that invertible Jones pairs are equivalent
to four-weight spin models.  In Sections~\ref{Odd_Gauge} and
\ref{Even_Gauge}, we extend Jaeger's notion of gauge equivalence to
one-sided Jones pairs.  
The results in this chapter are joint work by Godsil, Munemasa and the author.

We examine spin models as Jones pairs  
in Chapter~\ref{SpinModels}. 
In Section~\ref{Mod_Inv}, we provide
Jaeger, Matsumoto and Nomura's derivation of the modular invariance
equation.  In Section~\ref{Hyper_Dual}, we extend
Curtin and Nomura's theorem to the strongly hyper-self-duality of
the Terwilliger algebra of the Nomura algebra of a spin model.
In Section~\ref{Two_Class}, we provide a new and shorter proof of one
direction of Jaeger's result on spin models from two class association
schemes.  The last section contains our new proof, using Jones pairs, 
to Jaeger and Nomura's construction of the symmetric and non-symmetric Hadamard spin models.

Chapter~\ref{InvJP} contains the main results of this thesis.
In Section~\ref{N_W}, we construct a type-II matrix $W$
from an invertible Jones pair $(A,B)$.
This type-II matrix gives a formally dual pair of Bose-Mesner algebras.
We determine the dimensions, the basis of Schur idempotents and 
the basis of principal idempotents for these algebras.
In Section~\ref{Ext}, we extend Nomura's construction \cite{N_Twist} to
build a pair of symmetric spin models, $V$ and $V'$, from 
each four-weight spin model.
We design an algorithm to exhaustively search for four-weight spin models, 
which is described in Section~\ref{JP_MI}.  
Sections~\ref{JP_Quotient} to \ref{JP_Subschemes} 
document our observation of the relations of 
the Nomura algebras of $A$, $W$, $V$ and $V'$.
Section~\ref{JP_Dim2} documents what we know about
the simplest class of Jones pairs.
Finally, we discuss several directions of future work.


\chapter{Nomura Algebras}
\label{Nomura}

In \cite{N_Alg}, Nomura constructed $n^2$ vectors from a symmetric
type-II matrix $A$.
He then considered the set of all matrices for which these vectors are 
eigenvectors,
and showed that this set forms a Bose-Mesner algebra.
This algebra is now called the Nomura algebra of $A$, and denoted by $\Nom{A}$.
Nomura went further and proved that if $A$ is a spin model, 
then $A$ belongs to its Nomura algebra.
Thus Nomura algebras play a significant role in the theory of spin
models.

During our investigation of Jones pairs, Godsil generalized Nomura's
construction.  Given a pair of matrices $A$ and $B$, we consider the set
of matrices of which for all $\range{i,j}{1}{n}$,
\begin{equation*}
\Evector{A}{i}{B}{j}
\end{equation*}
are eigenvectors.
We call this set of matrices the Nomura algebra of $(A,B)$, and denote it
by $\N{A}{B}$.
In general, the Nomura algebra of $(A,B)$ is not a Bose-Mesner algebra.
However it will become clear in subsequent chapters that $\N{A}{B}$ is a 
powerful tool in the study of Jones pairs.

Sections~\ref{N_AB} and \ref{N_AB_XDX} introduce the Nomura algebra of
$(A,B)$ and a special function on $\N{A}{B}$ called the duality map.
The results in these two sections provide technical background for
Chapters~\ref{IINom} to \ref{InvJP}.
They also appear in the paper \cite{CGM} by Godsil, Munemasa and
the author.

In Sections~\ref{N_AB_Transform} and \ref{N_AB_Tensor}, we
generalize several results due to Jaeger, Matsumoto and Nomura in
\cite{JMN}, which are in turn extensions of Nomura's results on
symmetric type-II matrices.

\section{Nomura Algebra of Two Matrices}
\label{N_AB}

We introduce the Nomura algebras of a pair of $\nbyn$ matrices
$(A,B)$, denoted by \SL{$\N{A}{B}$}\ and \SL{$\dualN{A}{B}$}.
We study the properties of these algebras and a map
\begin{equation*}
\T{A}{B}:\N{A}{B} \rightarrow \dualN{A}{B}
\end{equation*}called the duality map.
In this section, we resist the temptation of adding dispensable
conditions on $A$ and $B$.

For any pair of $\nbyn$ matrices $A$ and $B$, 
their \SL{Schur product $A \schur B$}\ is the entry-wise product of $A$ and $B$.  
That is, for $\range{i,j}{1}{n}$,
\begin{equation*}
(A \schur B)_{i,j} = A_{i,j} B_{i,j}.
\end{equation*}
The $\nbyn$ matrix of all ones, denoted by $J$, is the 
identity with respect to the Schur product.
If all entries of $A$ are non-zero, then the matrix $\Sinv{A}$
with
\begin{equation*}
\Sinv{A}_{i,j} = \frac{1}{A_{i,j}}
\end{equation*}
satisfies $A \schur \Sinv{A} = J$.
We say that $A$ is \SL{Schur invertible}\ and $\Sinv{A}$ is 
the \SL{Schur inverse of $A$}.
Note that $\mat{n}{\com}$ forms a commutative algebra with
respect to the Schur product.

Let $\{\seq{e}{1}{n}\}$ denote the standard basis for $\com^n$.
Given a pair of $\nbyn$ matrices $A$ and $B$,
we obtain the following set of $n^2$ vectors
\begin{equation*}
\{\Evector{A}{i}{B}{j} : \range{i,j}{1}{n}\}.
\end{equation*}
The vector $\Evector{A}{i}{B}{j}$ equals the Schur product of 
the $i$-th column of $A$ with the $j$-th column of $B$.
In most cases we encounter, these $n^2$ vectors are not distinct.

We define $\N{A}{B}$ to be the set of $\nbyn$ matrices for which 
$\Evector{A}{i}{B}{j}$ is an eigenvector, for all $\range{i,j}{1}{n}$.

\begin{lemma}
\label{lem_N_AB}
Let $A$ and $B$ be $\nbyn$ matrices.  
Then the set $\N{A}{B}$ is a vector space which
is closed under multiplication, and it contains the identity matrix.
\end{lemma}
{\sl Proof.}
Let $M, N \in \N{A}{B}$.
For $\range{i,j}{1}{n}$, there exist $m_{ij}$ and $n_{ij}$ such that
$M \ \Evector{A}{i}{B}{j} = m_{ij} \ \Evector{A}{i}{B}{j}$ and 
$N \ \Evector{A}{i}{B}{j} = n_{ij} \ \Evector{A}{i}{B}{j}$.
Now,
\begin{equation*}
(MN) \ \Evector{A}{i}{B}{j} = m_{ij} n_{ij} \ \Evector{A}{i}{B}{j}.
\end{equation*}
So we have $MN \in \N{A}{B}$.
It is immediate that $\Evector{A}{i}{B}{j}$ is an eigenvector for $I$,
for all $\range{i,j}{1}{n}$.  So $I\in\N{A}B{}$.
\eop

We call $\N{A}{B}$ the \SL{Nomura algebra for the pair $(A,B)$}.
Now for each matrix $M$ in $\N{A}{B}$, we use $\T{A}{B}(M)$ to denote
the $\nbyn$ matrix whose $ij$-entry is the eigenvalue of $M$
with respect to the eigenvector $\Evector{A}{i}{B}{j}$. 
That is, 
\begin{equation*}
M \ \Evector{A}{i}{B}{j} = \T{A}{B}(M)_{i,j} \ \Evector{A}{i}{B}{j},
\end{equation*}
for $\range{i,j}{1}{n}$.
We call \SL{$\T{A}{B}$}\ the \SL{duality map}. We use $\dualN{A}{B}$ to denote
the image of $\N{A}{B}$ under $\T{A}{B}$.
The following two lemmas tell us that $\dualN{A}{B}$ is an algebra
with respect to the Schur product.

\begin{lemma}
\label{lem_T_AB_Swap}
Let $A$ and $B$ be $\nbyn$ matrices.
Then for all $M$ and $N$ in $\N{A}{B}$, we have
\begin{equation*}
\T{A}{B}(MN) = \T{A}{B}(M) \schur \T{A}{B}(N).
\end{equation*}
\end{lemma}
{\sl Proof.}
If $M$ and $N$ lie in $\N{A}{B}$ then
\begin{eqnarray*}
(M N) \ \Evector{A}{i}{B}{j} &=& 
\T{A}{B}(M)_{i,j} \T{A}{B}(N)_{i,j}  \ \Evector{A}{i}{B}{j} \\
&=& (\T{A}{B}(M) \schur  \T{A}{B}(N))_{i,j} \ \Evector{A}{i}{B}{j}.
\end{eqnarray*} 
So the result follows.
\eop

\begin{lemma}
\label{lem_dualN_AB}
Let $A$ and $B$ be $\nbyn$ matrices.  
Then $\dualN{A}{B}$ is a vector space in $\mat{n}{\com}$
which is closed under the Schur product, and it contains $J$. 
\end{lemma}  
{\sl Proof.}
It follows from Lemma~\ref{lem_T_AB_Swap} that $\dualN{A}{B}$
is closed under the Schur product.
Moreover, since $I \ \Evector{A}{i}{B}{j} = 1 \ \Evector{A}{i}{B}{j}$, we have
$\T{A}{B}(I)= J$.
\eop

The following standard result (Theorem~2.6.1~\cite{BCN})
says that there is a basis of Schur idempotents for $\dualN{A}{B}$.

\begin{lemma}
\label{lem_SchurBasis}
Let ${\cal N} \subseteq \mat{n}{\com}$ be a vector space that is closed 
under the Schur product. Assume further that $J\in {\cal N}$.
Then ${\cal N}$ has a basis of Schur
idempotents.
\end{lemma}
{\sl Proof.}
Choose $M$ to be a matrix in ${\cal N}$ with the maximal number of distinct 
non-zero entries.  
We can write $M=\sum_{i=1}^d \alpha_i A_i$,
where $\seq{\alpha}{1}{d}$ are distinct and $\seq{A}{1}{d}$ are $01$-matrices
that sum to $J$ and satisfy $A_i \schur A_j = \delta_{ij} A_i$.  
Then the $A_i$'s are linearly independent.
For any $\range{i}{1}{d}$, the following product
\begin{equation*}
\left(\prod_{j=1,j\neq i}^d \frac{1}{\alpha_i - \alpha_j}\right) 
(M-\alpha_1 J) \schur \ldots (M-\alpha_{i-1} J) \schur
(M-\alpha_{i+1} J) \schur \ldots (M-\alpha_d J)
\end{equation*}
equals $A_i$ and 
belongs to ${\cal N}$.
So the dimension of ${\cal N}$ is no less than $d$.

Suppose $N \in {\cal N}$ is not a linear combination of 
$\seq{A}{1}{d}$.  Let $N = \sum_{j=1}^{k} \beta_j B_j$ where
$\seq{\beta}{1}{k}$ are distinct and $\seq{B}{1}{k}$ are $01$-matrices
that are mutually orthogonal with respect to the Schur product.

Since $\{\seq{B}{1}{k}\} \neq \{\seq{A}{1}{d}\}$ and ${\cal N}$ is
closed under the Schur product, 
\begin{equation*}
\{A_i \schur B_j : \range{i}{1}{d};\quad \range{j}{1}{k}\}
\end{equation*}
is a set of more than $d$ Schur idempotents in ${\cal N}$.
As a result, any linear combination of the matrices in this set with
distinct non-zero coefficients contradicts the choice of $M$.
We conclude that $\{ \seq{A}{1}{d}\}$ spans ${\cal N}$.   
\eop

For example, when $A=I$ and $B=J$, 
the set of eigenvectors is $\{\seq{e}{1}{n}\}$.
Therefore $\N{A}{B}$ equals 
the set of $\nbyn$ diagonal matrices.
Since for $\range{i,j}{1}{n}$,
\begin{equation*}
\Evector{I}{i}{J}{j} = e_i, 
\end{equation*}
we get $\T{I}{J}(D) = DJ$ for any diagonal matrix $D$.
The basis of Schur idempotents of $\dualN{I}{J}$ equals
$\{e_i \mathbf{1}^T : \range{i}{1}{n}\}$.
In this case, both algebras have dimension $n$.

We can say more about these algebras if we assume $A$ is invertible and
$B$ is Schur invertible.  More importantly, we will see in
Chapter~\ref{JP} that these conditions hold when $(A,B)$ is a Jones
pair.
\begin{theorem}
\label{thm_N_AB_Dim}
Let $A$ be an invertible matrix and let $B$ be a Schur invertible
matrix.
Then the duality map $\T{A}{B}$ is an isomorphism from 
the algebra $\N{A}{B}$ (with respect to matrix multiplication) to
$\dualN{A}{B}$ (with respect to the Schur product).
Moreover, $\N{A}{B}$ is a commutative algebra.
\end{theorem}
{\sl Proof.}
For each $r$, the set of vectors
\begin{equation*}
\{ \Evector{A}{i}{B}{r}: \text{ for } \range{i}{1}{n}\}
\end{equation*}
is linearly independent and hence it is a basis for $\com^n$.
As a result, 
for every $M$ in $\N{A}{B}$,  
each column of $\T{A}{B}(M)$ contains all $n$ eigenvalues of $M$.
Now we conclude that $\T{A}{B}(M) = \T{A}{B}(M')$ if and only
if $M=M'$.  So the map $\T{A}{B}$ is a bijection.
By Lemma~\ref{lem_T_AB_Swap}, if $M, N \in \N{A}{B}$, then
\begin{eqnarray*}
\T{A}{B}(MN)&=& \T{A}{B}(M) \schur \T{A}{B}(N)\\
&=& \T{A}{B}(N) \schur \T{A}{B}(M)\\
&=& \T{A}{B}(NM).
\end{eqnarray*}
It follows that $MN=NM$, for all $M, N \in \N{A}{B}$.
\eop

\section{The Duality Map $\T{A}{B}$}
\label{N_AB_XDX}

We now define two types of \SL{endomorphisms of $\mat{n}{\com}$}.
They are used by Jones \cite{VJ_Knot} to construct
braid group representations from spin models.
More importantly for us, 
Jones pairs are defined using these endomorphisms, so they
show up everywhere in this thesis.

Let $C$ be an $\nbyn$ matrix.  
We define two endomorphisms of $\mat{n}{\com}$: 
\SL{$X_C$}\ and \SL{$\D_C$}, as
follows
\begin{equation*}
X_C(M)=CM, \quad \D_C=C \schur M.
\end{equation*}
Now $X_C$ is invertible if and only if $C$ is invertible. 
Similarly, $\D_C$ is invertible if and only if $C$ is Schur invertible.
If $A$ is invertible and $B$ is Schur invertible, we have
\begin{equation*}
X_A^{-1} = X_{\inv{A}}, \quad 
\D_{B}^{-1} = \D_{\Sinv{B}}.
\end{equation*}
It is worth noting that $\D_B$ and $\D_{C}$ commute for all $\nbyn$
matrices $B$ and $C$.  Moreover, if $D$ is a diagonal matrix, then
$X_D = \D_{DJ}$. Thus $X_D$ and $\D_C$ commute.

Let $M$ and $N$ be two $\nbyn$ matrices. 
Then $\tr(M^T N) = \elsm(M \schur N)$ is a non-degenerate
bilinear form on $\mat{n}{\com}$.  
If $Y$ is an endomorphism of $\mat{n}{\com}$, 
we use $Y^T$ to denote the adjoint of $Y$ relative to this bilinear form.
We call it the \SL{transpose of $Y$}.
It is straightforward to verify that
\begin{equation*}
X_C^T=X_{C^T}, \quad \D_C^T = \D_C.
\end{equation*}

In Section~\ref{N_AB}, $\T{A}{B}(R)$ is 
defined as a store of eigenvalues of $R \in \N{A}{B}$.
The following theorem gives an equivalent definition of the duality map
using the endomorphisms $X_A$ and $\D_B$.
This definition helps us identify when a matrix belongs to $\N{A}{B}$
and it is used repeatedly in the rest of the thesis.
\begin{theorem}
\label{thm_T_AB_XDX}
Let $A, B \in \mat{n}{\com}$.
Then $R \in \N{A}{B}$ and $S = \T{A}{B}(R)$ if and only if
\begin{equation}
\label{eqn_T_AB_XDX}
\XDX{R}{B}{A}=\DXD{B}{A}{S}.
\end{equation}
\end{theorem}
{\sl Proof.}
Let $\Eij=e_i e_j^T$ be the $\nbyn$ matrix with one in its $ij$-entry
and zero elsewhere.
The set $\{\Eij: \range{i,j}{1}{n}\}$ is a basis for $\mat{n}{\com}$.
So (\ref{eqn_T_AB_XDX}) holds if and only if 
\begin{equation*}
\XDX{R}{B}{A}(\Eij) = \DXD{B}{A}{S}(\Eij), 
\end{equation*}
for all $\range{i,j}{1}{n}$.
Now the left-hand side equals 
\begin{eqnarray*}
\XDX{R}{B}{A} (e_i e_j^T) &=& R\ ( B \schur A(e_i e_j^T))\\
&=&R\ ( \Evector{B}{j}{A}{i}) e_j^T,
\end{eqnarray*}
and the right-hand side equals 
\begin{eqnarray*}
\DXD{B}{A}{S}(e_i e_j^T) &=& B \schur (A (S \schur e_i e_j^T))\\
&=& S_{i,j}\ ( \Evector{B}{j}{A}{i}) e_j^T.
\end{eqnarray*}
Both sides are equal if and only if
\begin{equation*}
R \ (\Evector{B}{j}{A}{i}) = S_{i,j} \ (\Evector{B}{j}{A}{i}).
\end{equation*}
We conclude that the relation (\ref{eqn_T_AB_XDX}) holds if and only if 
$S=\T{A}{B}(R)$.   
\eop

\section{The Exchange Lemma}
\label{Exchange}
The Exchange Lemma results from the trivial idea that
the set 
\begin{equation*}
\{ \Evector{A}{i}{B}{j} : \range{i,j}{1}{n}\}
\end{equation*}
does not change by swapping $A$ and $B$.  
Therefore $\N{A}{B} = \N{B}{A}$.
If $S=\T{A}{B}(R)$ then 
\begin{equation*}
R \ \Evector{B}{i}{A}{j} = S_{j,i} \ \Evector{B}{i}{A}{j},
\end{equation*}
which is the same as writing
$\T{B}{A}(R)=S^T$ and $\dualN{B}{A}={\dualN{A}{B}}^T$.
As a result, we have 
\begin{equation*}
\XDX{R}{B}{A}=\DXD{B}{A}{S}
\end{equation*}
if and only if 
\begin{equation*}
\XDX{R}{A}{B}=\DXD{A}{B}{S^T}.
\end{equation*}

The lemma below, called the \SL{Exchange Lemma}, gives a more general 
form of the above equivalence.
It was first discovered by Munemasa during our investigation of Jones pairs.
This result is of much greater importance than it may first appear.
It is applied at numerous places in this thesis.

\begin{lemma}
\label{lem_Exchange}
If $A,B,C,Q,R,S \in \mat{n}{\com}$ then
\begin{equation*}
\XDX{A}{B}{C}=\DXD{Q}{R}{S}
\end{equation*}
if and only if 
\begin{equation*}
\XDX{A}{C}{B}=\DXD{R}{Q}{S^T}
\end{equation*}
\end{lemma}
{\sl Proof.}
Pick any $i,j \in \{1,\ldots,n\}$. 
Applying both sides of the first relation to $\Eij$ gives
\begin{equation*}
A\ (\Evector{B}{j}{C}{i})e_j^T = S_{i,j} \ (\Evector{Q}{j}{R}{i})e_j^T.
\end{equation*}
By multiplying each side by $e_j e_i^T$, we get
\begin{equation*}
A \ (\Evector{C}{i}{B}{j})e_i^T=(S^T)_{j,i} \ (\Evector{R}{i}{Q}{j})e_i^T
\end{equation*}
which is the same as
\begin{equation*}
\XDX{A}{C}{B}(E_{ji}) = \DXD{R}{Q}{S^T}(E_{ji}).
\end{equation*}
So the result follows.   
\eop

Using the Exchange Lemma, we obtain several equivalent forms of
(\ref{eqn_T_AB_XDX}).
\begin{corollary}
\label{cor_T_AB_XDX}
If $A$ is invertible and $B$ is Schur invertible, then
the following are equivalent:
\begin{enumerate}
\item
\label{eqn_T_AB_XDX_a}
$R\in \N{A}{B}$ and $S=\T{A}{B}(R)$
\item
\label{eqn_T_AB_XDX_b}
$\XDX{B^T}{A}{R^T}=\DXD{S^T}{B^T}{A}$
\item
\label{eqn_T_AB_XDX_c}
$\DXD{R}{\Sinv{B}}{B^T}=\XDX{A}{\inv{A}}{S}$
\item
\label{eqn_T_AB_XDX_d}
$\DXD{B^T}{\SinvT{B}}{R}=\XDX{S^T}{\inv{A}}{A^T}$
\end{enumerate}
\end{corollary}
{\sl Proof.}
By Theorem~\ref{thm_T_AB_XDX} and
the Exchange Lemma to (\ref{eqn_T_AB_XDX}), Condition~(\ref{eqn_T_AB_XDX_a})
holds if and only if 
\begin{equation*}
\XDX{R}{A}{B}=\DXD{A}{B}{S^T}.
\end{equation*}
Taking the transpose of both sides yields (\ref{eqn_T_AB_XDX_b}).
So (\ref{eqn_T_AB_XDX_a}) and (\ref{eqn_T_AB_XDX_b}) are equivalent.
Moreover we can write (\ref{eqn_T_AB_XDX}) as
\begin{equation*}
\DXD{\Sinv{B}}{R}{B}=\XDX{A}{S}{\inv{A}}.
\end{equation*}
Applying the Exchange Lemma gives (\ref{eqn_T_AB_XDX_c}).
Lastly, we obtain (\ref{eqn_T_AB_XDX_d}) by taking the 
transpose of each side of (\ref{eqn_T_AB_XDX_c}).   
As a result, Conditions~(\ref{eqn_T_AB_XDX_a}), (\ref{eqn_T_AB_XDX_c}) 
and (\ref{eqn_T_AB_XDX_d}) are equivalent.

\eop

Assuming all diagonal entries of $A$ are non-zero, we use
Corollary~\ref{cor_T_AB_XDX}~(\ref{eqn_T_AB_XDX_b}) to get an explicit
formula of $\T{A}{B}$.  
All equivalent forms of Equation~(\ref{eqn_T_AB_XDX}) give formulae of
$\T{A}{B}$.  
We choose this particular one because we will use it in the proof of the 
modular invariance equation in Section~\ref{Mod_Inv}.
\begin{lemma}
\label{lem_T_AB_Formula}
Let $A$ and $B$ be $\nbyn$ matrices, and let $R \in \N{A}{B}$.
If $A$ and $A \schur I$ are invertible and $B$ is Schur invertible, then
\begin{equation*}
\T{A}{B}(R) = \Sinv{B} \schur (\inv{(A\schur I)}(A^T \schur R) B).
\end{equation*}
\end{lemma}
{\sl Proof.}
Let $S=\T{A}{B}(R)$.
Applying Corollary~\ref{cor_T_AB_XDX}~(\ref{eqn_T_AB_XDX_b}) to $I$, we get
\begin{equation*}
B^T\ (A \schur (R^T I)) = S^T \schur (B^T\ (A\schur I)).
\end{equation*}
By taking the transpose of each side, we have
\begin{equation*}
S \schur ((A^T \schur I)B) = (A^T \schur R)B.
\end{equation*}
Since $A^T\schur I = A \schur I$ is a diagonal matrix, 
the left-hand side equals $(A\schur I) (S \schur B)$.
Hence we get
\begin{equation*}
S = \Sinv{B} \schur \left(\inv{(A\schur I)}(A^T \schur R) B \right).
\end{equation*}  
\eop

\section{Transformations of $(A,B)$}
\label{N_AB_Transform}

Two matrices $A$ and $C$ are \SL{monomially equivalent}\ if 
$C=MAN$, where $M$ and $N$ are products of permutation matrices and
diagonal matrices.
When $A$ is a type-II matrix, Jaeger, Matsumoto and Nomura examined the 
the relation between the Nomura algebras of the pairs
$(A,\Sinv{A})$ and $(C,\Sinv{C})$, 
(see Proposition~2 and 3 in \cite{JMN}).
In this section, we extend their result to $\N{A}{B}$ 
for any pair of $\nbyn$ matrices $(A,B)$.

As we will see in Sections~\ref{Odd_Gauge} and \ref{Even_Gauge},
if $(A_1,B_1)$ and $(A_2,B_2)$ are guage equivalent invertible Jones Pairs,
then $A_1$ and $A_2$ are monomially equivalent, and so are $B_1$ and $B_2$.
So the following lemmas tell us how the Nomura algebras of gauge-equivalent
invertible Jones pairs relate to each other.

\begin{lemma}
\label{lem_Scaling}
If $D,E$ and $F$ are invertible diagonal $\nbyn$ matrices, then
\begin{equation*}
\N{DAE}{\inv{D}B F}=\N{A}{B},
\quad \text{and} \quad
\dualN{DAE}{\inv{D}B F}=
\dualN{A}{B}. 
\end{equation*}
\end{lemma}
{\sl Proof.}
Since $D\ (B\schur C) = (DB \schur C)$, we have
\begin{equation*}
(DAEe_i) \schur (\inv{D}BFe_j) = 
(E)_{i,i} (F)_{j,j} \ \Evector{A}{i}{B}{j}
\end{equation*}
which implies the lemma.   
\eop

\begin{lemma}
\label{lem_Perm}
If $P,Q$ and $R$ are $\nbyn$ permutation matrices, then
\begin{equation*}
\N{PAQ}{PB R}= 
P\ \N{A}{B}\inv{P}, 
\quad \text{and} \quad
\dualN{PAQ}{PB R}=
\inv{Q}\ \dualN{A}{B}R.
\end{equation*}
\end{lemma}
{\sl Proof.}
Let $M \in \N{A}{B}$. Then for all $\range{i,j}{1}{n}$,
\begin{eqnarray*}
(P M \inv{P}) \ (\Evector{PAQ}{i}{PBR}{j}) &=&
PM \inv{P} P \ (\Evector{AQ}{i}{BR}{j})\\
&=&PM\ (\Evector{AQ}{i}{BR}{j}).
\end{eqnarray*}
If $Qe_i=e_{i'}$ and $Re_j=e_{j'}$, we get
\begin{eqnarray*}
(P M \inv{P}) \ (\Evector{PAQ}{i}{PBR}{j}) 
&=& PM  \ (\Evector{A}{i'}{B}{j'})\\
&=& {\T{A}{B}(M)}_{i',j'}  \ (\Evector{PAQ}{i}{PBR}{j}),
\end{eqnarray*}
So the first equality of the lemma follows.
Moreover,  
${\T{A}{B}(M)}_{i',j'}= (Q^T \ \T{A}{B}(M) R)_{i,j}, $ 
and the second equality holds.
\eop

\section{Tensor Products}
\label{N_AB_Tensor}

We use $\otimes$ to denote both the Kronecker product
of two matrices and the tensor product of two algebras.
The next lemma generalizes Proposition~7 in \cite{JMN}, due to
Jaeger, Matsumoto and Nomura.
In particular, we will use this lemma in Section~\ref{1JP} to show
that the tensor product of two Jones pairs is also a Jones pair.

\begin{lemma}
\label{lem_N_AB_Tensor}
Let $A_1, B_1 \in \mat{m}{\com}$ and 
$A_2, B_2 \in \mat{n}{\com}$.
If $A_1, A_2$ are invertible, and $B_1, B_2$ are
Schur invertible, then
\begin{equation*}
\N{A_1 \otimes A_2}{B_1 \otimes B_2} = 
\N{A_1}{B_1} \otimes \N{A_2}{B_2}, 
\end{equation*}
and
\begin{equation*}
\dualN{A_1 \otimes A_2}{B_1 \otimes B_2} = 
\dualN{A_1}{B_1} \otimes \dualN{A_2}{B_2}.
\end{equation*}
\end{lemma}
{\sl Proof.}
Since $\{e_i \otimes e_h : \range{i}{1}{m} \text{ and }
\range{h}{1}{n}\}$ form a basis for $\com^{mn}$,
we can write the eigenvectors for the matrices in
$\N{A_1 \otimes A_2}{B_1 \otimes B_2}$ as 
\begin{eqnarray*}
&&
\left((A_1 \otimes A_2) (e_{i} \otimes e_{h})\right) \schur
\left((B_1 \otimes B_2) (e_{j} \otimes e_{k})\right) \\
&=&
(A_1 e_i \otimes A_2 e_h) \schur (B_1 e_j \otimes B_2 e_k)\\
&=&
(\Evector{A_1}{i}{B_1}{j}) \otimes \ (\Evector{A_2}{h}{B_2}{k}),
\end{eqnarray*}
for all $\range{i,j}{1}{m}$ and $\range{h,k}{1}{n}$.
This implies if $M \in \N{A_1}{B_1}$ and $N \in \N{A_2}{B_2}$,
then $M\otimes N \in \N{A_1 \otimes A_2}{B_1 \otimes B_2}$.
Consequently, 
\begin{equation*}
\N{A_1}{B_1} \otimes \N{A_2}{B_2} \subseteq
\N{A_1 \otimes A_2}{B_1 \otimes B_2}. 
\end{equation*}

Now $A_1 \otimes A_2$ is invertible and $B_1 \otimes B_2$ is Schur invertible.
By Theorem~\ref{thm_N_AB_Dim}, we have for $i=1,2$,
$\dim{\N{A_i}{B_i}} = \dim{\dualN{A_i}{B_i}}$
and
$\dim{\N{A_1 \otimes A_2}{B_1 \otimes B_2}}=
\dim{\dualN{A_1 \otimes A_2}{B_1 \otimes B_2}}$.

It is sufficient to show that
${\dualN{A_1 \otimes A_2}{B_1 \otimes B_2}}$ and 
${\dualN{A_1}{B_1} \otimes \dualN{A_2}{B_2}}$ have the same dimension.
If the Hermitian product of the two eigenvectors
\begin{equation*}
\left((A_1 \otimes A_2) (e_{i_1} \otimes e_{h_1})\right) \schur
\left((B_1 \otimes B_2) (e_{j_1} \otimes e_{k_1})\right),
\end{equation*}
and
\begin{equation*}
\left((A_1 \otimes A_2) (e_{i_2} \otimes e_{h_2})\right) \schur
\left((B_1 \otimes B_2) (e_{j_2} \otimes e_{k_2})\right)
\end{equation*}
is non-zero,
then they belong to the same eigenspace of $M$, 
for all $M$ in $\N{A_1 \otimes A_2}{B_1 \otimes B_2}$.
Thus
\begin{equation*}
{\T{A_1\otimes A_2}{B_1 \otimes B_2}(M)}_{((i_1,h_1),(j_1,k_1))} =
{\T{A_1\otimes A_2}{B_1 \otimes B_2}(M)}_{((i_2,h_2),(j_2,k_2))}.
\end{equation*}

But the Hermitian product of those two eigenvectors equals
\begin{equation*}
\ip{\Evector{A_1}{i_1}{B_1}{j_1}}{\Evector{A_1}{i_2}{B_1}{j_2}}
\ip{\Evector{A_2}{h_1}{B_2}{k_1}}{\Evector{A_2}{h_2}{B_2}{k_2}},
\end{equation*}
and so the $((i_1,h_1),(j_1,k_1))$ and $((i_2,h_2),(j_2,k_2))$ entries
of a Schur idempotent of $\dualN{A_1 \otimes A_2}{B_1 \otimes B_2}$ 
equal to one if and only if  
there exist a Schur idempotent $F$ in $\dualN{A_1}{B_1}$ 
and a Schur idempotent $G$ in $\dualN{A_2}{B_2}$ such that
\begin{equation*}
F_{i_1,j_1}=F_{i_2,j_2}=1,\quad \text{and }
G_{h_1,k_1}=G_{h_2,k_2}=1.
\end{equation*}

In other words, the set of matrices $F_r \otimes G_s$,
for all Schur idempotents $F_r$ of $\dualN{A_1}{B_1}$ and
all Schur idempotents $G_s$ of $\dualN{A_2}{B_2}$,
forms a basis of Schur idempotents for 
$\dualN{A_1 \otimes A_2}{B_1 \otimes B_2}$.  So the result follows.  
\eop


\chapter{Type-II Matrices and Nomura Algebras}
\label{IINom}

Jaeger, Matsumoto and Nomura showed in \cite{JMN} that
the type-II condition of a spin model is sufficient for the existence of
a Bose-Mesner algebra containing the spin model.  
Further, they proved that a type-II matrix belongs to its Nomura algebra 
if and only if it is a spin model up to scalar multiplication.
If a matrix satisfies the type-II condition of a spin model, we call it
a type-II matrix.

Since the type-II condition plays a crucial role in the connection of
spin models to Bose-Mesner algebras, we feel obliged to give a
detailed treatment to type-II matrices and their Nomura algebras,
see Sections~\ref{Type_II}, \ref{N_A} and \ref{TypeII_Families}.
Most results in these sections are originally due to Jaeger,
Matsumoto and Nomura in \cite{JMN}.  
However, we present new proofs, using the tools developed in the
previous chapter, given by Godsil, Munemasa and the author in \cite{CGM}.
Sections~\ref{Scheme_BM_Alg} to \ref{Dual_Scheme} 
consist of the standard theory of 
Bose-Mesner algebras and association schemes.
Section~\ref{N_II} examines the Nomura algebra of a pair of type-II matrices
$A$ and $B$.  This section gives the foundation for the theory of
invertible Jones pairs in Chapters~\ref{JP} and \ref{InvJP}.
The results in Section~\ref{N_II} also appear in \cite{CGM}.

\section{Type-II Matrices}
\label{Type_II}

If an $\nbyn$ Schur-invertible matrix $A$ satisfies
\begin{equation}
\label{eqn_TypeII}
A \SinvT{A} = nI,
\end{equation}
we call it a \SL{type-II matrix}.  
In other words, $A$ is a type-II matrix if and only if 
\begin{equation*}
\inv{A}=n^{-1} \SinvT{A},
\end{equation*}
or equivalently,
for $\range{i,j}{1}{n}$,
\begin{equation*}
\sum_{k=1}^n \frac{A_{i,k}}{A_{j,k}} = \delta_{i,j} n
\quad \text{and} \quad
\sum_{k=1}^n \frac{A_{k,i}}{A_{k,j}} = \delta_{i,j} n.
\end{equation*}
Note that Equation~(\ref{eqn_TypeII}) is equivalent to the definition of 
type-II matrix in Section~\ref{Intro_Intro}.
Also note that $A$ is a type-II matrix if and only if $A^T$ is also type II.

Suppose $A$ is a type-II matrix.
If $D_1$ and $D_2$ are invertible diagonal matrices, then
\begin{equation*}
(D_1AD_2)\SinvT{(D_1AD_2)} = (D_1AD_2)(\inv{D_2} \SinvT{A} \inv{D_1}) = nI.
\end{equation*}
Therefore $D_1AD_2$ is also a type-II matrix.
Similarly, if $P_1$ and $P_2$ are permutation matrices, then
\begin{equation*}
(P_1AP_2)\SinvT{(P_1AP_2)} = (P_1AP_2)({P_2}^T \SinvT{A} P_1^T) = nI,
\end{equation*}
and $P_1AP_2$ is a type-II matrix.
Note that $A$, $P_1AP_2$ and $D_1AD_2$ are monomially equivalent.
The relation between the Nomura algebras of two monomially equivalent
matrices is discussed in Section~\ref{N_AB_Transform}.

\begin{example}
\label{ex_TypeII_2to5}
We list all type-II matrices of orders two to five, 
up to monomially equivalence.  
For details, please see \cite{JMN} and \cite{N_II5}.
\begin{enumerate}
\item
$n=2$:
$ \begin{pmatrix}
1 & 1 \\ 1 & -1
\end{pmatrix} $
\item
$n=3$:
$ \begin{pmatrix}
1 & 1 & \omega \\ \omega & 1 & 1 \\ 1 & \omega & 1
\end{pmatrix} $, 
where $\omega$ is a cube root of unity.
\item
$n=4$:
$ \begin{pmatrix}
1&1&1&1\\1&1&-1&-1\\1&-1&\lambda&-\lambda\\1&-1&-\lambda&\lambda
\end{pmatrix} $, 
for any non-zero complex number $\lambda$.
\item 
$n=5$: 
for $\eta$ a fifth-root of unity,
\begin{equation*}
\begin{pmatrix}
1&1&1&1&1\\1&\eta&\eta^2&\eta^3&\eta^4\\1&\eta^2&\eta^4&\eta&\eta^3\\
1&\eta^3&\eta&\eta^4&\eta^2\\1&\eta^4&\eta^3&\eta^2&\eta
\end{pmatrix}; 
\end{equation*}
and
\begin{equation*}
\left( \frac{-5-\sqrt{5}}{2}\right) I+J; \quad 
\left( \frac{-5+\sqrt{5}}{2}\right) I+J.
\end{equation*}
\end{enumerate}
\end{example}

The Kronecker product of two type-II matrices
is also type II.  So there exist infinitely many type-II matrices.
Some significant examples are spin models and four-weight spin models.
Furthermore, if a matrix is unitary and all its entries are roots of unity, 
then it is type II.
These matrices are objects of interest in the theory of Von Neumann algebras
\cite{JS}.

\section{Nomura Algebras of a Type-II Matrix}
\label{N_A}

Godsil \cite{G_II} observed the following condition on $\N{A}{\Sinv{A}}$
that is equivalent to $A$ being type II.
\begin{lemma}
\label{lem_J_NA}
Let $A$ be both Schur invertible and invertible.
Then $A$ is type II if and only if $J \in \N{A}{\Sinv{A}}$.
Moreover, $\T{A}{\Sinv{A}}(J)=nI$.
\end{lemma}
{\sl Proof.}
Assume $J \in \N{A}{\Sinv{A}}$.
For each $\range{j}{1}{n}$, the set 
\begin{equation*}
\{\Evector{A}{1}{\Sinv{A}}{j}, \ldots, \ \Evector{A}{n}{\Sinv{A}}{j}\}
\end{equation*}
is a basis for $\com^n$ consisting eigenvectors of $J$.
Now $\Evector{A}{j}{\Sinv{A}}{j}$ equals $\mathbf{1}$,
the vector of all ones, which is the only eigenvector of
$J$ with non-zero eigenvalue. So we conclude that
\begin{equation*}
J \ \Evector{A}{i}{\Sinv{A}}{j}=\delta_{i,j}n \ \Evector{A}{i}{\Sinv{A}}{j}.
\end{equation*}
That is, for all $\range{i,j}{1}{n}$,
\begin{equation*}
\sum_{k=1}^n \frac{A_{k,i}}{A_{k,j}} = \delta_{i,j} n, \quad
\end{equation*}
and $A$ is type II.

The converse is straightforward.  
\eop

Following existing conventions, we use \SL{$\Nom{A}$}, \SL{$\dualNom{A}$}\ 
and \SL{$\Theta_{A}$}\ to stand for $\N{A}{\Sinv{A}}$, $\dualN{A}{\Sinv{A}}$ and 
$\T{A}{\Sinv{A}}$ respectively.

\begin{theorem}
\label{thm_T_A}
Let $A$ be a type-II matrix.
If $R \in \Nom{A}$ then $\Theta_A(R) \in \Nom{A^T}$
and $\Theta_{A^T}(\Theta_{A} (R)) = nR^T$.
\end{theorem}
{\sl Proof.}
Letting $B=\Sinv{A}$ in Corollary~\ref{cor_T_AB_XDX}~(\ref{eqn_T_AB_XDX_c}), 
we have $R \in \Nom{A}$ and $S=\Theta_A(R)$ if and only if 
\begin{equation*}
\DXD{R}{A}{\SinvT{A}} = \XDX{A}{\inv{A}}{S},
\end{equation*}
which is the same as
\begin{equation*}
\DXD{\Sinv{(\inv{A})}}{\inv{A}}{R} = \XDX{S}{A^T}{\inv{A}}.
\end{equation*}
Since $\inv{A} = n^{-1} \SinvT{A}$,
we can rewrite the above as
\begin{equation*}
n \DXD{A^T}{\SinvT{A}}{R} = \XDX{S}{A^T}{\SinvT{A}}.
\end{equation*}
Applying the Exchange Lemma, we get
\begin{equation*}
\DXD{\SinvT{A}}{A^T}{nR^T} = \XDX{S}{\SinvT{A}}{A^T}.
\end{equation*}
By Theorem~\ref{thm_T_AB_XDX}, we conclude that $S$ belongs to $\Nom{A^T}$
and $\Theta_{A^T}(S) = n R^T$.  
\eop

\begin{corollary}
\label{cor_T_A}
If $A$ is a type-II matrix, then $\dualNom{A} = \Nom{A^T}$
and $\dualNom{A^T} = \Nom{A}$.
Moreover, $R \in \Nom{A}$ if and only if $R^T \in \Nom{A}$.
\end{corollary}
{\sl Proof.}
Applying Theorem~\ref{thm_T_A} to the type-II matrix $A$,
we get $\dualNom{A} \subseteq \Nom{A^T}$.
Since $A^T$ is also type II, so $\dualNom{A^T} \subseteq \Nom{A}$.
By Theorem~\ref{thm_N_AB_Dim},
\begin{equation*}
\dim{\Nom{A}} = \dim{\dualNom{A}}, \quad \text{and }
\dim{\Nom{A^T}} = \dim{\dualNom{A^T}}.
\end{equation*}
As a result, we get $\dualNom{A}=\Nom{A^T}$ and $\dualNom{A^T}=\Nom{A}$.

The second part of  the corollary follows from Theorem~\ref{thm_T_A}
and the fact that $\dualNom{A^T}=\Nom{A}$.
\eop

We deduce from Theorem~\ref{thm_N_AB_Dim} and this corollary 
that $\Nom{A}$ is a commutative algebra 
(with respect to matrix multiplication) containing
$I$ and $J$ which is 
closed under the Schur product and the transpose.
An algebra with all these properties is called a \SL{Bose-Mesner algebra}.  
We introduce the theory of Bose-Mesner algebras in Section~\ref{Scheme_BM_Alg}.

\begin{theorem}
\label{thm_N_A_BM}
If $A$ is a type-II matrix, then 
both $\Nom{A}$ and $\dualNom{A}$ are Bose-Mesner algebras.
\eop
\end{theorem}

\begin{lemma}
\label{lem_T_A_Swap}
If $A$ is a type-II matrix and if $R_1, R_2 \in \Nom{A}$ then
\begin{equation*}
\Theta_A(R_1 R_2) = \Theta_A(R_1) \schur \Theta_A(R_2)
\end{equation*}
and
\begin{equation*}
\Theta_A(R_1 \schur R_2) = n^{-1} \Theta_A(R_1) \Theta_A(R_2).
\end{equation*}
\end{lemma}
{\sl Proof.}
The first equation follows directly from Lemma~\ref{lem_T_AB_Swap}.
From the equality $\Nom{A}=\dualNom{A^T}$,
there exist $S_1$ and $S_2$ in $\Nom{A^T}$ such that 
$\Theta_{A^T}(S_i)=R_i$, for $i=1,2$.
Then by Theorem~\ref{thm_T_A}
\begin{eqnarray*}
\Theta_A(R_1 \schur R_2) &=& \Theta_A(\Theta_{A^T}(S_1) \schur \Theta_{A^T}(S_2))\\
&=& \Theta_A(\Theta_{A^T}(S_1 S_2))\\
&=& n S_2^T S_1^T.
\end{eqnarray*}
The second equation follows from Theorem~\ref{thm_T_A} and the commutativity of $\Nom{A^T}$.
\eop

We have shown that the duality map $\Theta_A$ interchanges the Schur product and 
the matrix multiplication.
Moreover, the following lemma shows that 
$\Theta_A$ and the transpose map commute.
\begin{lemma}
\label{lem_T_A_Transpose}
Let $A$ be a type-II matrix.
If $R \in \Nom{A}$, then $\Theta_A(R^T) = \Theta_A(R)^T$.
\end{lemma}
{\sl Proof.}
Let $S = \Theta_A(R)$.  
By setting $B=\Sinv{A}$ in 
Corollary~\ref{cor_T_AB_XDX}~(\ref{eqn_T_AB_XDX_b}), we get 
\begin{equation*}
\XDX{\SinvT{A}}{A}{R^T} = \DXD{S^T}{\SinvT{A}}{A},
\end{equation*}
which can be rewritten as 
\begin{equation*}
\XDX{R^T}{\Sinv{A}}{\inv{(\SinvT{A})}} =
\DXD{\Sinv{A}}{\inv{(\SinvT{A})}}{S^T}.
\end{equation*}
Since $\SinvT{A} = n \inv{A}$, we have
\begin{equation*}
\XDX{R^T}{\Sinv{A}}{A} = \DXD{\Sinv{A}}{A}{S^T}.
\end{equation*}
So by Theorem~\ref{thm_T_AB_XDX},
we get $\Theta_A(R^T) = \Theta_A(R)^T$.  
\eop

\section{Association Schemes and Bose-Mesner Algebras}
\label{Scheme_BM_Alg}

A \SL{Bose-Mesner algebra}\ is a finite dimensional vector space
of $\nbyn$ matrices that is closed under the transpose,
the Schur product and the matrix multiplication.
It is commutative with respect to the matrix multiplication
and it contains $I$ and $J$.
The Nomura algebras of type-II matrices are Bose-Mesner algebras.
As we see in this section,
Bose-Mesner algebras are equivalent to association schemes.
Association schemes can be viewed as partitions of the complete
graph on $n$ vertices into directed graphs that satisfy
some regular conditions.
Sections~\ref{ex_Schemes} to \ref{Dual_Scheme} serve as an
introduction to the theory of association schemes.
For further information on association schemes, please refer to \cite{BCN}.

In this thesis, we choose to give the definition of association schemes
in terms of matrices.
An \SL{association scheme on $n$ elements with $d$ classes}\ 
is a set of $\nbyn$ $01$-matrices $\A=\{A_0,\ldots,A_d\}$
that satisfies the following conditions:
\begin{enumerate}
\item
$A_0=I$
\item
$\sum_{i=0}^d A_i = J$
\item
$A_i^T = A_{i'}$, for some $i' \in \{0,\ldots,d\}$
\item
There exist non-negative integers $p_{ij}^k$ such that
for all $\range{i,j}{0}{d}$,
\begin{equation*}
A_i A_j  = \sum_{k=0}^d p_{ij}^k A_k.
\end{equation*}
\item
$A_i A_j = A_j A_i$, for all $\range{i,j}{0}{d}$.
\end{enumerate}
If $A_i^T=A_i$ for $\range{i}{0}{d}$, we say that $\A$ is 
a symmetric association scheme.
Let $\BM$ be the span of $\{ \seq{A}{0}{d}\}$.
Condition~(b) says that $A_i \schur A_j = \delta_{i,j}A_i$ for all
$\range{i,j}{1}{n}$.  Hence $\BM$ is closed under the Schur product.
Conditions~(d) and (e) tell us that $\BM$ is closed under and
commutative with respect to the matrix multiplication.
Now the set $\{\seq{A}{0}{d}\}$ is a basis of Schur idempotents for
$\BM$.  With the first three conditions, we know that $\BM$ is
also closed under the transpose and it contains $I$ and $J$.
Consequently, we get a Bose-Mesner algebra of dimension $d+1$
for each association scheme with $d$ classes.

Conversely, by Lemma~\ref{lem_SchurBasis}, we know that any 
Bose-Mesner algebra of dimension $m$ has a basis of Schur idempotents
$\seq{A}{0}{m-1}$.  
It is standard result that the properties of a Bose-Mesner algebra enforce
conditions~(a) to (e) to hold for $\seq{A}{0}{m-1}$.

\section{Examples of Association Schemes}
\label{ex_Schemes}

We list examples of association schemes that often appear
in the context of spin models and four-weight spin models.

The simplest example is the \SL{trivial scheme}\ $\A=\{I,J-I\}$.
Its Bose-Mesner algebra contains the Potts model.

This family of association schemes has the most number of
classes possible.
Let $X$ be a finite Abelian group and $n=|X|$.
For each $z \in X$, define the $\nbyn$ $01$-matrix $A_z$ by
\begin{equation*}
(A_z)_{x,y} = \delta_{y-x,z}.
\end{equation*}
These $n$ permutation matrices form an association scheme
with $n-1$ classes, called the \SL{Abelian group scheme of $X$}.
Bannai, Bannai and Jaeger showed in \cite{BBJ} that 
the Bose-Mesner algebra of an Abelian group scheme always contains a spin model.

The last examples are two association schemes with four classes
constructed from a Hadamard matrix.
Let $H$ be an $\nbyn$ Hadamard matrix.  Then
\begin{eqnarray*}
&
\begin{pmatrix}
I&0&0&0\\0&I&0&0\\0&0&I&0\\0&0&0&I
\end{pmatrix},
\begin{pmatrix}
0&I&0&0\\I&0&0&0\\0&0&0&I\\0&0&I&0
\end{pmatrix},
\begin{pmatrix}
J-I&J-I&0&0\\J-I&J-I&0&0\\0&0&J-I&J-I\\0&0&J-I&J-I
\end{pmatrix}, &\\
&
\begin{pmatrix}
0&0&\frac{J+H}{2}&\frac{J-H}{2}\\0&0&\frac{J-H}{2}&\frac{J+H}{2}\\
\frac{J+H^T}{2}&\frac{J-H^T}{2}&0&0\\\frac{J-H^T}{2}&\frac{J+H^T}{2}&0&0
\end{pmatrix},
\begin{pmatrix}
0&0&\frac{J-H}{2}&\frac{J+H}{2}\\0&0&\frac{J+H}{2}&\frac{J-H}{2}\\
\frac{J-H^T}{2}&\frac{J+H^T}{2}&0&0\\\frac{J+H^T}{2}&\frac{J-H^T}{2}&0&0
\end{pmatrix} &
\end{eqnarray*}
form a symmetric association scheme on $4n$ elements with four classes.

Replacing the last two matrices above by
\begin{equation*}
\begin{pmatrix}
0&0&\frac{J+H}{2}&\frac{J-H}{2}\\0&0&\frac{J-H}{2}&\frac{J+H}{2}\\
\frac{J-H^T}{2}&\frac{J+H^T}{2}&0&0\\\frac{J+H^T}{2}&\frac{J-H^T}{2}&0&0
\end{pmatrix},
\begin{pmatrix}
0&0&\frac{J-H}{2}&\frac{J+H}{2}\\0&0&\frac{J+H}{2}&\frac{J-H}{2}\\
\frac{J+H^T}{2}&\frac{J-H^T}{2}&0&0\\\frac{J-H^T}{2}&\frac{J+H^T}{2}&0&0
\end{pmatrix},
\end{equation*}
we obtain a non-symmetric association scheme with four classes.
In \cite{JN}, Jaeger and Nomura constructed symmetric and
non-symmetric Hadamard spin models.  They are contained in the
Bose-Mesner algebras of the above schemes, respectively.

\section{Idempotents and Eigenvalues}
\label{Idem_Eigen}

Let $\BM$ be a Bose-Mesner algebra of an association scheme with $d$ classes.
Since the matrices in $\BM$ are normal and they commute with respect to 
the matrix multiplication, they are simultaneously diagonalizable.
Therefore $\BM$ has a basis $\{ \seq{E}{0}{d}\}$ such that
$E_i$ is the orthogonal projection onto the $i$-th common eigenspace 
of the matrices in $\BM$.  So we have 
\begin{equation*}
\sum_{k=0}^d E_k = I, \quad \text{and} \quad  E_iE_j=\delta_{i,j}E_i,
\end{equation*}
for all $\range{i,j}{1}{n}$.
The $E_i$'s are called the \SL{principal idempotents}\ 
of the association scheme.

By the definition of the principal idempotents, there exist complex
numbers $p_i(j)$'s such that
\begin{equation*}
A_i E_j = p_i(j) E_j,
\end{equation*}
for all $\range{i,j}{0}{d}$.
The numbers $p_i(j)$'s are the eigenvalues for the Schur idempotents.
Define $P$ to be the matrix whose $ji$-entry equals $p_i(j)$.
We call $P$ the \SL{matrix of eigenvalues}\  of the association scheme.
Now, for each $\range{i}{0}{d}$, we can write
\begin{equation*}
A_i = \sum_{j=0}^d P_{j,i} E_j.
\end{equation*}
Similarly, if we let $Q=n\inv{P}$, then for each $\range{i}{0}{d}$,
\begin{equation*}
E_i = n^{-1} \sum_{j=0}^d Q_{j,i} A_j.
\end{equation*}
Note that $E_i \schur A_j = n^{-1}Q_{j,i} A_j$.
The entries of $Q$ act like the eigenvalues of $E_i$ with respect to
the Schur product.
Hence we name $Q$ the \SL{matrix of dual eigenvalues}\ 
of the association scheme.

\section{Dualities of Association Schemes}
\label{Dual_Scheme}

By Theorem~\ref{thm_N_A_BM}, if $A$ is a type-II matrix, then
$\Nom{A}$ and $\Nom{A^T}$ are Bose-Mesner algebras.
Moreover, the map $\Theta_A:\Nom{A} \rightarrow \Nom{A^T}$
satisfies 
\begin{eqnarray*}
\Theta_A(MN) &=& \Theta_A(M) \schur \Theta_A(N), \quad \text{and} \\
\Theta_A(M \schur N) &=& n^{-1} \Theta_A(M) \Theta_A(N).
\end{eqnarray*}

In general, a \SL{duality}\  between two Bose-Mesner algebras 
$\BM_1$ and $\BM_2$ 
is an invertible linear map $\Psi : \BM_1 \rightarrow \BM_2$ 
that satisfies
\begin{equation*}
\Psi(MN)=\Psi(M) \schur \Psi(N), \quad \text{and} \quad
\Psi(M \schur N) = n^{-1} \Psi(M) \Psi(N).
\end{equation*}
We say that the $\BM_1$ and $\BM_2$ (or their corresponding association
schemes) are \SL{formally dual}\  to each other.
As we have seen in Section~\ref{N_A},
Nomura's construction provide an abundant source of formally dual pairs
of association schemes.

When $\BM_1 = \BM_2$ and $\Psi^2(M)=nM^T$, 
we say that $\BM_1$ is \SL{formally self-dual}.
In this case, $\Psi$ maps the basis of Schur idempotents of $\BM_1$ to
its basis of principal idempotents.
It is possible to order the Schur idempotents and the principal idempotents
so that the matrix of eigenvalues $P$ is the matrix of $\Psi$ with respect to
the basis of Schur idempotents.  If $T$ is the matrix of the transpose
with respect to the basis of Schur idempotents, then we have $P^2 = nT$.

The duality map and the matrix of eigenvalues of the Nomura algebra of
a spin model are crucial in the derivation of the modular invariance
equation.  We present this in Section~\ref{Mod_Inv}.

\section{Infinite Families of Type-II Matrices}
\label{TypeII_Families}

We record four families of type-II matrices in this section.
Each of these type-II matrices is monomially equivalent to some spin model.

The first example is the simplest family of type-II matrices.
Let $A=tI+(J-I)$ for some non-zero scalar $t$.
Then
\begin{eqnarray*}
A \SinvT{A} &=&
(tI + (J-I))(t^{-1} I + (J-I)) \\
&=& (-t-t^{-1}+2)I+ (t+t^{-1}-2+n)J.
\end{eqnarray*}
Therefore $A$ is type II if and only if $t+t^{-1}-2 +n=0$.
For each $n\geq 2$, the solutions to the quadratic equation 
give two type-II matrices contained in 
the Bose-Mesner algebra of the trivial scheme on $n$ elements.

Now for any distinct $i$ and $j$ not equal to $1$,
\begin{eqnarray*}
\ip{\Evector{A}{1}{\Sinv{A}}{i}}{\Evector{A}{1}{\Sinv{A}}{j}} &=& t^2+2t^{-1}+n-3\\
&=&-n(t+1),
\end{eqnarray*}
which is non-zero when $n\neq 4$.
That is, when $n\neq 4$, no two vectors from 
\begin{equation*}
(\Evector{A}{1}{\Sinv{A}}{2}),\ldots, \ (\Evector{A}{1}{\Sinv{A}}{n})
\end{equation*}
are orthogonal to each other.  So, they all lie in the same eigenspace.
Hence the matrices in $\Nom{A}$ have at most two eigenspaces, one spanned by
the $n-1$ vectors listed above and the other spanned by the vector
$\Evector{A}{1}{\Sinv{A}}{1}$.
This implies that all matrices in $\dualNom{A}$ have at most two 
distinct entries and therefore $\dualNom{A} = \spn(I,J)$.
By Corollary~\ref{cor_T_A}, we have $\dualNom{A} = \Nom{A^T}$.
But $A$ is symmetric, so we conclude that $\Nom{A}=\Nom{A^T}$ is the Bose-Mesner 
algebra of the trivial scheme.

When $n=4$, we have $A=-I+(J-I)$ and it is shown in Section~{5.3} of \cite{JMN}
that $\Nom{A}$ is the Abelian group scheme for $\mathbb{Z}_2 \times \mathbb{Z}_2$.
Finally, the \SL{Potts model}\  is defined to be $u^{-1}A$ where
$u^4=-t$, and it is easy to see that $A \in \Nom{A}$.

The second example comes from finite Abelian groups.
It is documented in Section~5.1 of \cite{JMN}.
Let X be a finite Abelian group and let $\A_X$ denotes its Abelian group scheme.
For any $x,y \in X$, we have $A_x A_y = A_{x+y}$.
Therefore the entries of the matrix of eigenvalues satisfy
$P_{z,x} P_{z,y} = P_{z,x+y}$.
Thus each row of $P$ is a character of $X$ which implies 
$P$ is a type-II matrix.
It is shown in Section~{5.1} of \cite{JMN} that $\dualNom{P} = \A_X$.
However, in general, $\Nom{P}$ may not equal $\dualNom{P}$ and 
it may not contain $P$.
Further, if $X$ is a cyclic group of order $n$, we define $W$ to have entries
$W_{x,y} = \omega^{(x-y)^2}$ for some $n$-th root of unity $\omega$.
Then $W$ is symmetric and $W \in \Nom{W} = A_X$.
This matrix is a spin model.

An $\nbyn$ $\{1,-1\}$-matrix $H$ is called a \SL{Hadamard matrix}\  if $H H^T=nI$.
Since $H^T = \SinvT{H}$, Hadamard matrices form an infinite family of
type-II matrices.  
Using easy counting argument, when $n\geq 12$ and $n \equiv 4 \mod{8}$,
$\Nom{H}$ is just the span of $\{I,J\}$ (see Section~{5.2} of \cite{JMN}).

In \cite{JN} and \cite{N_Had}, Jaeger and Nomura constructed 
two $4n \times 4n$ spin models (hence type-II matrices) 
from each $\nbyn$ Hadamard matrix $H$.
We will provide a new proof that they are spin models in Section~\ref{Had_Spin}.
For now, we focus on their type-II property.
Let $A$ be an $\nbyn$ Potts model, that is, $A=-u^3I+u^{-1}(J-I)$ with
$(u^2+u^{-2})^2=n$.
For each $\epsilon \in \{1,-1\}$ and $\omega$ being a fourth root of 
$\epsilon$, the matrix
\begin{equation*}
W_{\epsilon} = 
\begin{pmatrix}
A & A & \omega H & -\omega H \\
A & A & -\omega H & \omega H \\
\epsilon \omega H^T & -\epsilon \omega H^T & A & A\\
-\epsilon \omega H^T & \epsilon \omega H^T & A & A\\
\end{pmatrix}
\end{equation*}
is a type-II matrix.

The Bose-Mesner algebras of the symmetric and non-symmetric 
association schemes in Section~\ref{ex_Schemes} contain $W_1$ and
$W_{-1}$, respectively.
In general, the Nomura algebras of these type-II matrices are
not equal to the Bose-Mesner algebras of the four-class association schemes
described in Section~\ref{ex_Schemes}.
For instance, when $n = 4$, $u=-i$ and 
\begin{equation*}
H=
\begin{pmatrix}
1&-1&1&1\\1&1&-1&1\\1&1&1&-1\\-1&1&1&1
\end{pmatrix},
\end{equation*}
$\Nom{W_{\epsilon}}$ has dimension $16$.
\section{Nomura Algebras of Two Type-II Matrices}
\label{N_II}

We examine the case where $A$ and $B$ are type-II matrices, and
study the interactions among the algebras $\N{A}{B}$,
$\dualN{A}{B}$, $\Nom{A}$, $\Nom{A^T}$, $\Nom{B}$ and $\Nom{B^T}$.
This section lays the groundwork for Chapters~\ref{JP} to \ref{InvJP}.

\begin{theorem}
\label{thm_T_ABC}
Let $A$, $B$ and $C$ be $\nbyn$ type-II matrices.
If $F \in \N{A}{B}$ and $G \in \N{\Sinv{B}}{C}$, then
$F \schur G \in \N{A}{C}$ and 
\begin{equation*}
\T{A}{C} (F \schur G) = n^{-1} \T{A}{B} (F)\ \T{\Sinv{B}}{C} (G).
\end{equation*}
\end{theorem}
{\sl Proof.}
Let $F'=\T{A}{B}(F)$ and $G'=\T{\Sinv{B}}{C} (G)$.
Applying Corollary~\ref{cor_T_AB_XDX}~(\ref{eqn_T_AB_XDX_d}) to
$F'=\T{A}{B}(F)$, we get
\begin{equation*}
\DXD{B^T}{\SinvT{B}}{F} = \XDX{(F')^T}{\inv{A}}{A^T},
\end{equation*}
which is equivalent to
\begin{equation*}
\DXD{F}{\invT{A}}{\Sinv{(\inv{A})}} = \XDX{\inv{(\SinvT{B})}}{\SinvT{B}}{(F')^T}.
\end{equation*}
Since $\Sinv{(\inv{A})}= nA^T$, $\invT{A}=n^{-1}\Sinv{A}$ and 
$\SinvT{B} = n \inv{B}$, the above equation equals
\begin{equation}
\label{eqn_T_ABC_1}
\DXD{F}{\Sinv{A}}{A^T} = \XDX{B}{\inv{B}}{(F')^T}.
\end{equation}
Similarly, applying Corollary~\ref{cor_T_AB_XDX}~(\ref{eqn_T_AB_XDX_d})
to $G'=\T{\Sinv{B}}{C} (G)$, we have
\begin{equation*}
\DXD{C^T}{\SinvT{C}}{G} = \XDX{(G')^T}{\inv{(\Sinv{B})}}{\SinvT{B}},
\end{equation*}
which becomes
\begin{equation*}
\DXD{C^T}{\SinvT{C}}{G} = \XDX{(G')^T}{B^T}{\inv{B}}.
\end{equation*}
Multiply the left-hand side of Equation~(\ref{eqn_T_ABC_1}) by
$\DXD{C^T}{\SinvT{C}}{G}$ and multiply its right-hand side by
$\XDX{(G')^T}{B^T}{\inv{B}}$. This gives
\begin{equation*}
\DXD{C^T}{\SinvT{C}}{G} \DXD{F}{\Sinv{A}}{A^T}
= \XDX{(G')^T}{B^T}{\inv{B}} \XDX{B}{\inv{B}}{(F')^T},
\end{equation*}
which simplifies to
\begin{equation*}
\D_{C^T} X_{\SinvT{C}} \D_{F \schur G} X_{\Sinv{A}} \D_{A^T} = 
n^{-1} X_{(F'G')^T}.
\end{equation*}
Hence
\begin{equation*}
\DXD{F \schur G}{\Sinv{A}}{A^T}=\XDX{C}{\inv{C}}{n^{-1} (F'G')^T}.
\end{equation*}
By Corollary~\ref{cor_T_AB_XDX}~(\ref{eqn_T_AB_XDX_c}),
we get
\begin{equation*}
\T{C}{A}(F \schur G) = n^{-1} (F'G')^T.
\end{equation*}
Therefore
\begin{equation*}
\T{A}{C} (F\schur G) = n^{-1} F'G'.  
\end{equation*}
\eop

Lemma~\ref{lem_T_A_Transpose} is a special case of the next lemma.
\begin{lemma}
\label{lem_Sinv_AB}
Let $A$ and $B$ be $\nbyn$ type-II matrices.
If $R \in \N{A}{B}$ then 
\begin{equation*}
R^T \in \N{\Sinv{A}}{\Sinv{B}}, \quad \text{and }
\T{\Sinv{A}}{\Sinv{B}}(R^T) = \T{A}{B}(R).
\end{equation*}
\end{lemma}
{\sl Proof.}
Let $S = \T{A}{B}(R)$.
By Corollary~\ref{cor_T_AB_XDX}~(\ref{eqn_T_AB_XDX_b}), we have
\begin{equation*}
\XDX{B^T}{A}{R^T} = \DXD{S^T}{B^T}{A}
\end{equation*}
which can be rewritten as
\begin{equation*}
\XDX{R^T}{\Sinv{A}}{\invT{B}} = \DXD{\Sinv{A}}{\invT{B}}{S^T}.
\end{equation*}
Now $\invT{B}=n^{-1} \Sinv{B}$ and using the Exchange Lemma, we obtain
\begin{equation*}
\XDX{R^T}{\Sinv{B}}{\Sinv{A}}= \DXD{\Sinv{B}}{\Sinv{A}}{S}.
\end{equation*}
So the result follows by Theorem~\ref{thm_T_AB_XDX}.  
\eop

The next two theorems are easy consequences of Theorem~\ref{thm_T_ABC}.
They describe some interactions among the maps $\T{A}{B}$, $\Theta_A$
and $\Theta_B$.
\begin{theorem}
\label{thm_T_AB}
Let $A$ and $B$ be $\nbyn$ type-II matrices.
If $F\in \Nom{A}$, $G \in \N{A}{B}$ and $H\in \Nom{B}$, then
$F \schur G$, and $G \schur H$ belong to $\N{A}{B}$ and 
\begin{eqnarray*}
\T{A}{B}(F \schur G) &=& n^{-1} \Theta_A(F)\ \T{A}{B}(G)\\
\T{A}{B}(G \schur H) &=& n^{-1} \T{A}{B}(G)\ \Theta_B(H)^T.
\end{eqnarray*}
\end{theorem}
{\sl Proof.}
The first equation results from applying Theorem~\ref{thm_T_ABC} to
the matrices $(A,\Sinv{A},B)$.

Applying the same theorem to $(A,B,B)$, we find that
\begin{eqnarray*}
\T{A}{B}(G \schur H) &=& n^{-1} \T{A}{B}(G)\ \T{\Sinv{B}}{B}(H)\\
&=& n^{-1} \T{A}{B}(G)\ \Theta_B(H)^T,
\end{eqnarray*}
because $\Theta_B(H) = \T{B}{\Sinv{B}}(H) = \T{\Sinv{B}}{B}(H)^T$.
\eop

\begin{theorem}
\label{thm_TA_TB}
Let $A$ and $B$ be $\nbyn$ type-II matrices.
If $F,G \in \N{A}{B}$, then $F \schur G^T \in \Nom{A} \cap \Nom{B}$ and 
\begin{eqnarray*}
\Theta_A(F \schur G^T) &=& n^{-1} \T{A}{B}(F)\ \T{A}{B}(G)^T\\
\Theta_B(F \schur G^T) &=& n^{-1} \T{A}{B}(F)^T\ \T{A}{B}(G).
\end{eqnarray*}
\end{theorem}
{\sl Proof.}
By Lemma~\ref{lem_Sinv_AB}, 
we know $G^T$ belongs to $\N{\Sinv{B}}{\Sinv{A}} = \N{\Sinv{A}}{\Sinv{B}}$.
By applying Theorem~\ref{thm_T_ABC} to the matrices $(A,B,\Sinv{A})$, we get
\begin{equation*}
\T{A}{\Sinv{A}}(F \schur G^T) = n^{-1} \T{A}{B}(F)\ \T{\Sinv{B}}{\Sinv{A}}(G^T).
\end{equation*}
Using Lemma~\ref{lem_Sinv_AB} again, we see that 
\begin{equation*}
\T{\Sinv{B}}{\Sinv{A}}(G^T)=\T{B}{A}(G) = \T{A}{B}(G)^T.
\end{equation*}
Hence the first equation holds.

We now apply Theorem~\ref{thm_T_ABC} to the matrices $(B,A,\Sinv{B})$, and
obtain
\begin{eqnarray*}
\Theta_B(F \schur G^T) &=& n^{-1} \T{B}{A}(F)\ \T{\Sinv{A}}{\Sinv{B}}(G^T)\\
&=& n^{-1} \T{A}{B}(F)^T\ \T{A}{B}(G).
\end{eqnarray*} 
\eop

The following is an important consequence of Theorems~\ref{thm_T_AB}
and \ref{thm_TA_TB}.  
It implies that if $\N{A}{B}$ contains a Schur-invertible matrix, then
$\N{A}{B}$, $\dualN{A}{B}$, $\Nom{A}$, $\Nom{B}$, $\Nom{A^T}$ and 
$\Nom{B^T}$ have the same dimension and $\Nom{A}=\Nom{B}$.
\begin{theorem}
\label{thm_Dim}
Let $A$ and $B$ be $\nbyn$ type-II matrices.
If $\N{A}{B}$ contains a Schur-invertible matrix $G$ and $H=\T{A}{B}(G)$, then
\begin{enumerate}
\item
$G \schur \Nom{A} = \N{A}{B}$ and $G^T \schur \N{A}{B} = \Nom{A}$.
\item
$\Nom{A^T}H=\dualN{A}{B}$ and $\dualN{A}{B} H^T = \Nom{A^T}$.
\item
$\Nom{B} = \Nom{A}$.
\item
$\Nom{B^T} = \inv{H} \Nom{A^T} H$.
\end{enumerate}
\end{theorem}
{\sl Proof.}
By Theorem~\ref{thm_T_AB}, we have $G \schur \Nom{A} \subseteq \N{A}{B}$
and $G \schur \Nom{B} \subseteq \N{A}{B}$.
Since $G$ is Schur invertible, the dimensions of $\Nom{A}$ and $\Nom{B}$
are less than or equal to $\dim{\N{A}{B}}$.
Similarly, Theorem~\ref{thm_TA_TB} implies 
$G^T \schur \N{A}{B}$ is a subset of  $\Nom{A}$ and $\Nom{B}$. 
The Schur-invertibility of $G^T$ implies that the dimension of
$\N{A}{B}$ is less than or equal to the dimensions of $\Nom{A}$ and
$\Nom{B}$.
As a result, we have $\Nom{A}=\Nom{B}=G^T \schur \N{A}{B}$ and
$G \schur \Nom{A} = \N{A}{B}$.

By the first equation of Theorem~\ref{thm_T_AB}, we have
$\Nom{A^T}H \subseteq \dualN{A}{B}$.
Since we have $\N{A}{B}=G \schur \Nom{A}$ for any $M \in \N{A}{B}$, 
there exists $M' \in \Nom{A}$ such that $M = M' \schur G$.
It follows from Theorem~\ref{thm_T_AB} that
\begin{eqnarray*}
\T{A}{B}(M) &=& \T{A}{B}(M' \schur G)\\
&=& n^{-1} \Theta_A(M') H
\end{eqnarray*}
As a result, we have $\dualN{A}{B} \subseteq \Nom{A^T} H$
and the first part of (b) follows.

Similarly, the first part of Theorem~\ref{thm_TA_TB} implies that 
$\dualN{A}{B} H^T \subseteq \Nom{A^T}$.
Because $\Nom{A}=G^T \schur \N{A}{B}$,
for all $N\in \Nom{A}$, there exists $N'\in \N{A}{B}$ such that
$N = N' \schur G^T$.  It follows from
the first equation of Theorem~\ref{thm_TA_TB} that
\begin{eqnarray*}
\Theta_A(N) &=& \Theta_A(N' \schur G^T)\\
&=&n^{-1}\T{A}{B}(N')H^T.
\end{eqnarray*}
So $\Nom{A^T} \subseteq \dualN{A}{B} H^T$ and we have
proved the rest of (b).

Using the same kind of argument,
 the second equations of Theorems~\ref{thm_T_AB} and 
\ref{thm_TA_TB}
imply that $\dualN{A}{B}=H\Nom{B^T}$ and $\Nom{B^T}=H^T \dualN{A}{B}$,
respectively.
Consequently, 
\begin{equation*}
\Nom{B^T}=\inv{H} \dualN{A}{B} =\inv{H} \Nom{A^T} H.
\end{equation*} 
\eop

Since $\Nom{A}$ is closed under the transpose, part~(a) of this theorem tells us
that if $\N{A}{B}$ contains a symmetric Schur-invertible matrix, then
$\N{A}{B}$ is also closed under the transpose.


\chapter{Jones Pairs}
\label{JP}

Given an $\nbyn$ symmetric spin model $W$, Jones \cite{VJ_Knot} 
defined endomorphisms of $\com^n \otimes \com^n$, 
$X_W$ and $\D_{\Sinv{W}}$, that provide a braid group representation.  
In Section~3.5 of the same paper,
Jones questioned the necessity of the type-II condition on $W$ for the
representation to give a link invariant.
This question motivates us to consider Jones' braid group 
representation without assuming type-II condition.
We extend Jones' idea to use endomorphisms $X_A$ and $\D_B$,
where $B$ may not equal $\Sinv{A}$.  
Jones pairs are defined in this process.

We devote this chapter to develop the theory of Jones pairs.
We introduce Jones pairs and a weaker version, called 
one-sided Jones pairs, in Sections~\ref{1JP} and \ref{1JP_Prop}.
In Section~\ref{JP_Braid}, we discuss the braid group representations 
obtained from Jones pairs.
In the remaining sections, we focus on the effect of the invertibility
of $B$ in a one-sided Jones pair $(A,B)$.
An important consequence is the equivalence of invertible Jones pairs
and four-weight spin models.
In the last section, we reproduce Jaeger's results on gauge equivalence with 
weaker assumptions.
Except for the three corollaries in Section~\ref{Nom_Inv_1JP}, due to
the author, all results in this chapter are joint work by Godsil,
Munemasa and the author \cite{CGM}.

\section{One-Sided Jones Pairs}
\label{1JP}

A pair of $\nbyn$ matrices $(A,B)$ is 
a \SL{one-sided Jones pair}\  if both $X_A$ and $\D_{B}$ are invertible 
and they satisfy
\begin{equation}
\label{eqn_ABA}
\XDX{A}{B}{A} = \DXD{B}{A}{B}.
\end{equation}
Equivalent to the definition in Section~\ref{Intro_Intro}, we say that
the pair $(A,B)$ is a \SL{Jones pair}\  if 
$(A,B)$ and $(A,B^T)$ are one-sided Jones pairs.

By Theorem~\ref{thm_T_AB_XDX}, an invertible matrix $A$ and
a Schur-invertible matrix $B$ form a one-sided Jones pair
if and only if 
\begin{equation*}
A \in \N{A}{B} \quad \text{and} \quad \T{A}{B}(A)=B.
\end{equation*}
The pair $(I,J)$ is an obvious example of one-sided Jones pair.
It is in fact a Jones pair because $J$ is symmetric.

Using the eigenvector approach, if $A$ is invertible and 
$B$ is Schur-invertible, then $(A,B)$ is 
a one-sided Jones pair if and only if 
\begin{equation*}
A \ (\Evector{A}{i}{B}{j})= B_{i,j} \ (\Evector{A}{i}{B}{j}),
\end{equation*}
for all $\range{i,j}{1}{n}$.
The $k$-th entry of both sides equal
\begin{equation}
\label{eqn_1JP_AB}
\sum_{x=1}^n A_{k,x} A_{x,i} B_{x,j} = B_{i,j} A_{k,i} B_{k,j}.
\end{equation}

Moreover, if $(A,B)$ is a Jones pair, replacing $B$ by $B^T$ in
Equation~(\ref{eqn_1JP_AB}) gives 
\begin{equation}
\label{eqn_1JP_ABT}
\sum_{x=1}^n A_{k,x} A_{x,i} B_{j,x} = B_{j,i} A_{k,i} B_{j,k}.
\end{equation}

Let $W$ be a spin model with loop variable $d$.
Then the type-III condition of $W$ is 
\begin{equation*}
\sum_{x=1}^n \frac{W_{k,x}W_{x,i}}{W_{j,x}}
= d \frac{W_{k,i}}{W_{j,i}W_{j,k}},
\end{equation*}
which is exactly Equation~(\ref{eqn_1JP_AB}) with $A=d^{-1}W$ and 
$B = \Sinv{W}$.
As a result, all spin models give one-sided Jones pairs.

Let $(W_1,W_2, W_3, W_4;d)$ be a four-weight spin model.
Then its type-III conditions are
\begin{eqnarray*}
\sum_{x=0}^n (W_1)_{k,x} (W_1)_{x,i} (W_4)_{x,j} &=&
d (W_4)_{i,j} (W_1)_{k,i} (W_4)_{k,j} \quad \text{and}\\
\sum_{x=0}^n (W_1)_{k,x} (W_1)_{x,i} (W_4)_{j,x} &=&
d (W_4)_{j,i} (W_1)_{k,i} (W_4)_{j,k}.
\end{eqnarray*}
These equations are the same as Equations~(\ref{eqn_1JP_AB}) and 
(\ref{eqn_1JP_ABT}) with $A=d^{-1}W_1$ and $B=W_4$, respectively.
So we obtain a Jones pair from every four-weight spin model.

In addition, we can build Jones pairs using Kronecker product.
Suppose $(A,B)$ and $(A',B')$ are one-sided Jones pairs.
By Lemma~\ref{lem_N_AB_Tensor}, we have
\begin{equation*}
A \otimes A' \in \N{A\otimes A'}{B \otimes B'}.
\end{equation*}
From the proof of this lemma, we see that the eigenvectors can be
written as $(\Evector{A}{i}{B}{j}) \otimes \ (\Evector{A'}{h}{B'}{k})$.
It follows that 
\begin{equation*}
\T{A\otimes A'}{B \otimes B'}(A\otimes A') =
\T{A}{B}(A) \otimes \T{A'}{B'}(A') = B \otimes B'.
\end{equation*}
Hence $(A\otimes A', B\otimes B')$ is also a one-sided Jones pair.
As a result, there exist infinitely many one-sided Jones
pairs.

Suppose $(A,B)$ is a one-sided Jones pair.
Fix any $j\in \{1,\ldots, n\}$, define $D_j$ to be the diagonal
matrix with its $ii$-entry equals $B_{i,j}$.
Then $(D_j J)_{h,k} =B_{h,j}$ 
and $(D_jJ) e_k = B e_j$,
for all $\range{k}{1}{n}$.
As a result, we get
\begin{eqnarray*}
A \ (\Evector{A}{h}{(D_jJ)}{k}) &=& A \ (\Evector{A}{h}{B}{j})\\
&=& B_{h,j} \ (\Evector{A}{h}{B}{j})\\
&=& (D_jJ)_{h,k} \ (\Evector{A}{h}{(D_jJ)}{k})
\end{eqnarray*}
and $(A,D_jJ)$ is a one-sided Jones pair.

Now $(I,J)$, $(A,D_jJ)$ and $(A\otimes I,B \otimes J)$
are the only known examples of one-sided Jones pair with
the second matrix being non-invertible.  
In Section~\ref{4wt}, we see that the Jones pairs with the second matrix
invertible are equivalent to four-weight spin models.
As a result, we would be very excited to see
any new example of one-sided Jones pairs with the second matrix
non-invertible.

\section{Properties of One-Sided Jones Pairs}
\label{1JP_Prop}

This section lists some useful properties of one-sided Jones pairs.

\begin{lemma}
\label{lem_1JP_Col_Sum}
If $(A,B)$ is a one-sided Jones pair, then $B^T J = \tr(A)J$.
Furthermore, if $(A,B)$ is a Jones pair, then the columns and the rows
of $B$ sum to $\tr(A)$.
\end{lemma}
{\sl Proof.}
Since $A$ is invertible and $B$ is Schur invertible,
the set
\begin{equation*}
\{\Evector{A}{1}{B}{r},
\Evector{A}{2}{B}{r},\ldots,\Evector{A}{n}{B}{r}\}
\end{equation*}
is a basis of $\com^n$,
for each $\range{r}{1}{n}$. 
As a result, each column of $B$ contains 
all eigenvalues of $A$ and $B^TJ = \tr(A)J$.

Similarly, if $(A,B^T)$ is also a one-sided Jones pair, then
$BJ=\tr(A)J$.  
So if $(A,B)$ is a Jones pair, then $B^T J = BJ =\tr(A)J$.
\eop

\begin{lemma}
\label{lem_1JP_Eq}
Suppose $(A,B)$ is a one-sided Jones pair.  So is
\begin{enumerate}
\item
$(A^T,B)$
\item
$(\inv{A},\Sinv{B})$
\item
$(\inv{D}AD,B)$, for any invertible diagonal matrix $D$
\item
$(A,BP)$, for any permutation matrix $P$
\item
$(PA\inv{P},PB\inv{P})$, for any permutation matrix $P$
\item
$(\lambda A, \lambda B)$, for any non-zero complex number $\lambda$
\end{enumerate}
\end{lemma}
{\sl Proof.}
\begin{enumerate}
\item
Taking the transpose of both sides of (\ref{eqn_ABA}), we get
$\XDX{A^T}{B}{A^T}=\DXD{B}{A^T}{B}$.
\item
Since $\inv{X_A} = X_{\inv{A}}$ and $\inv{\D_B}=\D_{\Sinv{B}}$,
inverting both sides of (\ref{eqn_ABA}) gives
$\XDX{\inv{A}}{\Sinv{B}}{\inv{A}} = \DXD{\Sinv{B}}{\inv{A}}{\Sinv{B}}$.
\item
Note that for any diagonal matrix $D$, $X_D=\D_{DJ}$ and 
$X_D\D_M = \D_M X_D$.
So
\begin{eqnarray*}
\XDX{\inv{D}AD}{B}{\inv{D}AD} &=&
X_{\inv{D}} X_A (\XDX{D}{B}{\inv{D}}) X_A X_D\\
&=& X_{\inv{D}} (\XDX{A}{B}{A}) X_D\\
&=& X_{\inv{D}} \DXD{B}{A}{B} X_D\\
&=& \D_B (X_{\inv{D}} X_A X_D) \D_B\\
&=& \DXD{B}{\inv{D}AD}{B}.
\end{eqnarray*}
The third equality results from $(A,B)$ being
a one-sided Jones pair.
\item
By Lemma~\ref{lem_Perm},
we have 
\begin{equation*}
A \in \N{A}{B} = \N{A}{BP} 
\end{equation*}
and
\begin{equation*}
\T{A}{BP}(A)=\T{A}{B}(A)P=BP.
\end{equation*}
\item
Using the same lemma with $Q=R=\inv{P}$, we have
\begin{equation*}
PA\inv{P} \in \N{PA\inv{P}}{PB\inv{P}}
\end{equation*}
and
\begin{equation*}
\T{PA\inv{P}}{PB\inv{P}}(PA\inv{P}) = P \T{A}{B}(A) \inv{P} = PB\inv{P}.
\end{equation*}
\item
Replacing $(A,B)$ by $(\lambda A, \lambda B)$ is equivalent to multiplying
both sides of (\ref{eqn_ABA}) by $\lambda^3$.
\end{enumerate} 
\eop

If $(A,B)$ is a one-sided Jones pair, then $A$ is invertible
and $B$ is Schur invertible and $B = \T{A}{B}(A)$.
So Corollary~\ref{cor_T_AB_XDX} provides three equivalent forms of
Equation~(\ref{eqn_ABA}).  
We give one more useful reformulation of Equation~(\ref{eqn_ABA}) below
when $A$ is also Schur invertible and $B$ is also invertible.
\begin{lemma}
\label{lem_XDX_A_B}
If $A$ and $B$ are both invertible and Schur-invertible, then
$(A,B)$ is a one-sided Jones pair if and only if
\begin{equation*}
\DXD{\Sinv{A}}{A}{A} = \XDX{B}{B^T}{\inv{B}}.
\end{equation*}
\end{lemma}
{\sl Proof.}
Applying the Exchange Lemma to Equation~(\ref{eqn_ABA}) yields
\begin{equation*}
\XDX{A}{A}{B} = \DXD{A}{B}{B^T}.
\end{equation*}
Since both $\Sinv{A}$ and $\inv{B}$ exist, 
we get the relation in the lemma immediately.  
\eop
\section{Braid Group Representations}
\label{JP_Braid}

Given a Jones pair $(A,B)$,
we demonstrate Jones' method of constructing braid group representations
from $X_A$ and $\D_B$.  
The \SL{braid group}\  $\Br{m}$ on $m$ strands is 
generated by $\seq{\sigma}{1}{m-1}$
satisfying 
\begin{enumerate}
\item
For all $\range{i}{1}{m-2}$,
\begin{equation}
\label{eqn_Br1}
\sigma_i \sigma_{i+1} \sigma_i = \sigma_{i+1} \sigma_i \sigma_{i+1}.
\end{equation}
\item
For all $|i-j| \geq 2$,
\begin{equation*}
\sigma_j \sigma_k = \sigma_k \sigma_j.  
\end{equation*}
\end{enumerate}

Let $k=\lceil \frac{m}{2} \rceil$ and 
let $V$ denote the vector space $\com^n$.
Given a pair of $\nbyn$ matrices $(A,B)$,  we map the generators of $\Br{m}$ to
the endomorphisms $\seq{g}{1}{m-1}$ of $V^{\otimes k}$ as follows:
\begin{eqnarray*}
g_{2h-1}(e_{r_1} \otimes \ldots \otimes e_{r_k}) &=& 
e_{r_1} \otimes \ldots \otimes (A e_{r_h}) \otimes \ldots \otimes e_{r_k},\\
g_{2h}(e_{r_1} \otimes \ldots \otimes e_{r_k}) &=& 
B_{{r_h},{r_{h+1}}}(e_{r_1}\otimes \ldots \otimes e_{r_k}).
\end{eqnarray*} 
When $|i-j|\geq 2$, $g_i$ and $g_j$ act on different tensor factors of 
$V^{\otimes k}$, so they commute.
Note that $g_{2h-1}$, $g_1$ and $g_3$ have the same action 
on the $h$-th, $1$-st and $2$-nd tensor factors of $V^{\otimes k}$,
respectively.
Similarly, the action of $g_{2h}$ on the $h$-th and the $(h+1)$-th 
tensor factors of $V^{\otimes k}$ is the same as the action of $g_2$
on the first and second tensor factors.
Consequently, showing $g_1 g_2 g_1= g_2 g_1 g_2$ and 
$g_2 g_3 g_2=g_3 g_2 g_3$ is sufficient to prove that (\ref{eqn_Br1})
holds for all $\range{i}{1}{m-2}$.
\begin{lemma}
\label{lem_Br_XDX}
The relation $g_1 g_2 g_1 = g_2 g_1 g_2$ holds if and only if 
$(A,B)$ is a one-sided Jones pair.
\end{lemma}
{\sl Proof.}
We use the isomorphism $\phi:V \otimes V \rightarrow \mat{n}{\com}$ 
which maps $e_i \otimes e_j$ to $e_i e_j^T$.
We get
\begin{eqnarray*}
\phi(g_1(e_i\otimes e_j)) &=& \phi(Ae_i \otimes e_j)\\
&=& Ae_i e_j^T\\
&=& X_A(\Eij).
\end{eqnarray*}
So we have $X_A = \phi g_1 \inv{\phi}$.
Similarly,
\begin{equation*}
\phi(g_2(e_i\otimes e_j)) = B_{i,j} e_i e_j^T = \D_B(\Eij)
\end{equation*}
and $\D_B = \phi g_2 \inv{\phi}$.
Consequently,
\begin{equation*}
\phi g_1 g_2 g_1 \inv{\phi} = \XDX{A}{B}{A} 
\end{equation*}
and
\begin{equation*}
\phi g_2 g_1 g_2 \inv{\phi} = \DXD{B}{A}{B}.
\end{equation*}
So the result follows.  
\eop

Since $\Br{3}$ has only two generators, every one-sided Jones pair gives
a representation of $\Br{3}$.

\begin{lemma}
\label{lem_Br_YDY}
The relation $g_2 g_3 g_2 = g_3 g_2 g_3$ holds if and only if 
$(A,B^T)$ is a one-sided Jones pair.
\end{lemma}
{\sl Proof.}
For any $\nbyn$ matrix $C$, let $Y_C$ be the endomorphism of $\mat{n}{\com}$
defined as $Y_C(M) = MC^T$.
Using $\phi$ as above, we have
\begin{equation*}
\phi(g_3(e_i \otimes e_j)) = \phi(e_i \otimes Ae_j) = e_i(Ae_j)^T = Y_A(\Eij)
\end{equation*}
and hence $Y_A = \phi g_3 \inv{\phi}$.
As a result, the relation $g_2 g_3 g_2 = g_3 g_2 g_3$ holds if and only if
$\D_B Y_A \D_B = Y_A \D_B Y_A$.
We can write this equation as
\begin{equation*}
(B \schur (\Eij A^T)) A^T = B \schur ((B \schur \Eij)A^T),
\end{equation*}
for all $\range{i,j}{1}{n}$.
Taking the transpose of each side, we get
\begin{equation*}
A(B^T \schur (A\Eij)) = B^T \schur(A(B^T \schur \Eij)),
\end{equation*}
which is equivalent to $(A,B^T)$ being a one-sided Jones pair.  
\eop

\begin{corollary}
\label{cor_Braid}
Every Jones pair gives a representation of $\Br{m}$.
\end{corollary}


Suppose $(A,B)$ is a Jones pair.
Let $\seq{g}{1}{m-1}$ be the representation of $\Br{m}$ described above.
We proved in \cite{CGM} that if $A$ and $B$ satisfy
\begin{eqnarray}
\label{eqn_Link}
A \schur I =& \inv{A} \schur I =& \frac{1}{\sqrt{n}} I, \quad \text{and}\\ 
BJ =& \Sinv{B} J = & \sqrt{n}J, \nonumber
\end{eqnarray}
then for any $h$ generated by $\seq{g}{1}{m-2}$.
We have
\begin{eqnarray*}
\tr(h g_{m-1}) &=& \frac{1}{n}\tr(h) \tr(A), \quad \text{and,}\\
\tr(h g_{m-1}^{-1}) &=& \frac{1}{n}\tr(h) \tr(A)^{-1}.\\
\end{eqnarray*}
These conditions are sufficient for a Jones pair to give a link
invariant in the form of the trace of the endomorphisms generated by
$\seq{g}{1}{m-1}$.
In the next section, we will see that if $B$ is invertible, 
then (\ref{eqn_Link}) holds.

\section{Invertibility}
\label{Invert}

By definition, if $(A,B)$ is a one-sided Jones pair then $A$ is invertible and
$B$ is Schur invertible.  We call $(A,B)$ 
an \SL{invertible one-sided Jones pair}\ 
if $\Sinv{A}$ and $\inv{B}$ also exist.

In this section, we show that for a one-sided Jones pair $(A,B)$,
the invertibility of $B$ implies the Schur-invertibility of $A$.
In fact, assuming $B$ is invertible is a more stringent condition
than one may expect at first.
We prove that the invertibility of $B$ implies the type-II condition on both 
$A$ and $B$;
This will be used in Section~\ref{4wt} to justify our assertion that
invertible Jones pairs and four-weight spin models are equivalent concepts.

\begin{theorem}
\label{thm_B_Inv}
Let $(A,B)$ be a one-sided Jones pair.
If $B$ is invertible then both $A$ and $B$ are type-II matrices.
Moreover, $A$ has constant diagonal and $B$ has constant row sums.
\end{theorem}
{\sl Proof.}
Since $\inv{A}$ has the same eigenvectors as $A$ and their eigenvalues 
with respect to the same eigenvector are reciprocal to each other, 
we have $\T{A}{B}(\inv{A}) = \Sinv{B}$.
Now applying Corollary~\ref{cor_T_AB_XDX}~(\ref{eqn_T_AB_XDX_b}) to
$\T{A}{B}(\inv{A})=\Sinv{B}$, we get
\begin{equation*}
\XDX{B^T}{A}{\invT{A}} = \DXD{\SinvT{B}}{B^T}{A}.
\end{equation*}
Evaluating this equation at $I$ yields
\begin{eqnarray*}
B^T(A \schur \invT{A}) &=& \SinvT{B} \schur (B^T (A\schur I))\\
&=& (\SinvT{B} \schur B^T) (A \schur I)\\
&=& J (A \schur I).
\end{eqnarray*}
Since $(A,B)$ is a one-side Jones pair, $B^TJ=\tr(A)J$.  So 
\begin{equation*}
A\schur \invT{A} = \invT{B} J (A\schur I) = \tr(A)^{-1} J (A \schur I).
\end{equation*}
The sum of the $i$-th column of $A\schur \invT{A}$ equals
\begin{equation*}
\sum_{k=1}^n A_{k,i} (\inv{A})_{i,k} = 1.
\end{equation*}
The $i$-th column of $\tr(A)^{-1} J (A \schur I)$ sums to 
$n \tr(A)^{-1}A_{i,i}$.
We have $A_{i,i} = n^{-1}\tr(A)$
and $A\schur I = n^{-1}\tr(A)I$.
Therefore $A \schur \invT{A} = n^{-1} J$ and 
it follows that $A$ is a type-II matrix.

By Lemma~\ref{lem_XDX_A_B}, we have
\begin{equation*}
\DXD{\Sinv{A}}{A}{A} = \XDX{B}{B^T}{\inv{B}}.
\end{equation*}
Evaluating both sides at $I$ gives
\begin{equation*}
\Sinv{A} \schur (A(A \schur I)) = B(B^T \schur \inv{B})
\end{equation*}
which leads to
\begin{equation*}
\inv{B}(\Sinv{A} \schur A) (A \schur I) = B^T \schur \inv{B},
\end{equation*}
and finally
\begin{equation*}
\frac{\tr(A)}{n} \inv{B} J = B^T \schur \inv{B}.
\end{equation*}
The sum of the $i$-th row of $B^T \schur \inv{B}$ equals
\begin{equation*}
\sum_{k=1}^n B_{k,i}(\inv{B})_{i,k}=1.
\end{equation*}
So we have $\tr(A) \sum_{k=1}^n (\inv{B})_{i,k} = 1$, 
or equivalently, $BJ=\tr(A)J$.
Consequently, we get
\begin{equation*}
B^T \schur \inv{B} = n^{-1} J,
\end{equation*}
which is equivalent to the type-II condition on $B$.  
\eop
  
\begin{theorem}
\label{thm_A_Sinv}
Let $(A,B)$ be a Jones pair.
If $A$ is Schur invertible then $B$ is invertible.
\end{theorem}
{\sl Proof.}
By Theorem~\ref{thm_T_AB_XDX},
$\T{A}{B}(\inv{A}) = \Sinv{B}$ is equivalent to
\begin{equation*}
\XDX{\inv{A}}{B}{A} = \DXD{B}{A}{\Sinv{B}}.
\end{equation*}
Applying the Exchange Lemma, we get
\begin{equation*}
\XDX{\inv{A}}{A}{B} = \DXD{A}{B}{\SinvT{B}}.
\end{equation*}
Evaluating both sides at $J$ yields
\begin{equation*}
\inv{A} (A \schur BJ) = A \schur (B(\SinvT{B} \schur J)).
\end{equation*}
Since $(A,B^T)$ is a one-sided Jones pair, 
we have $BJ=\tr(A)J$, by Lemma~\ref{lem_1JP_Col_Sum}. 
So
\begin{equation*}
\tr(A)\ \inv{A} (A \schur J) = A \schur (B \SinvT{B}).
\end{equation*}
The left-hand side equals $\tr(A) I$, so 
\begin{equation*}
B \SinvT{B} = \tr(A)\ (\Sinv{A} \schur I)
\end{equation*}
and since $\inv{(\Sinv{A} \schur I)} = A \schur I$,
\begin{equation*}
\inv{B} = \tr(A)^{-1} \SinvT{B} (A \schur I).
\end{equation*} 
\eop

These two theorems tell us that if $(A,B)$ is a Jones pair,
then $A$ is Schur invertible if and only if $B$ is invertible.
In this case, both $A$ and $B$ are type-II matrices with
$A \schur I = \tr(A) n^{-1} I$ and $BJ = \tr(A) J = B^T J$.
So the Jones pair
\begin{equation*}
\left( \frac{\sqrt{n}}{\tr(A)} A, \frac{\sqrt{n}}{\tr(A)} B \right)
\end{equation*}
satisfies the conditions in (\ref{eqn_Link}).
Thus it gives a link invariant.

If $(A,B)$ is a one-sided Jones pair, then it follows immediately from
the definition of one-sided Jones pair that $A \in \N{A}{B}$.
The following theorem investigates the opposite direction.
It extends Jaeger, Matsumoto and Nomura's result \cite{JMN}
which says that if $A\in \Nom{A}$ then $A$ is a spin model up to scalar
multiplication.
\begin{theorem}
\label{thm_A_in_N_AB}
Let $A$ and $B$ be $\nbyn$ type-II matrices.
Suppose $A\schur I = aI$ and $B^T J =bJ$ for some non-zero $a,b \in \com$.
If $A \in \N{A}{B}$ then $(A,ab^{-1}nB)$ is a one-sided Jones pair.
\end{theorem}
{\sl Proof.}
Since $A$ and $\inv{A}$ share the same set of eigenvectors,
the matrix $A$ belongs to $\N{A}{B}$ if and only if 
$\inv{A} = n^{-1} \SinvT{A}$ belongs to $\N{A}{B}$.
Using the formula of $\T{A}{B}$ in Lemma~\ref{lem_T_AB_Formula}, we get
\begin{eqnarray*}
\T{A}{B}(\inv{A}) &=& 
\Sinv{B} \schur (\inv{(A \schur I)} (A^T \schur n^{-1} \SinvT{A})B)\\
&=& \frac{1}{n} \Sinv{B} \schur ((a^{-1}I)JB)\\
&=& \frac{b}{an} \Sinv{B} \schur J\\
&=& \frac{b}{an} \Sinv{B}.
\end{eqnarray*}
Now an eigenvalue of $A$ is the reciprocal of the eigenvalue of
$\inv{A}$ with respect to the same eigenvector, we conclude that
\begin{equation*}
\T{A}{B}(A) = \frac{an}{b} B.
\end{equation*}
Hence
\begin{equation*}
\XDX{A}{\frac{an}{b}B}{A} = \DXD{\frac{an}{b}B}{A}{\frac{an}{b}B}.
\end{equation*} 
\eop

\section{Nomura Algebras of Invertible Jones Pairs}
\label{Nom_Inv_1JP}

If $(A,B)$ is an invertible one-sided Jones pair, then $A$ and $B$ are
type II.  From $A$ and $B$, we construct several Nomura algebras such as
$\Nom{A}$, $\Nom{B}$, $\N{A}{B}$, $\N{A}{B^T}$, et~cetera.
In this section, we investigate the relations among these algebras.

\begin{theorem}
\label{thm_Inv_1JP}
Let $(A,B)$ be an invertible one-sided Jones pair.
Then 
\begin{equation*}
\Nom{A} = \Nom{A^T} = \Nom{B} \quad \text{and} \quad
\dualN{A}{B}=\dualN{A^T}{B}.
\end{equation*}
\end{theorem}
{\sl Proof.}
Since $A\in \N{A}{B}$ is Schur invertible, 
by Theorem~\ref{thm_Dim}~(c), we have
\begin{equation*}
\Nom{A} = \Nom{B}.
\end{equation*}
By Lemma~\ref{lem_1JP_Eq}~(a), $(A^T,B)$ is also a one-sided Jones pair.
Applying Theorem~\ref{thm_Dim}~(c) on $A^T \in \N{A^T}{B}$ leads to
$\Nom{B}=\Nom{A^T}$.  

Since both $(A,B)$ and $(A^T,B)$ are one-sided Jones pairs, we have
\begin{equation*}
\T{A}{B}(A)=B=\T{A^T}{B}(A^T).
\end{equation*}
We deduce from Theorem~\ref{thm_Dim}~(b) that
\begin{equation*}
\dualN{A}{B}=\Nom{A^T}B, \quad \text{and}\quad \dualN{A^T}{B}=\Nom{A}B.
\end{equation*} 
But $\Nom{A}=\Nom{A^T}$ and so $\dualN{A}{B}=\dualN{A^T}{B}$.
\eop

Note that Theorem~\ref{thm_TA_TB} applied to $(A,B)$ yields
\begin{equation*}
\Theta_A(A\schur A^T) = n^{-1}BB^T,
\end{equation*}
which indicates that $\Nom{A}$ is non-trivial in most cases.

In the rest of this section, we examine the relations among
different Nomura algebras and their duality maps constructed from 
an invertible Jones pair.
\begin{theorem}
\label{thm_Inv_JP}
Let $(A,B)$ be an invertible Jones pair.
Then
\begin{equation*}
\Nom{A} = \Nom{A^T} = \Nom{B} = \Nom{B^T}.
\end{equation*}
Moreover, the dualities $\Theta_A$, $\Theta_B$ satisfy
$\Theta_B(F)^T = \inv{B} \Theta_A(F) B$ for all $F \in \Nom{A}$.
\end{theorem}
{\sl Proof.}
Since $(A,B)$ is an invertible one-sided Jones pair, the previous theorem
gives $\Nom{A}=\Nom{A^T}=\Nom{B}$.
Similarly, since $(A,B^T)$ is also an invertible one-sided Jones pair, we get
$\Nom{A}=\Nom{B^T}$ using Theorem~\ref{thm_Dim}~(c).  
So the first part of the theorem holds.
For any $F \in \Nom{A}=\Nom{B}$, Theorem~\ref{thm_T_AB} gives
\begin{eqnarray*}
\T{A}{B}(F\schur A) &=& n^{-1} \Theta_A(F) B,\\
\T{A}{B}(A \schur F) &=& n^{-1} B \Theta_B(F)^T.
\end{eqnarray*}
Hence we get $\Theta_B(F)^T = \inv{B}\ \Theta_A(F) B$.  
\eop

\begin{corollary}
\label{cor_N_AB_ABT}
If $(A,B)$ is an invertible Jones pair, then
\begin{equation*}
\N{A}{B^T}=\N{A}{B} \quad \text{and} \quad
\dualN{A}{B^T} = \inv{B} \dualN{A}{B}B^T.
\end{equation*}
\end{corollary}
{\sl Proof.}
Theorem~\ref{thm_Dim}~(a) applied to 
$A$ in $\N{A}{B}$ and $A$ in $\N{A}{B^T}$ gives
\begin{equation*}
\N{A}{B}=A \schur \Nom{A}=\N{A}{B^T}.
\end{equation*} 

We now prove the second equality.
Now both $A$ and $B$ are type-II matrices, by Corollary~\ref{cor_T_A},
$\dualNom{A}=\Nom{A^T}$ and $\dualNom{B} = \Nom{B^T}$.
So the second part of Theorem~\ref{thm_Inv_JP} says that
\begin{equation*}
\Nom{B^T} = \inv{B} \Nom{A^T} B.
\end{equation*}
But $(A,B)$ is an invertible Jones pair. So using the same theorem, 
we have $\Nom{A^T}=\Nom{B^T}$.  
Hence $\Nom{A^T}=\inv{B}\Nom{A^T}B$, or equivalently 
$B \Nom{A^T} = \Nom{A^T}B$.
Applying Theorem~\ref{thm_Dim}~(b) to $\T{A}{B}(A)=B$, we get
\begin{equation*}
\dualN{A}{B} = \Nom{A^T}B=B\Nom{A^T}.
\end{equation*}
Similarly, applying the same theorem to $\T{A}{B^T}(A)=B^T$ yields
\begin{eqnarray*}
\dualN{A}{B^T} &=& \Nom{A^T}B^T\\
&=& \inv{B} \dualN{A}{B} B^T.
\end{eqnarray*}
\eop

We will need the following lemma in the computation in
Chapter~\ref{SpinModels}.
\begin{lemma}
\label{lem_AB_ABT}
If $(A,B)$ is an invertible Jones pair and $R \in \N{A}{B}$,
then
\begin{equation*}
\T{A}{B}(R)B^T = \T{A}{B^T}(R)B.
\end{equation*}
\end{lemma}
{\sl Proof.}
By Theorem~\ref{thm_TA_TB}, for any $R$ in $\N{A}{B}$, 
we have 
\begin{equation*}
\Theta_A(R\schur A^T) = n^{-1} \T{A}{B}(R) B^T.
\end{equation*}
Since $\N{A}{B}=\N{A}{B^T}$, 
$R$ also belongs to $\N{A}{B^T}$.
The same theorem tells us that
\begin{equation*}
\Theta_A(R \schur A^T)= n^{-1} \T{A}{B^T}(R)B.
\end{equation*}
Thus we have the equation in the lemma.
\eop

\section{Four-Weight Spin Models}
\label{4wt}

Bannai and Bannai \cite{BB_4wt} generalized spin models to four-weight
spin models, and showed that the partition functions of
four-weight spin models provide invariants for oriented links.
In this section, we show that four-weight spin models are equivalent to
invertible Jones pairs.

As defined in Section~\ref{Intro_Spin},
a \SL{four-weight spin model}\  is a $5$-tuple
$(W_1,W_2,W_3,W_4;d)$ with $d=\pm \sqrt{n}$ satisfying
\begin{enumerate}[(a)]
\item
There exists non-zero scalar $a$ such that
\begin{eqnarray*}
W_1 \schur I = aI, & \quad & W_3 \schur I = a^{-1}I,\\
W_2 J = W_2^T J = da^{-1}J, & \quad & W_4 J = W_4^TJ=daJ.
\end{eqnarray*}
\item
The matrices $W_1, W_2, W_3$ and $W_4$ are type II and
\begin{equation*}
W_3=\SinvT{W_1}, \quad W_2=\SinvT{W_4}.
\end{equation*}
\item
\begin{eqnarray}
\label{eqn_4wt_1}
\sum_{x=1}^n (W_1)_{a,x} (W_1)_{x,b} (W_4)_{c,x} &=& 
d (W_1)_{a,b} (W_4)_{c,a} (W_4)_{c,b},\\
\label{eqn_4wt_2}
\sum_{x=1}^n (W_1)_{x,a} (W_1)_{b,x} (W_4)_{x,c} &=& 
d (W_1)_{b,a} (W_4)_{a,c} (W_4)_{b,c}.
\end{eqnarray}
\end{enumerate}

Using Theorems~\ref{thm_B_Inv} and \ref{thm_A_Sinv}, it is almost
immediate that invertible Jones pairs are the same as four-weight spin
models.
\begin{theorem}
\label{thm_4wt}
Let $A,B \in \mat{n}{\com}$ and $d^2=n$.
Then the following are equivalent.
\begin{enumerate}
\item
$(A,B)$ is an invertible Jones pair.
\item
$(dA,n\inv{B},d\inv{A},B;d)$ is a four-weight spin model.
\end{enumerate}
\end{theorem}
{\sl Proof.}
If $(dA, n\inv{B}, d \inv{A},B;d)$ is a four-weight spin model, then as
we have seen in Section~\ref{1JP}, $(A,B)$ is a Jones pair where $A$ and
$B$ are both type II.

Conversely, suppose $(A,B)$ is an invertible Jones pair.
By Theorem~\ref{thm_B_Inv}, both $A$ and $B$ are type-II matrices.
Moreover, $A$ has constant diagonal and $B$ has constant row sum and
column sum.
Let $W_1=dA$ and $W_4=B$.
Considering Equation~(\ref{eqn_1JP_AB}) with $k=b$, $i=a$, and $j=c$,
we see that $(A,B)$ is a one-sided Jones pair 
if and only if Condition~(\ref{eqn_4wt_2}) holds.
Similarly, using the same equation, $(A^T,B^T)$ is a one-sided Jones pair 
if and only if Condition~(\ref{eqn_4wt_1}) holds.
By Lemma~\ref{lem_1JP_Eq}~(a), we find that
$(A^T,B^T)$ is a one-sided Jones pair if and only if $(A,B^T)$ is also 
a one-sided Jones pair.
As a result, $(dA, n \inv{B}, d\inv{A},B;d)$ is a four-weight spin
model.
\eop

We will see in Theorem~\ref{thm_JP_Spin} that $W$ is a spin model 
if and only if $(d^{-1}W,\SinvT{W})$ is an invertible Jones pair.  
Therefore it is equivalent to
$(W,W,\SinvT{W},\SinvT{W};d)$ being a four-weight spin model.

In \cite{BB_4wt}, Bannai and Bannai studied three types of four-weight spin models:
Jones type, pseudo-Jones type and Hadamard type.
They correspond respectively to the following types of Jones pairs
\begin{enumerate}
\item
$(A,d\SinvT{A})$, where $A$ is a spin model,
\item
$(A,dA)$,
\item
$(A,B)$, where one of $A, B$ is a Hadamard matrix.
\end{enumerate}

Using $\nbyn$ Hadamard matrices satisfying certain conditions,
Yamada constructed $n^2 \times n^2$ four-weight spin models of 
pseudo-Jones type and symmetric Hadamard type. 
Using tensor products, her construction gives infinite families of four-weight spin models
of both types.  For details, see \cite{Yamada1}.

In \cite{Ban1}, Bannai extended work by Guo and proved that if 
$(W_1,W_2,W_3,W_4;d)$ is a four-weight spin model, then for
$i=1,\ldots,4$,
\begin{equation*}
\Nom{W_1}=\Nom{W_i}=\Nom{W_i^T}.
\end{equation*}
Note that by Theorem~\ref{thm_4wt}, the pair $(d^{-1} W_1,W_4)$ is
an invertible Jones pair.  
As shown in Theorem~\ref{thm_Inv_1JP}, we are able to obtain part of the
above equations, assuming a weaker condition on $(d^{-1} W_1,W_4)$, that
is, it is a one-sided Jones pair.

\section{Odd Gauge Equivalence}
\label{Odd_Gauge}

Two four-weight spin models 
$(W_1,W_2,W_3,W_4;d)$ and $(W'_1,W'_2,W'_3,W'_4;d)$ 
are \SL{gauge equivalent}\  if there exist an invertible diagonal matrix $D$, 
a permutation matrix $P$ and a non-zero scalar $c$ such that
\begin{eqnarray*}
W'_1 &= cDW_1\inv{D}, \quad W'_3&=c^{-1} DW_3\inv{D},\\
W'_2 &= c^{-1} \inv{P}W_2, \quad  W'_4&=cW_4P.
\end{eqnarray*}

Jaeger proved that gauge-equivalent spin models give the same invariant
(see Proposition~11 in \cite{J_4wt}).
In the same paper, he showed that $W'_2=W_2$ and $W'_4=W_4$ if and only if
$W'_1 = D W_1 \inv{D}$ and $W'_3=D W_3 \inv{D}$ for some invertible diagonal matrix $D$.
In this case, we say that the two four-weight spin models are related by 
an \SL{odd gauge transformation}.
He also proved that $W'_1=W_1$ and $W'_3=W_3$ if and only if $W'_2 = \inv{P}W_2$ and 
$W'_4=W_4P$ for some permutation matrix $P$.  The two four-weight spin models are said
to be related by an \SL{even gauge transformation}.


In Sections~\ref{Odd_Gauge} and \ref{Even_Gauge}, we extend 
Jaeger's result on gauge equivalence of four-weight spin models to
invertible one-sided Jones pairs.
In particular, the last lemma in this section allows us to consider only 
the invertible Jones pairs with their first matrix symmetric.
This will simplify the computation in Chapter~\ref{InvJP} immensely.

We say that the invertible one-sided Jones pairs $(A,B)$ and $(C,B)$ are 
\SL{odd gauge equivalent}\  if
$A=DC\inv{D}$, for some invertible diagonal matrix $D$.
In the following, we examine the odd gauge equivalence of 
one-sided invertible Jones pairs.
\begin{lemma}
\label{lem_Odd_Gauge}
Let $A,C,M$ be Schur-invertible matrices.
If $X_A \D_M = \D_M X_C$, then there exists invertible diagonal matrix $D$ such that
\begin{equation*}
\Sinv{C} \schur A = D J \inv{D}.
\end{equation*}
\end{lemma}
{\sl Proof.}
Since $\D_J = X_I$ is the identity endomorphism of $\mat{n}{\com}$, we have
\begin{equation*}
\DXD{J}{A}{M} = \XDX{I}{M}{C}.
\end{equation*}
Applying the Exchange Lemma, we obtain
\begin{eqnarray*}
\DXD{A}{J}{M^T} &=& \XDX{I}{C}{M}\\
&=&\D_C X_M
\end{eqnarray*}
which gives
\begin{equation*}
\D_{\Sinv{C} \schur A} X_J = X_M \D_{\SinvT{M}}.
\end{equation*}
Evaluating the left-hand side at $E_{1j}=e_1 e_j^T$, we get
\begin{eqnarray*}
(\Sinv{C} \schur A) \schur (Je_1 e_j^T) &=& 
((\Sinv{C} \schur A)e_j \schur \mathbf{1})e_j^T\\
&=& (\Sinv{C} \schur A) e_j e_j^T.
\end{eqnarray*}
The right-hand side evaluated at $E_{1j}$ equals $\inv{(M_{j,1})} Me_1e_j^T$.
Therefore the $ij$-entry of $\Sinv{C} \schur A$ equals $\inv{(M_{j,1})}M_{i,1}$.
So if $D$ is the diagonal matrix with $D_{i,i} = M_{i,1}$ then
\begin{equation*}
\Sinv{C} \schur A = D J \inv{D}.
\end{equation*} 
\eop

\begin{corollary}
\label{cor_Odd_Gauge}
If both $(A,B)$ and $(C,B)$ are invertible one-sided Jones pairs, 
then $A=D C \inv{D}$ for some invertible diagonal matrix $D$.
\end{corollary}
{\sl Proof.}
By Lemma~\ref{lem_XDX_A_B}, we have
\begin{equation*}
\DXD{\Sinv{A}}{A}{A} = \XDX{B}{B^T}{\inv{B}} = \DXD{\Sinv{C}}{C}{C}
\end{equation*}
which gives
\begin{equation*}
X_A \D_{A\schur \Sinv{C}} =  \D_{A\schur \Sinv{C}} X_C.
\end{equation*}
Applying Lemma~\ref{lem_Odd_Gauge} with $M= A \schur \Sinv{C}$, we get
$\Sinv{C} \schur A = DJ \inv{D}$ for some invertible diagonal matrix $D$.
Now 
\begin{equation*}
A= C \schur (DJ \inv{D}) = D (C \schur J) \inv{D} = DC\inv{D}.
\end{equation*} 
\eop

Combining this with Lemma~\ref{lem_1JP_Eq}~(c), we see that $(A,B)$ and $(C,B)$ 
are invertible one-sided Jones pairs
if and only if $A =DC\inv{D}$ for some invertible diagonal matrix $D$.

By Lemma~\ref{lem_1JP_Eq}~(a), 
if $(A,B)$ is an invertible one-sided Jones pair, then so is $(A^T,B)$.
Therefore there exists an invertible diagonal matrix $D$ that satisfies 
$A = DA^T\inv{D}$.
Since the diagonal entries of $D$ are non-zero complex numbers, there exists
diagonal matrix $D_1$ satisfying $D_1^2=D$.
Then the matrix
\begin{equation*}
\inv{D_1}AD_1 = D_1 A^T \inv{D_1} = (\inv{D_1}A D_1)^T
\end{equation*}
is symmetric and $(\inv{D_1}A D_1,B)$ is odd gauge equivalent to $(A,B)$.
So we get the following result,
which generalizes Proposition~7~(ii) from Jaeger \cite{J_4wt}.

\begin{lemma}
\label{lem_Odd_Sym}
Let $(A,B)$ be an invertible one-sided Jones pair. 
Then there exists a symmetric matrix $A'$ such that 
the invertible one-sided Jones pair $(A',B)$
is odd gauge equivalent to $(A,B)$.
\eop
\end{lemma} 

\section{Even Gauge Equivalence}
\label{Even_Gauge}

We say that the invertible one-sided Jones pairs $(A,B)$ and $(A,C)$ are 
\SL{even gauge equivalent}\  if
$C=BP$, for some permutation matrix $P$.
Now we extend Jaeger's result on even gauge equivalence to invertible
one-sided Jones pairs.
\begin{lemma}
\label{lem_Even_Gauge}
Let $F,G$, and $M$ be invertible matrices.
If $\D_F X_M = X_M \D_G$ then $G\inv{F}$ is a permutation matrix.
\end{lemma}
{\sl Proof.}
Multiplying both sides by $\D_J=X_I$, we get
\begin{equation*}
\DXD{F}{M}{J}=\XDX{M}{G}{I}
\end{equation*}
with the Exchange Lemma,
\begin{equation*}
\DXD{M}{F}{J}=\XDX{M}{I}{G}
\end{equation*}
and
\begin{equation*}
\D_M = \XDX{M}{I}{G\inv{F}}.
\end{equation*}
Evaluating both sides at $\Eij$ yields
\begin{equation*}
M_{i,j} e_ie_j^T = M (Ie_j \schur G\inv{F} e_i)e_j^T.
\end{equation*}
Then for all $\range{i}{1}{n}$,
\begin{equation*}
M_{i,j} e_i = (G \inv{F})_{j,i} M e_j.
\end{equation*}
But $M$ is invertible. So for each $j$, there exists a unique $r$ such that 
$Me_j$ is a scalar multiple of $e_r$.
That is, $(G\inv{F})_{j,r}$ is the only non-zero entry in the 
$j$-th row of $G\inv{F}$.
Since $M_{r,j} = (G\inv{F})_{j,r} M_{r,j} \neq 0$, 
we conclude that $(G\inv{F})_{j,r}=1$.
Therefore $G\inv{F}$ is a permutation matrix.  
\eop

\begin{corollary}
\label{cor_Even_Gauge}
If $(A,B)$ and $(A,C)$ are invertible one-sided Jones pairs,
then $C=BP$ for some permutation matrix $P$.
\end{corollary}
{\sl Proof.}
By Lemma~\ref{lem_XDX_A_B}, we have
\begin{equation*}
\XDX{B}{B^T}{\inv{B}} = \DXD{\Sinv{A}}{A}{A} = \XDX{C}{C^T}{\inv{C}}
\end{equation*}
which gives
\begin{equation*}
\D_{C^T}X_{\inv{C}B} = X_{\inv{C}B} \D_{B^T}.
\end{equation*}
Applying Lemma~\ref{lem_Even_Gauge} with $F=C^T$, $G=B^T$ and $M=\inv{C}B$,
there exists a permutation matrix $P$ such that $B^T \invT{C}=P$ which leads to
$C=BP$.  
\eop

Together with Lemma~\ref{lem_1JP_Eq}~(d), $(A,B)$ and $(A,C)$ are invertible
one-sided Jones pairs if and only if $C=BP$ for some permutation matrix $P$.

If $(A,B)$ is an invertible Jones pair, then both $(A,B)$ and $(A,B^T)$ are 
invertible one-sided Jones pairs.
So there exists permutation matrix $P$ such that $B^T=BP$.
Now
\begin{equation*}
B = (B^T)^T=(BP)^T=P^TB^T=P^TBP
\end{equation*}
and hence $B$ and $P$ commute.
We focus on the case where $P$ has order $2r-1$.
Let $Q=P^r$.  Note that $Q^T=P^{r-1}$.
Then we get
\begin{equation*}
(BQ)^T=Q^TB^T=P^{r-1}BP=BQ.
\end{equation*}
Therefore $(A,B)$ is even-gauge equivalent to $(A,BQ)$, where $BQ$ is symmetric.
This is a proof of Proposition~10~(ii) in \cite{J_4wt}.
In summary, if $(A,B)$ is an invertible one-sided Jones pair with
$P=\inv{B}B^T$ having odd order, then $(A,B)$ is gauge equivalent to 
some invertible one-sided Jones pair whose matrices are symmetric.


\chapter{Spin Models}
\label{SpinModels}

In Section~6 of \cite{J_Class}, 
Jaeger proposed to study the properties
of association schemes that contain spin models.
We use this chapter to survey some classical results in this area.
We first examine spin models from the point of view of Jones pairs.
In Section~\ref{Mod_Inv}, 
we present the derivation of the modular invariance equation given by
Jaeger, Matsumoto, and Nomura in \cite{JMN}.
In Section~\ref{Hyper_Dual}, we provide a new and shorter proof
of Curtin and Nomura's result, which states that the Nomura algebra of a
spin model is strongly hyper-self-dual \cite{CN_hyper}.
Section~\ref{Two_Class} contains a shorter proof of Jaeger's 
characterization of two-class association schemes that contain spin models.
This was the first connection between spin models and association schemes 
discovered \cite{J_SRG}.
In the last section, we give a new proof using Jones pairs
of Jaeger and Nomura's result on the symmetric and non-symmetric
Hadamard spin models \cite{JN}.

\section{The Jones Pair $(d^{-1}W,\Sinv{W})$}
\label{Spin_JP}

As mentioned in Section~\ref{Intro_Spin}, Kawagoe, Munemasa and Watatani
defined a spin model with loop variable $d=\pm \sqrt{n}$ to be an $\nbyn$
matrix $W$ that satisfies
\begin{enumerate}[(I)]
\item
There exists some non-zero scalar $a$ such that 
\begin{equation*}
W \schur I = aI, \quad \text{and} \quad WJ=W^T J =d a^{-1}J.
\end{equation*}
\item
$W$ is a type-II matrix.
\item
For all $\range{i,j,k}{1}{n}$,
\begin{equation}
\label{eqn_Type_III}
\sum_{x=1}^n  \frac{W_{k,x}W_{x,i}}{W_{j,x}} = d\frac{W_{k,i}}{W_{j,i}W_{j,k}}.
\end{equation}
\end{enumerate}

An interesting example is the Higman-Sims model discovered by Jaeger
\cite{J_SRG}.
The Higman-Sims graph is a strongly regular graph with parameters $(100,22,0,6)$,
defined from the unique $3-(22,6,1)$ design.
Let $A_1$ and $A_2$ be adjacency matrices of the Higman-Sims graph and its complement
respectively.  If $t$ satisfies $t^2+t^{-2}=-3$, then
\begin{equation*}
W = (5t-3)I+tA_1+t^{-1}A_2
\end{equation*}
is a spin model.
The Nomura algebra $\Nom{W}$ equals the span of $\{I,A_1,A_2\}$, see \cite{JMN}.

In Section~\ref{1JP}, we see that if $W$ is a spin model with loop
variable $d$ then $(d^{-1}W,\SinvT{W})$ is a one-sided Jones pair.
It turns out that $W$ is a spin model if and only if
$(d^{-1}W,\SinvT{W})$ is an invertible Jones pair.
We will prove this statement below.

\begin{lemma}
\label{lem_W_JP}
Let $W$ be a type-II matrix.  
Then $(d^{-1}W,\SinvT{W})$ is a one-sided Jones pair if and only if
$(d^{-1}W,\Sinv{W})$ is also a one-sided Jones pair.
\end{lemma}
{\sl Proof.}
We apply the Exchange Lemma to 
\begin{equation*}
\XDX{d^{-1}W}{\SinvT{W}}{d^{-1}W} = \DXD{\SinvT{W}}{d^{-1}W}{\SinvT{W}}
\end{equation*}
to get
\begin{equation*}
d^{-1}\XDX{W}{W}{\SinvT{W}} = \DXD{W}{\SinvT{W}}{\Sinv{W}}.
\end{equation*}
Since $\SinvT{W} = n \inv{W}$, taking the inverse of each side gives
\begin{equation*}
d \XDX{W}{\Sinv{W}}{\inv{W}} = \DXD{W}{W}{\Sinv{W}},
\end{equation*}
which equals
\begin{equation*}
d \DXD{\Sinv{W}}{W}{\Sinv{W}} = \XDX{W}{\Sinv{W}}{W}.
\end{equation*} 
Since every step above is reversible, the converse is also true.
\eop

If $W$ is a spin model then it is also a type-II matrix.
So the above lemma implies that $(d^{-1}W, \SinvT{W})$ is an invertible
Jones pair.

\begin{lemma}
\label{lem_JP_Spin}
If $(d^{-1}W, \SinvT{W})$ is an invertible Jones pair, then $W$ is a
spin model with loop variable $d$.
\end{lemma}
{\sl Proof.}
By Lemma~\ref{thm_B_Inv}, the invertibility of $\SinvT{W}$ implies that
$W$ is a type-II matrix and $W \schur I = aI$ where $a = \tr(W)n^{-1}$.
Applying Theorem~\ref{lem_1JP_Col_Sum} to the one-sided Jones pairs
$(d^{-1}W, \SinvT{W})$ and $(d^{-1}W,\Sinv{W})$ yields
\begin{equation*}
\Sinv{W}J = da J \quad \text{and} \quad \SinvT{W} J =daJ,
\end{equation*}
respectively.

It follows from $\inv{W} = n^{-1}\SinvT{W}$ that $W^TJ=WJ=da^{-1}J$.
Lastly, using Equation~(\ref{eqn_1JP_AB}) on the one-sided Jones pair
$(d^{-1}W,\SinvT{W})$ gives the type-III condition of spin models.
\eop

\begin{theorem}
\label{thm_JP_Spin}
Let $W$ be an $\nbyn$ matrix.
Then $W$ is a spin model if and only if $(d^{-1}W,\SinvT{W})$ is an
invertible Jones pair.
\eop
\end{theorem}

In the following, we use the theory developed in the previous chapters
to reproduce two existing results about spin models.

Lemma~\ref{lem_W_JP} together with Theorem~\ref{thm_A_in_N_AB},
we have a proof of a result due to Jaeger, Matsumoto and Nomura,
Proposition~9 in \cite{JMN}.

\begin{theorem}
\label{thm_W_in_NW}
Suppose $W$ is a type-II matrix.
Then $W\in \Nom{W}$ if and only if
$cW$ is a spin model for some non-zero scalar $c$.
\eop
\end{theorem} 
We now prove Proposition~2 in \cite{JN} due to Jaeger and Nomura.
\begin{theorem}
\label{thm_JN}
If $W$ is a spin model then there exist a diagonal matrix $D$ and 
a permutation matrix $P$ such that
\begin{equation*}
\SinvT{W} \schur W = DJ\inv{D},
\end{equation*}
and 
\begin{equation*}
n^{-1} \Sinv{W} W =P.
\end{equation*}
\end{theorem}
{\sl Proof.}
By Theorem~\ref{thm_JP_Spin}, if $W$ is a spin model with loop variable
$d$, then
\begin{equation*}
(d^{-1}W, \Sinv{W}), (d^{-1}W^T, \Sinv{W}),
(d^{-1}W, \SinvT{W}), (d^{-1}W^T, \SinvT{W})
\end{equation*} 
are invertible one-sided Jones pairs.
The first equation follows from applying Corollary~\ref{cor_Odd_Gauge}
to the first two invertible one-sided Jones pairs listed above.
Applying Corollary~\ref{cor_Even_Gauge} to the first and the third pairs
above yields the second equality.
\eop

Jaeger and Nomuar called the order of $P$ the \SL{index of the spin model}.
In the same paper, they also proved that any spin model of index two has to take
the following form
\begin{equation*}
\begin{pmatrix}
A&A&B&-B\\
A&A&-B&B\\
-B^T&B^T&C&C\\
B^T&-B^T&C&C\\
\end{pmatrix}
\end{equation*}
where $A$ and $C$ are symmetric (Proposition~7 in \cite{JN}).
As we will see in Section~\ref{Had_Spin}, the construction of the
non-symmetric Hadamard spin models is very similar to the above form.

\section{Duality and Modular Invariance Equation}
\label{Mod_Inv}

Let $P$ be the matrix of eigenvalues for a formally self-dual Bose-Mesner
algebra of dimension $m+1$. 
Let $d=\pm \sqrt{n}$ and
$D$ be a diagonal matrix with $D_{i,i} = t_i$ for $\range{i}{0}{m}$.
We call the equation 
\begin{equation*}
(PD)^3=t_0 d^{3} I
\end{equation*}
the \SL{modular invariance equation}.
We say that $P$ satisfies the \SL{modular invariance property}\ 
if there exists a diagonal matrix $D$ such that the modular invariance
equation holds.
Bannai, Bannai and Jaeger (\cite{BBJ}) first discovered that 
this property of $P$
is a necessary condition for a formally self-dual
Bose-Mesner algebra to contain a spin model. 
Using this equation, 
they provided a method to exhaustively search for all spin models contained 
in a formally self-dual Bose-Mesner algebra.
In Chapter~\ref{InvJP}, we use this equation to design a search for
four-weight spin models.

Suppose $W$ is a spin model.
Then $(W,\Sinv{W})$ is an invertible Jones pair and 
by Theorem~\ref{thm_Inv_JP}, we know that $\Nom{W}=\Nom{W^T}$.
For all $M \in \Nom{W^T}$, we have by Lemma~(\ref{lem_AB_ABT}) 
\begin{eqnarray*}
\T{W^T}{\SinvT{W}}(M) &=&
\T{W^T}{\Sinv{W}}(M) \inv{W} W^T\\
&=& (\T{\Sinv{W}}{W^T}(M))^T \inv{W} W^T.
\end{eqnarray*}
Applying Lemma~(\ref{lem_AB_ABT}) again, the above becomes
\begin{eqnarray*}
\Theta_{W^T}(M)
&=& (\T{\Sinv{W}}{W}(M)W^T\inv{W})^T \inv{W}W^T \\
&=& \invT{W} W(\T{\Sinv{W}}{W}(M))^T \inv{W}W^T \\
&=& \invT{W} W\Theta_W(M) \inv{W}W^T .
\end{eqnarray*}
Now $\Nom{W^T}$ is commutative and it contains $W$, $\invT{W}$, and $\Theta_W(M)$.
So we have $\Theta_{W^T}(M)=\Theta_W(M)$, for all $M$ in $\Nom{W}$. 
By Theorem~\ref{thm_T_A}, we have $\Theta_W^2(M) = nM^T$.
Hence $\Nom{W}$ is formally self-dual.
By Lemma~\ref{lem_T_AB_Formula}, the duality map $\Theta_W$ is expressed
explicitly as
\begin{equation}
\label{eqn_T_W_Formula}
\Theta_W(M) = \frac{n}{\tr(W)} W \schur ((W^T \schur M) \Sinv{W}),
\end{equation}
for all $M \in \Nom{W}$.

In the following, we present Jaeger, Matsumoto and Nomura's proof that
the modular invariance property is a necessary condition of $\Nom{W}$
for $W$ to be a spin model.

\begin{theorem}
\label{thm_Mod_Inv}
Suppose $W$ is a spin model with loop variable $d=\pm \sqrt{n}$ and
$\{\seq{A}{0}{m}\}$ is the basis of Schur idempotents of $\Nom{W}$.
If $W = \sum_{i=0}^m t_iA_i^T$ then the the diagonal matrix $D$
with $D_{i,i}=t_i$ satisfies the modular invariance equation
\begin{equation*}
(PD)^3=t_0 d^{3} I.
\end{equation*}
\end{theorem}
{\sl Proof.}
Let ${\cal E} = \{ \seq{E}{0}{m}\}$ be
the basis of the principal idempotents of $\Nom{W}$ 
such that $\Theta_W(E_i)=A_i$ and $\Theta_W(A_i) = nE_i^T$.
The matrix of eigenvalues $P$ is the matrix of $\Theta_W$ with respect to ${\cal E}$.
However, $P$ is also the transition matrix from $\A$ to ${\cal E}$.
Therefore the matrix of $\Theta_W$ with respect to $\A$ is $\inv{P}PP=P$.
Since $\Theta_W^2(A_i) = nA_i^T$, we have $P^2=nT$ where $T$ represents the transpose
map with respect to $\A$.

Suppose $W=\sum_{i=0}^m t_iA_i^T$.
Since $(d^{-1}W,\Sinv{W})$ is a one-sided Jones pair,
we have $\Sinv{W}=\Theta_W(d^{-1}W)=d\sum_{i=0}^m t_i E_i$.
Now $D$ is the matrix representing the map $M \rightarrow W^T \schur M$
with respect to $\A$.
Similarly $TDT$ represents the map $M\rightarrow W\schur M$ with respect to $\A$.
Moreover, the matrix $d\inv{P}DP$ represents the map $M \rightarrow M\Sinv{W}$ 
with respect to $\A$.
So Equation~(\ref{eqn_T_W_Formula}) holds if and only if
\begin{equation*}
P= t_0^{-1} (TDT)(d\inv{P}DP)D.
\end{equation*}
Since $P^2=nT$, we have $T\inv{P}=\inv{P}T=n^{-1}P$ and 
\begin{eqnarray*}
I &=& t_0^{-1} d (\inv{P}T)D(T\inv{P})DPD\\
&=& t_0^{-1}d^{-3} PDPDPD.
\end{eqnarray*}
Hence
$(PD)^3=t_0 d^{3} I$.  
\eop

\section{Strongly Hyper-Self-Duality}
\label{Hyper_Dual}

In this section, we extend a result due to Curtin and Nomura,
(Theorem~5.5 in \cite{CN_hyper}).

Suppose $\A = \{\seq{A}{0}{d}\}$ is an association scheme with its Bose-Mesner
algebra denoted by $\BM$.
The \SL{Terwilliger algebra of $\BM$}\  can be defined as
\begin{equation*}
\Ter{\BM} = \{ \D_M, X_M : M\in \BM\}.
\end{equation*}
Now consider the subspace $S_p$ of $\mat{n}{\com}$ spanned by 
$\{E_{ip} : \range{i}{1}{n}\}$.
This space is isomorphic to $\com^n$.
For each endomorphism $Y$ of $\mat{n}{\com}$, 
we use $(Y)_p$ to denote $Y$ restricted to $S_p$.
The \SL{Terwilliger algebra of $\BM$ with respect to $p$}\  is defined as
\begin{equation*}
\Ter{\BM,p} = \{ (\D_M)_p, (X_M)_p: M\in \BM\}.
\end{equation*}

A \SL{hyper-duality of $\Ter{\BM,p}$}\  is an automorphism $\Psi_p$ that
swaps the sets 
\begin{equation*}
\{\D_M : M\in \Nom{W}\}\quad \text{and}\quad \{X_M: M \in \Nom{W}\}.
\end{equation*}
and satisfies $\Psi_p^2((X_M)_p) = (X_M)_p^T$ for all $M \in \BM$.
Furthermore, $\Ter{\BM,p}$ is \SL{strongly hyper-self-dual}\ if there exists a hyper-duality
that can be expressed as a conjugation of some invertible element of $\Ter{\BM,p}$.
That is, there exists some $Y \in \Ter{\BM,p}$ such that 
$\Psi_p(Z) = \inv{Y}ZY$, for all $Z \in \Ter{\BM,p}$.

Theorem~5.5 of \cite{CN_hyper} states that if $W$ is a spin model then $\Ter{\Nom{W},p}$
is strongly hyper-self-dual for all $\range{p}{1}{n}$.
In the following, we extend this result to $\Ter{\Nom{W}}$
by showing that there exists $\Lambda \in \Ter{\Nom{W}}$ such that 
the map $\Psi : Z \mapsto \inv{\Lambda} Z \Lambda$ interchanges 
the sets $\{\D_M:M\in \BM\}$ and $\{X_M: M \in \BM\}$, 
and it satisfies $\Psi^2(X_M)=X_M^T$ for all $M \in \Nom{W}$.

Suppose $W$ is a spin model, so $\Theta_W(W)=d \Sinv{W}$.
By Lemma~\ref{lem_T_A_Transpose},
\begin{equation*}
\Theta_W(W^T) = \Theta_W(W)^T = d \SinvT{W}.
\end{equation*}
Hence we have
\begin{equation*}
\XDX{W^T}{\Sinv{W}}{W}=\DXD{\Sinv{W}}{W}{d\SinvT{W}}.
\end{equation*}
We use $\Lambda$ to denote this operator.

\begin{lemma}
\label{lem_Hyper}
Let $W$ be a spin model.
\begin{enumerate}
\item
If $R_1 \in \N{W}{\Sinv{W}}$ and $S_1 = \T{W}{\Sinv{W}}(R_1)$ then
\begin{equation*}
\inv{\Lambda} X_{R_1} \Lambda = \D_{S_1}.
\end{equation*}
\item
If $R_2 \in \N{\Sinv{W}}{W}$ and $S_2 = \T{\Sinv{W}}{W}(R_2)$ then
\begin{equation*}
\inv{\Lambda} \D_{S_2} \Lambda = X_{R_2}.
\end{equation*}
\end{enumerate}
\end{lemma}
{\sl Proof.}
First we prove (a).
Using $\Lambda = d \DXD{\Sinv{W}}{W}{\SinvT{W}}$, we get
\begin{equation*}
\inv{\Lambda} X_{R_1} \Lambda 
= \DXD{W^T}{\inv{W}}{W} (\XDX{R_1}{\Sinv{W}}{W}) \D_{\SinvT{W}},
\end{equation*}
since $\Theta_W(R_1) = S_1$, 
applying Theorem~\ref{thm_T_AB_XDX} gives
\begin{eqnarray*}
\inv{\Lambda} X_{R_1} \Lambda 
&=& \DXD{W^T}{\inv{W}}{W} (\DXD{\Sinv{W}}{W}{S_1}) \D_{\SinvT{W}}\\
&=& \D_{W^T} \D_{S_1} \D_{\SinvT{W}}\\
&=& \D_{S_1}.
\end{eqnarray*}
So part~(a) holds.

Now we prove (b).
Using $\Lambda = \XDX{W^T}{\Sinv{W}}{W}$, we get
\begin{equation*}
\inv{\Lambda} \D_{S_2} \Lambda 
= \XDX{\inv{W}}{W}{\invT{W}} \D_{S_2} \XDX{W^T}{\Sinv{W}}{W}.
\end{equation*}
Since $\invT{W} = n^{-1} \Sinv{W}$, 
\begin{eqnarray*}
\inv{\Lambda} \D_{S_2} \Lambda 
&=& \frac{1}{n} X_{\inv{W}} (\DXD{W}{\Sinv{W}}{S_2}) \XDX{W^T}{\Sinv{W}}{W}\\
&=& \frac{1}{n} X_{\inv{W}} (\XDX{R_2}{W}{\Sinv{W}})\XDX{W^T}{\Sinv{W}}{W}\\
&=& X_{\inv{W}R_2 W}.
\end{eqnarray*}
Since both $R_2$ and $W$ belong to $\Nom{W}$, they commute and
$\inv{W}R_2W=R_2$.  Therefore
\begin{equation*}
\inv{\Lambda} \D_{S_2} \Lambda = X_{R_2}.
\end{equation*} 
\eop

Define an isomorphism $\Psi$ of $\Ter{\Nom{W}}$ as
\begin{equation*}
\Psi(Y) = \inv{\Lambda} Y \Lambda,
\end{equation*}
for all $Y \in \Ter{\Nom{W}}$.
Note that $\Psi$ is expressed as a conjugation of $\Lambda$ in
$\Ter{\Nom{W}}$.  
Now we show that $\Psi$ acts as a hyper-duality for $\Ter{\Nom{W}}$.
\begin{theorem}
\label{thm_Hyper}
If $W$ is a spin model, then the map $\Psi$ defined above interchanges
the sets 
\begin{equation*}
\{\D_M : M\in \Nom{W}\}\quad \text{and}\quad \{X_M: M \in \Nom{W}\}.
\end{equation*}
Moreover, it satisfies $\Psi^2(X_M) = X_M^T$, for all $M \in \Nom{W}$.
\end{theorem}
{\sl Proof.}
Since $\Nom{W} = \Nom{W^T}$,
Lemma~\ref{lem_Hyper} says that $\Psi$ interchanges the sets 
\begin{equation*}
\{\D_M : M\in \Nom{W}\}\quad \text{and}\quad \{X_M: M \in \Nom{W}\}.
\end{equation*}
Now consider ${\Psi}^2(X_R)$.
It follows from Lemma~\ref{lem_Hyper} that
\begin{equation*}
{\Psi}^2(X_R) = \Psi (\D_{\T{W}{\Sinv{W}}(R)}) = X_{R'},
\end{equation*}
for some $R'$ satisfying $\T{W}{\Sinv{W}}(R) = \T{\Sinv{W}}{W}(R')$.
But 
$\T{\Sinv{W}}{W}(R') = \T{W}{\Sinv{W}}({R'}^T)$
by Lemma~\ref{lem_Sinv_AB}. 
Since $\T{W}{\Sinv{W}}$ is an isomorphism, we have $R = {R'}^T$ and
\begin{equation*}
\Psi^2(X_R) = X_{R^T} = {X_R}^T.
\end{equation*}
\eop

\section{Spin Models in Two-class Association Schemes}
\label{Two_Class}

We give a new proof of one direction of Jaeger's result 
about the triply-regularity of 
two-class association schemes that contain a spin model \cite{J_SRG}.
His result was the first indication that spin models have strong
combinatorial properties.

Suppose $\A = \{ \seq{A}{0}{d} \}$ is an association scheme.
It is \SL{triply-regular}\  if for all $i,j,k,r,s,t \in \{0,1,\ldots,d\}$,
the cardinality of the set
\begin{equation*}
\{w : (A_i)_{w,x}=(A_j)_{w,y}=(A_k)_{w,z}=1 \}
\end{equation*}
depends only on $(i,j,k,r,s,t)$, where $(A_r)_{x,y}=(A_s)_{x,z}=(A_t)_{y,z}=1$.

The following Lemma is a direct translation of Lemma~4 in Munemasa's notes
\cite{Munemasa} into the language of endomorphisms.
\begin{lemma}
\label{lem_Triply}
If for all $\range{i,j,k}{0}{d}$,
\begin{equation*}
\XDX{A_i}{A_j}{A_k} \in
\spn( \DXD{A_r}{A_s}{A_t} : r,s,t = 0,1,\ldots,d)
\end{equation*}
then $\A$ is triply-regular.
\end{lemma}
{\sl Proof.}
For each $i,j,k$, the operator $\XDX{A_i^T}{A_j}{A_k}$ lies in the span
of 
\begin{equation*}
\{\DXD{A_r}{A_s}{A_t} : r,s,t = 0,1,\ldots,d\}=
\{\DXD{A_r}{A_s}{A_t^T} : r,s,t = 0,1,\ldots,d\}.
\end{equation*}
So there exists scalars $\kappa(ijk|rst)$ such that
\begin{equation*}
\XDX{A_i^T}{A_j}{A_k} =
\sum_{r,s,t=0}^d \kappa(ijk|rst) \DXD{A_r}{A_s}{A_t^T}.
\end{equation*}
Consider the left-hand side,
\begin{eqnarray*}
\XDX{A_i^T}{A_j}{A_k}(E_{zy})
&=& A_i^T (A_j e_y \schur A_ke_z) e_y^T\\
&=& \sum_x \left(A_i^T (A_j e_y \schur A_ke_z)\right)_x e_xe_y^T\\
&=& \sum_x \left(\sum_w (A_i)_{w,x}(A_j)_{w,y}(A_k)_{w,z}\right) e_x e_y^T\\
&=& \sum_x |\{w : (A_i)_{w,x}=(A_j)_{w,y}=(A_k)_{w,z}=1 \}| E_{xy}.
\end{eqnarray*}
Consider
\begin{eqnarray*}
\DXD{A_r}{A_s}{A_t^T}(E_{zy})
&=&  (A_t)_{y,z} (A_r e_y \schur A_s e_z) e_y^T\\
&=& \sum_x (A_r)_{x,y} (A_s)_{x,z} (A_t)_{y,z} E_{xy}.
\end{eqnarray*}
So we get
\begin{equation*}
|\{w : (A_i)_{w,x}=(A_j)_{w,y}=(A_k)_{w,z}=1 \}|
= \sum_{r,s,t=0}^d \kappa(ijk|rst) (A_r)_{x,y}(A_s)_{x,z} (A_t)_{y,z}
\end{equation*}
and $\A$ is triply-regular.  
\eop

\begin{theorem}
\label{thm_Triply}
Suppose $\A = \{I,A_1,A_2\}$ is a two-class association scheme with
Bose-Mesner algebra $\BM$.
If there exists a spin model $W=t_0 I + t_1 A_1 + t_2 A_2$ with $t_1 \neq t_2$,
then $\A$ is triply-regular.
\end{theorem}
{\sl Proof.}
It is well known that all two-class association schemes are symmetric, so $W=W^T$.
Let $A_0=I$ and 
\begin{equation*}
S = \spn(\DXD{A_r}{A_s}{A_t} : r,s,t=0,1,2).
\end{equation*}
Since $\{A_0,A_1,A_2\}$ is a basis for $\BM$,
\begin{equation*}
S = \spn( \DXD{F'}{G'}{H'} : F',G',H' \in \BM).
\end{equation*}
If we can show that $\XDX{F}{G}{H} \in S$ for all $F,G,H$ in $\BM$, then 
$\A$ is triply-regular by Lemma~\ref{lem_Triply}. 

Now the two sets $\{I,W,J\}$ and $\{I,\Sinv{W},J\}$ are bases for $\BM$.
So it is sufficient to show that for all 
$F,H \in \{I,W,J\}$ and $G\in \{I,\Sinv{W},J\}$
\begin{equation*}
\XDX{F}{G}{H} \in S.
\end{equation*}
When $H=I$, $F \in \{I,W,J\}$ and $G \in \{I,\Sinv{W},J\}$, we get 
\begin{equation*}
\XDX{F}{G}{I}=X_F\D_G=\DXD{J}{F}{G} \in S.
\end{equation*} 
Secondly, when $H=J$, $F \in \{I,W,J\}$ and 
$G \in \{I, \Sinv{W},J\}$,
we can apply the Exchange Lemma to
\begin{equation*}
\XDX{F}{J}{G} = X_{FG} = \DXD{J}{FG}{J},
\end{equation*} 
and get
\begin{equation*}
\XDX{F}{G}{J} = \DXD{FG}{J}{J} 
\end{equation*}
which belongs to $S$ because $FG \in \BM$.

Thirdly, when $H=W$, we enumerate the cases where $G =J$, $I$ or $\Sinv{W}$.
When $H=W$ and $G=J$, we obtain for all $F \{I,W,J\}$ that
\begin{equation*}
\XDX{F}{J}{W} = X_{FW}=\DXD{J}{FW}{J} \in S.
\end{equation*}
When $H=I$, $G=I$ and $F \in \{I,W,J\}$,
we apply the Exchange Lemma to
\begin{equation*}
\XDX{F}{W}{I} = X_F\D_W = \DXD{J}{F}{W},
\end{equation*}
and get
\begin{equation*}
\XDX{F}{I}{W} = \DXD{F}{J}{W} \in S.
\end{equation*}
When $H=W$, $G=\Sinv{W}$ and $F=W$,
we get
\begin{equation*}
\XDX{W}{\Sinv{W}}{W} = d \DXD{\Sinv{W}}{W}{\Sinv{W}} \in S.
\end{equation*}
because $(d^{-1}W, \Sinv{W})$ is a one-sided Jones pair.
When $H=W$, $G=\Sinv{W}$ and $F=I$, $\Theta_W(I)=J$ gives
\begin{equation*}
\XDX{I}{\Sinv{W}}{W}=\DXD{\Sinv{W}}{W}{J} \in S.
\end{equation*}
Lastly when $H=W$, $G=\Sinv{W}$ and $F=J$, 
\begin{equation*}
\Theta_W(J) = \Theta_{W^T}(\Theta_W(I))=nI
\end{equation*}
yields
\begin{equation*}
\XDX{J}{\Sinv{W}}{W}=n \DXD{\Sinv{W}}{W}{I}\in S.
\end{equation*}
By Lemma~\ref{lem_Triply}, we conclude that $\A$ is triply-regular.  
\eop

Suppose $\A=\{I,A_1,A_2\}$ is a triply-regular two-class association scheme.
Let $G$ be the strongly regular graph whose adjacency matrix is $A_1$.
The fact that $\A$ is triply-regular implies that 
the neighborhoods of any vertex in $G$ and its complements $\overline{G}$ 
induce strongly regular graphs, see \cite{J_Triply}.
In this case, we say that both $G$ and $\overline{G}$ are 
\SL{locally strongly regular}.
Now we have proved one direction of the following result due to Jaeger.
For the proof of the converse, please see \cite{J_SRG}.
\begin{theorem}
\label{thm_LSRG}
If $G$ is a strongly regular graph and $A_1$ is its adjacency matrix, then
there exist $t_0,t_1,t_2$ with $t_1 \neq t_2$ such that 
$W=t_0 I +  t_1 A_1 + t_2 (J-A_1-I)$ is a spin model
if and only if both $G$ and $\overline{G}$ are locally strongly regular.
\eop
\end{theorem} 
\section{Symmetric and Non-Symmetric Hadamard Spin Models}
\label{Had_Spin}

In \cite{JN}, Jaeger and Nomura constructed two $4n \times 4n$ spin
models from an $\nbyn$ Hadamard matrix.  They are called the
symmetric and non-symmetric Hadamard spin model.
They are one of the three known infinite families of spin models
that do not result from tensor product.
The other two families are the Potts model and the spin models that
come from finite Abelian groups, see Section~\ref{TypeII_Families}.
We present a new proof of Jaeger and Nomura's result here.

\begin{lemma}
\label{lem_SymH}
Let $A,B,C \in \mat{n}{\com}$ and let W be the following $4n \times 4n$
matrix with $\epsilon = \pm 1$
\begin{equation*}
W = 
\begin{pmatrix}
A & A & B & -B\\
A & A & -B & B\\
\epsilon B^T & -\epsilon B^T & C & C\\
-\epsilon B^T & \epsilon B^T & C & C\\
\end{pmatrix}.
\end{equation*}
Then $W$ is a spin model with loop variable $2d$ with $d=\pm \sqrt{n}$
if and only if the following conditions hold:
\begin{enumerate}
\item
\label{symH_a}
$B$ is a type II matrix; 
\item
\label{symH_b}
$A$ and $C$ are symmetric spin models with loop variable $d$;
\item
\label{symH_c}
$\XDX{C}{\SinvT{B}}{B^T}=d\DXD{\SinvT{B}}{B^T}{\Sinv{A}}$
\item
\label{symH_d}
$\XDX{C}{\SinvT{B}}{\SinvT{B}}=\epsilon d\DXD{B^T}{B^T}{\Sinv{A}}$
\end{enumerate}
\end{lemma}
{\sl Proof.}
The $4n \times 4n$ matrix $W$ is a spin model if and only if
$(d^{-1}W,\Sinv{W})$ is an invertible one-sided Jones pairs.
By construction, $W$ is type II if and only if $A,B$ and $C$ are type II.
The type-III condition is equivalent to 
\begin{equation*}
W \ \Evector{W}{h}{\Sinv{W}}{k} = \frac{2d}{W_{h,k}} \ \Evector{W}{h}{\Sinv{W}}{k},
\quad \text{for } \range{h,k}{1}{4n}.
\end{equation*}
The eigenvector constructed from column $i$ and $j$ of $W$ with $\range{i,j}{1}{n}$ is
\begin{equation*}
\begin{pmatrix}
\Evector{A}{i}{\Sinv{A}}{j} \\\Evector{A}{i}{\Sinv{A}}{j} \\
\Evector{B^T}{i}{\SinvT{B}}{j} \\ \Evector{B^T}{i}{\SinvT{B}}{j}
\end{pmatrix},
\end{equation*}
and the corresponding eigenvalue of $W$ is ${2d}{(A_{i,j})}^{-1}$.
Similarly, for $\range{i,j}{1}{n}$, the following lists all other eigenvectors
for $W$:
\begin{eqnarray*}
&
\begin{pmatrix}
\Evector{A}{i}{\Sinv{A}}{j} \\\Evector{A}{i}{\Sinv{A}}{j} \\
-\Evector{B^T}{i}{\SinvT{B}}{j} \\ -\Evector{B^T}{i}{\SinvT{B}}{j}
\end{pmatrix},
\begin{pmatrix}
\Evector{B}{i}{\Sinv{B}}{j} \\ \Evector{B}{i}{\Sinv{B}}{j} \\
\Evector{C}{i}{\Sinv{C}}{j} \\\Evector{C}{i}{\Sinv{C}}{j} \\
\end{pmatrix},
\begin{pmatrix}
-\Evector{B}{i}{\Sinv{B}}{j} \\ -\Evector{B}{i}{\Sinv{B}}{j} \\
\Evector{C}{i}{\Sinv{C}}{j} \\\Evector{C}{i}{\Sinv{C}}{j} \\
\end{pmatrix},
& \\
&
\begin{pmatrix}
\Evector{A}{i}{\Sinv{B}}{j} \\ -\Evector{A}{i}{B}{j} \\
\epsilon\ \Evector{B^T}{i}{\Sinv{C}}{j} \\ -\epsilon\ \Evector{B^T}{i}{\Sinv{C}}{j}
\end{pmatrix},
\begin{pmatrix}
\Evector{A}{i}{\Sinv{B}}{j} \\ -\Evector{A}{i}{B}{j} \\
-\epsilon\ \Evector{B^T}{i}{\Sinv{C}}{j} \\ \epsilon\ \Evector{B^T}{i}{\Sinv{C}}{j}
\end{pmatrix},
\begin{pmatrix}
\Evector{B}{i}{\Sinv{A}}{j} \\ -\Evector{B}{i}{\Sinv{A}}{j} \\
\epsilon\ \Evector{C}{i}{\SinvT{B}}{j} \\
-\epsilon\ \Evector{C}{i}{\SinvT{B}}{j}
\end{pmatrix},
\begin{pmatrix}
\Evector{B}{i}{\Sinv{A}}{j} \\ -\Evector{B}{i}{\Sinv{A}}{j} \\
-\epsilon\ \Evector{C}{i}{\SinvT{B}}{j} \\
\epsilon\ \Evector{C}{i}{\SinvT{B}}{j}
\end{pmatrix},
&
\end{eqnarray*}
with eigenvalues
\begin{equation*}
\frac{2d}{A_{i,j}},\
\frac{2d}{C_{i,j}},\ \frac{2d}{C_{i,j}},\
\frac{2d}{B_{i,j}},\ -\frac{2d}{B_{i,j}},\
\frac{2\epsilon d}{B_{j,i}},\ -\frac{2\epsilon d}{B_{j,i}},
\end{equation*}
respectively.  

Multiplying the first eigenvector by $W$ gives
\begin{eqnarray*}
\begin{pmatrix}
A & A & B & -B\\
A & A & -B & B\\
\epsilon B^T & -\epsilon B^T & C & C\\
-\epsilon B^T & \epsilon B^T & C & C\\
\end{pmatrix}
\begin{pmatrix}
\Evector{A}{i}{\Sinv{A}}{j} \\\Evector{A}{i}{\Sinv{A}}{j} \\
\Evector{B^T}{i}{\SinvT{B}}{j} \\ \Evector{B^T}{i}{\SinvT{B}}{j}
\end{pmatrix}
&=& 2
\begin{pmatrix}
A \ ( \Evector{A}{i}{\Sinv{A}}{j}) \\ A \ ( \Evector{A}{i}{\Sinv{A}}{j}) \\
C \ (\Evector{B^T}{i}{\SinvT{B}}{j}) \\ C \ (\Evector{B^T}{i}{\SinvT{B}}{j})
\end{pmatrix}\\
&=&\frac{2d}{A_{i,j}}
\begin{pmatrix}
\Evector{A}{i}{\Sinv{A}}{j} \\\Evector{A}{i}{\Sinv{A}}{j} \\
\Evector{B^T}{i}{\SinvT{B}}{j} \\ \Evector{B^T}{i}{\SinvT{B}}{j}
\end{pmatrix}
\end{eqnarray*}
if and only if 
\begin{equation*}
\XDX{A}{\Sinv{A}}{A} = d\DXD{\Sinv{A}}{A}{\Sinv{A}}
\end{equation*}
and 
\begin{equation*}
\XDX{C}{\SinvT{B}}{B^T} = d\DXD{\SinvT{B}}{B^T}{\Sinv{A}}.
\end{equation*}
Repeating the computation on other eigenvectors,
the type-III condition on $W$ is equivalent to the following set of relations
\begin{eqnarray}
\XDX{A}{\Sinv{A}}{A} &=& d \DXD{\Sinv{A}}{A}{\Sinv{A}} 
\label{symH_1} \\
\XDX{C}{\SinvT{B}}{B^T} &=& d \DXD{\SinvT{B}}{B^T}{\Sinv{A}} 
\label{symH_2} \\
\XDX{A}{\Sinv{B}}{B} &=& d \DXD{\Sinv{B}}{B}{\Sinv{C}} 
\label{symH_3} \\
\XDX{C}{\Sinv{C}}{C} &=& d \DXD{\Sinv{C}}{C}{\Sinv{C}}
\label{symH_4} \\
\XDX{B}{\Sinv{C}}{B^T} &=& \epsilon d \DXD{\Sinv{B}}{A}{\Sinv{B}} 
\label{symH_5} \\
\XDX{B^T}{\Sinv{B}}{A} &=& d \DXD{\Sinv{C}}{B^T}{\Sinv{B}} 
\label{symH_6} \\
\XDX{B}{\SinvT{B}}{C} &=& d \DXD{\Sinv{A}}{B}{\SinvT{B}} 
\label{symH_7} \\
\XDX{B^T}{\Sinv{A}}{B} &=& \epsilon d \DXD{\SinvT{B}}{C}{\SinvT{B}}. 
\label{symH_8}
\end{eqnarray}

Now (\ref{symH_1}) and (\ref{symH_4}) hold if and only if $A$ and $C$ are spin models
with loop variable $d$.
If we take the transpose of each side of (\ref{symH_6}) and compare it
with (\ref{symH_3}), then we see that 
these two equations hold if and only if $A$ is symmetric.
Similarly,  
(\ref{symH_7}) and (\ref{symH_2}) hold simultaneously
if and only if $C$ is symmetric. 

It remains to show the equivalence of (\ref{symH_2}) and (\ref{symH_3}),
and the equivalence of (\ref{symH_5}), (\ref{symH_8}) and (\ref{symH_d}).

We can rewrite (\ref{symH_3}) as
\begin{equation*}
\XDX{B}{\Sinv{C}}{\inv{B}} = \frac{1}{d} \DXD{B}{A}{\Sinv{B}},
\end{equation*}
by taking the inverse of each side, we get
\begin{equation*}
\XDX{B}{C}{\inv{B}} = d \DXD{B}{\inv{A}}{\Sinv{B}}. 
\end{equation*}
Taking the transpose of each side gives
\begin{equation*}
\XDX{\invT{B}}{C}{B^T} = d \DXD{\Sinv{B}}{\invT{A}}{B}.
\end{equation*}
Now $n \invT{M}= \Sinv{M}$ for any type II matrix $M$.
After applying the Exchange Lemma, we get
\begin{equation*}
\XDX{\invT{B}}{B^T}{C} = d \DXD{\Sinv{A}}{\invT{B}}{B^T}
\end{equation*}
and whence (\ref{symH_c}), which is identical to (\ref{symH_2}), holds.

By taking the inverse of each side of (\ref{symH_5}), we have
\begin{equation*}
\XDX{\invT{B}}{C}{\inv{B}} = \frac{\epsilon}{d} \DXD{B}{\inv{A}}{B}.
\end{equation*}
Now apply the exchange lemma to get
\begin{equation*}
\XDX{\invT{B}}{\inv{B}}{C} = \frac{\epsilon}{d} \DXD{\inv{A}}{B}{B^T}.
\end{equation*}
By taking the transpose of each side and using the symmetry of $A$ and $C$,
we have
\begin{equation*}
\XDX{C}{\SinvT{B}}{\SinvT{B}} = \epsilon d\DXD{B^T}{B^T}{\Sinv{A}}
\end{equation*}
which equals to (\ref{symH_d}) and can be easily rewritten as
(\ref{symH_8}).  
\eop

The following construction gives the symmetric Hadamard spin models described
in \cite{JN} and \cite{N_Twist} when $\epsilon = 1$. 
When $\epsilon=-1$, it gives the non-symmetric Hadamard spin models.
It is an easy consequence of Lemma~\ref{lem_SymH}.
\begin{corollary}
\label{cor_SymH}
Let $H$ be an $\nbyn$ Hadamard matrix.
Let $A=-u^3I+u^{-1}(J-I)$, with $-u^2-u^{-2}=d$, be the Potts model.
For $\epsilon = \pm 1$ and $\omega$ such that $\omega^4=\epsilon$, 
the $4n \times 4n$ matrix 
\begin{equation*}
W=
\begin{pmatrix}
A & A & \omega H & -\omega H \\
A & A & -\omega H & \omega H \\
\epsilon \omega H^T & -\epsilon \omega H^T & A & A\\
-\epsilon \omega H^T & \epsilon \omega H^T & A & A\\
\end{pmatrix}
\end{equation*}
is a spin model.
\end{corollary}
{\sl Proof.}
We use the construction in Lemma~\ref{lem_SymH} with $B=\omega H$ and $C=A$.
The matrix $B$ is  type II, $\Nom{B^T}=\Nom{H^T}$ is
a Bose Mesner algebra and therefore contains any linear combinations of $I$ and $J$.  
Since
\begin{eqnarray*}
\Theta_{B^T}(A) &=&
\Theta_{H^T}(-u^3 I + u^{-1} (J-I)) \\
&=& -u^3 J + u^{-1} n I - u^{-1} J\\
&=& d (-u^{-3} I + u(J-I)) \\
&=& d\Sinv{A},
\end{eqnarray*}
Condition~(\ref{symH_c}) holds.
Moreover, since $B^T = \omega^2 \SinvT{B}$, Condition~(\ref{symH_d}) is
equivalent to Condition~(\ref{symH_c}).
As a result, Lemma~\ref{lem_SymH} implies that $W$ is a spin model.  
\eop


\def\v#1#2#3{\mathbf{v}^{#1}_{#2#3}}   
\def\I{{\cal I}}                       
\def\J{{\cal J}}                       
\def\1{{\mathbf{1}}}                 
\def\2{{\mathbf{2}}}                 
\def\3{{\mathbf{3}}}                 
\def\4{{\mathbf{4}}}                 
\def\y#1#2{\mathbf{Y}^{#1,#2}_{i,j}}   
\def\Evone{
\begin{pmatrix}
\Evector{A}{i}{\Sinv{A}}{j}\\ \Evector{A}{i}{\Sinv{A}}{j}\\ 
\Evector{\SinvT{B}}{i}{B^T}{j} \\ \Evector{\SinvT{B}}{i}{B^T}{j} \\
\end{pmatrix}
}    
\def\Evtwo{
\begin{pmatrix}
-\Evector{A}{i}{\Sinv{A}}{j}\\ -\Evector{A}{i}{\Sinv{A}}{j}\\ 
\Evector{\SinvT{B}}{i}{B^T}{j} \\ \Evector{\SinvT{B}}{i}{B^T}{j} \\
\end{pmatrix}
}    
\def\Evthree{
\begin{pmatrix}
\Evector{\Sinv{B}}{i}{B}{j} \\ \Evector{\Sinv{B}}{i}{B}{j} \\
\Evector{A}{i}{\Sinv{A}}{j}\\ \Evector{A}{i}{\Sinv{A}}{j}\\ 
\end{pmatrix}
}    
\def\Evfour{
\begin{pmatrix}
\Evector{\Sinv{B}}{i}{B}{j} \\ \Evector{\Sinv{B}}{i}{B}{j} \\
-\Evector{A}{i}{\Sinv{A}}{j}\\ -\Evector{A}{i}{\Sinv{A}}{j}\\ 
\end{pmatrix}
}    
\def\Evfive{
\begin{pmatrix}
d \ \Evector{A}{i}{B}{j} \\ -d \ \Evector{A}{i}{B}{j} \\
d^{-1} \ \Evector{\SinvT{B}}{i}{\Sinv{A}}{j} \\ -d^{-1} \ \Evector{\SinvT{B}}{i}{\Sinv{A}}{j} \\
\end{pmatrix}
}    
\def\Evsix{
\begin{pmatrix}
d \ \Evector{A}{i}{B}{j} \\ -d \ \Evector{A}{i}{B}{j} \\
-d^{-1} \ \Evector{\SinvT{B}}{i}{\Sinv{A}}{j} \\ d^{-1} \ \Evector{\SinvT{B}}{i}{\Sinv{A}}{j} \\
\end{pmatrix}
}    
\def\Evseven{
\begin{pmatrix}
d^{-1} \ \Evector{\Sinv{B}}{i}{\Sinv{A}}{j} \\ -d^{-1} \ \Evector{\Sinv{B}}{i}{\Sinv{A}}{j} \\
d \ \Evector{A}{i}{B^T}{j} \\ -d \ \Evector{A}{i}{B^T}{j} \\
\end{pmatrix}
}    
\def\Eveight{
\begin{pmatrix}
d^{-1} \ \Evector{\Sinv{B}}{i}{\Sinv{A}}{j} \\ -d^{-1} \ \Evector{\Sinv{B}}{i}{\Sinv{A}}{j} \\
-d \ \Evector{A}{i}{B^T}{j} \\ d \ \Evector{A}{i}{B^T}{j} \\
\end{pmatrix}
}    

\chapter{Association Schemes}
\label{InvJP}

We present our main results in this thesis.
In Sections~\ref{N_W} and \ref{Ext},
we construct an $2n \times 2n$ type-II matrix $W$ and 
a pair of $4n \times 4n$ symmetric spin models, $V$ and $V'$,
from each $\nbyn$ invertible Jones pair $(A,B)$. 
We exhibit the intricate relations among the Nomura algebras of these
matrices in Sections~\ref{JP_Quotient} to \ref{JP_Subschemes}.
In Section~\ref{JP_MI},  
we design a strategy that allows us to find invertible
Jones pairs, or equivalently four-weight spin models, up to odd-gauge
equivalence.

The constructions of the type-II matrix and the spin models provide
three new Bose-Mesner algebras attached to a four-weight spin model.
In particular, we get a formally dual pair of Bose-Mesner algebras from
$W$ and a formally self-dual Bose-Mesner algebra from $V$.
So our constructions extend the existing theory of Bose-Mesner algebras
associated with four-weight spin models, which only concerns $\Nom{A}$.
In addition, these algebras form an interesting web of relations.
So we do not have just four Bose-Mesner algebras, 
we have a structured set of four Bose-Mesner algebras.

Our construction of the pair of symmetric spin models generalizes
Nomura's in \cite{N_Twist}.  
It places $A$ and $\Sinv{B}$ as submatrices of $V$ and $V'$.
Since both spin models and four-weight spin models are invertible Jones pairs,
they become submatrices of a pair of symmetric spin models four times their
sizes.  
Hence if we can enumerate all symmetric spin models, then we will
have found all spin models and all four-weight spin models.
This observation leads us to the only known strategy of finding 
four-weight spin models described in Section~\ref{JP_MI},
and answers Bannai's request \cite{Ban1} for such a method.

\section{A Dual Pair of Association Schemes}
\label{N_W}

Suppose $A$ and $B$ are $\nbyn$ matrices and $(A,B)$ is an 
invertible Jones pair.
We will use these matrices to define an $2n \times 2n$ type-II matrix $W$.
Consequently, we get two new Nomura algebras, $\Nom{W}$ and $\Nom{W^T}$,
associated to an invertible Jones pair, or equivalently a four-weight
spin model.
In this section, we show that the dimension of $\Nom{W}$ is 
twice that of $\Nom{A}$.
We also exhibit the basis of Schur idempotents of $\Nom{W}$ and 
the basis of principal idempotents of $\Nom{W^T}$.
Understanding these algebras allows us to see their connections 
to the other two Bose-Mesner algebras associated with the same
invertible Jones pair, in Sections~\ref{JP_Quotient} to
\ref{JP_Dim2}.

Godsil constructed the $2n \times 2n$ matrix $W$ mentioned above
and he hypothesized correctly about the dimension of
$\Nom{W}$.  The author proved his conjecture 
and we present this proof here.
(Subsequently Godsil found a shorter proof, but this assumed $A$ is
symmetric, and gives less information.)

Let $(A,B)$ be an invertible Jones pair.
Recall that by Theorem~\ref{thm_Inv_JP},
\begin{equation*}
\Nom{A}=\Nom{A^T}=\Nom{B}=\Nom{B^T}.
\end{equation*}
Moreover, Theorem~\ref{thm_Dim} tells us that $\N{A}{B}$, $\dualN{A}{B}$
and $\Nom{A}$ have the same dimension.
As discussed in Section~\ref{Odd_Gauge}, there exists an invertible diagonal
matrix $C$ such that $A^T = C^{-2} A C^2$.
We define an $2n \times 2n$ matrix $W$ by
\begin{equation}
\label{eqn_W}
W=
\begin{pmatrix}
A^T & -A^T\\
\SinvT{B}C & \SinvT{B}C\\
\end{pmatrix}.
\end{equation}
It is easy to check that $W$ is type II.
In the following, we show that 
\begin{equation*}
\dim \Nom{W} = 2 \dim \Nom{A},
\end{equation*}
and we find bases for $\Nom{W}$ and $\Nom{W^T}$.

\begin{lemma}
\label{lem_NW1}
We have 
\begin{equation*}
\begin{pmatrix}
J_n & 0 \\ 0 & J_n\\
\end{pmatrix}
\in \Nom{W},
\quad \text{and} \quad
\begin{pmatrix}
I_n & I_n \\ I_n & I_n \\
\end{pmatrix}
\in \Nom{W^T}.
\end{equation*}
\end{lemma}
{\sl Proof.}
The eigenvectors for the matrices in $\Nom{W}$ are
\begin{equation*}
\begin{pmatrix}
\pm \Evector{A^T}{i}{\SinvT{A}}{j} \\
\Evector{\SinvT{B}C}{i}{B^T \inv{C}}{j}\\
\end{pmatrix}
=
\begin{pmatrix}
\pm \Evector{A^T}{i}{\SinvT{A}}{j} \\
\frac{C_{i,i}}{C_{j,j}} \ \Evector{\SinvT{B}}{i}{B^T}{j}\\
\end{pmatrix}.
\end{equation*}
Applying Lemma~\ref{lem_J_NA} to the type-II matrices $A^T$ and
$\SinvT{B}C$, we get
\begin{equation*}
\Theta_{A^T}(J) = \Theta_{\SinvT{B}C}(J) = nI.
\end{equation*}
So for $\range{i,j}{1}{n}$,
\begin{equation*}
\begin{pmatrix}
J & 0 \\ 0 & J\\
\end{pmatrix}
\begin{pmatrix}
\pm \Evector{A^T}{i}{\SinvT{A}}{j} \\
\Evector{\SinvT{B}C}{i}{B^T \inv{C}}{j}\\
\end{pmatrix}
=
\delta_{ij}n
\begin{pmatrix}
\pm \Evector{A^T}{i}{\SinvT{A}}{j}\\
\Evector{\SinvT{B}C}{i}{B^T\inv{C}}{j}\\
\end{pmatrix}.
\end{equation*}
As a result, 
\begin{equation*}
\begin{pmatrix}
J & 0 \\ 0 & J\\
\end{pmatrix}
\in \Nom{W},
\end{equation*}
and its image under $\Theta_W$ is
\begin{equation*}
n
\begin{pmatrix}
I&I\\I&I\\
\end{pmatrix}.
\end{equation*} 
\eop

Now the Schur idempotents of $\Nom{W}$ that sum to 
\begin{equation*}
\begin{pmatrix}
J & 0 \\ 0 & J\\
\end{pmatrix}
\end{equation*}
have the form 
\begin{equation*}
\begin{pmatrix}
C_1 & 0 \\ 0 & C_2 \\
\end{pmatrix},
\end{equation*}
and the rest have the form
\begin{equation*}
\begin{pmatrix}
0 & C_1\\ C_2 & 0 \\
\end{pmatrix}.
\end{equation*}
We examine the two types of Schur idempotents below.

The following two lemmas analyze the Schur idempotents with zero
diagonal blocks.
\begin{lemma}
\label{lem_NW2}
If 
$
\begin{pmatrix}
0 & H_1 \\ H_2^T & 0\\
\end{pmatrix}
\in \Nom{W},
$
then $H_1, H_2 \in \dualN{A}{B}$.
\end{lemma}
{\sl Proof.}
Since
\begin{equation*}
\begin{pmatrix}
\Evector{A^T}{i}{\SinvT{A}}{j}\\
\Evector{\SinvT{B}C}{i}{B^T\inv{C}}{j}\\
\end{pmatrix}
\end{equation*}
is an eigenvector, there exists some matrix $S$ such that
\begin{eqnarray*}
H_1 \ (\Evector{\SinvT{B}C}{i}{B^T \inv{C}}{j}) &=&
S_{i,j} \ (\Evector{A^T}{i}{\SinvT{A}}{j}),\\
H_2^T \ (\Evector{A^T}{i}{\SinvT{A}}{j}) &=&
S_{i,j} \ (\Evector{\SinvT{B}C}{i}{B^T \inv{C}}{j}).
\end{eqnarray*}
This is equivalent to
\begin{eqnarray}
\label{eqn_NW2_1}
\XDX{H_1}{B^T \inv{C}}{\SinvT{B}C} &=& \DXD{\SinvT{A}}{A^T}{S},\\
\label{eqn_NW2_2}
\XDX{H_2^T}{\SinvT{A}}{A^T} &=& \DXD{B^T\inv{C}}{\SinvT{B}C}{S}.
\end{eqnarray}
Since $\SinvT{B}=n\inv{B}$, Equation~(\ref{eqn_NW2_1}) can be rewritten as
\begin{equation*}
\DXD{A^T}{H_1}{B^T\inv{C}} = \XDX{A^T}{S}{n^{-1} \inv{C}B},
\end{equation*}
and applying the Exchange Lemma, we get
\begin{equation*}
\DXD{H_1}{A^T}{\inv{C}B} = \frac{1}{n} \XDX{A^T}{\inv{C}B}{S}.
\end{equation*}
Recall that $\inv{C}$ is a diagonal matrix, so 
$\D_{\inv{C}B}$ equals $\D_B X_{\inv{C}}$ 
and the above equation becomes
\begin{equation*}
\DXD{H_1}{A^T}{B}X_{\inv{C}} = \frac{1}{n} \XDX{A^T}{B}{\inv{C}}X_{S}.
\end{equation*}
Taking the transpose of each side gives
\begin{equation*}
X_{\inv{C}} \DXD{B}{A}{H_1} = \frac{1}{n} X_{S^T} \XDX{\inv{C}}{B}{A}.
\end{equation*}
Therefore
\begin{equation*}
\frac{1}{n}\XDX{CS^T\inv{C}}{B}{A} = \DXD{B}{A}{H_1},
\end{equation*}
and $H_1 = n^{-1} \T{A}{B}(CS^T\inv{C}) \in \dualN{A}{B}$.

Equation~(\ref{eqn_NW2_2}) is equivalent to
\begin{equation*}
\DXD{\SinvT{B}C}{H_2^T}{\SinvT{A}} = \XDX{\SinvT{B}C}{S}{\inv{(A^T)}}.
\end{equation*}
Now replacing $\inv{(A^T)}$ by $n^{-1} \Sinv{A}$ and 
applying the Exchange Lemma gives,
\begin{equation*}
\DXD{H_2^T}{\SinvT{B}C}{\Sinv{A}} = \XDX{\SinvT{B}C}{n^{-1}\Sinv{A}}{S}.
\end{equation*}
Taking the transpose of each side yields
\begin{equation*}
\DXD{\Sinv{A}}{C \Sinv{B}}{H_2^T} = \frac{1}{n} \XDX{S^T}{\Sinv{A}}{C\Sinv{B}}.
\end{equation*}
Again, $C$ is a diagonal matrix, so
$\D_{\Sinv{A}} X_C = X_C \D_{\Sinv{A}}$ and
\begin{equation*}
X_C \DXD{\Sinv{A}}{\Sinv{B}}{H_2^T} = \frac{1}{n} X_{S^T} \XDX{C}{\Sinv{A}}{\Sinv{B}}.
\end{equation*}
Consequently,
\begin{equation*}
\frac{1}{n}\XDX{\inv{C}S^TC}{\Sinv{A}}{\Sinv{B}} = \DXD{\Sinv{A}}{\Sinv{B}}{H_2^T},
\end{equation*}
and
\begin{equation*}
H_2^T = \frac{1}{n} \T{\Sinv{B}}{\Sinv{A}}(\inv{C}S^TC).
\end{equation*}
This is the same as
\begin{eqnarray*}
H_2 &=& \frac{1}{n} \T{\Sinv{A}}{\Sinv{B}}(\inv{C}S^TC)\\
&=& \frac{1}{n}\T{A}{B}(CS\inv{C}),
\end{eqnarray*}
by Lemma~\ref{lem_Sinv_AB}. So $H_2 \in \dualN{A}{B}$.  
\eop

Now if $F \in \N{A}{B}$, then by Theorem~\ref{thm_Dim}~(a),
there exists a matrix $H \in \Nom{A}$ such that $F = H \schur A$.
Taking the transpose of each side, we have 
\begin{eqnarray*}
F^T &=& H^T \schur A^T \\
&=&  H^T \schur (C^{-2} A C^2)\\
&=& C^{-2} (H^T \schur A) C^2.
\end{eqnarray*}
So $C^2 F^T C^{-2} = H^T \schur A$.  Since $H^T \in \Nom{A}$, 
we know that $C^2 F^T C^{-2} \in \N{A}{B}$ by Theorem~\ref{thm_Dim}~(a). 
If $\N{A}{B}$ has dimension $r$, the following lemma gives $r$ 
Schur idempotents of $\Nom{W}$.

\begin{lemma}
\label{lem_NW3}
Let $\{\seq{F}{0}{r-1}\}$ be the basis of the principal idempotents of 
$\N{A}{B}$.  Then for $\range{k}{0}{r-1}$, the matrix 
\begin{equation}
\label{eqn_NW3}
\widehat{F_k}=
\begin{pmatrix}
0 & \T{A}{B}(F_k)\\
\T{A}{B}(C^2 F_k^T C^{-2})^T & 0\\
\end{pmatrix}
\end{equation}
is a Schur idempotent of $\Nom{W}$, and 
\begin{equation*}
\frac{1}{2}
\begin{pmatrix}
CF_k^T \inv{C} & -CF_k^T \inv{C} \\
-CF_k^T \inv{C} & CF_k^T \inv{C} \\
\end{pmatrix} 
\end{equation*}
is a principal idempotent of $\Nom{W^T}$.
\end{lemma}
{\sl Proof.}
Let $S = C F_k^T\inv{C}$.
It follows from the proof of the previous lemma that
the following two equations hold:
\begin{eqnarray}
\label{eqn_NW3_1}
\XDX{\T{A}{B}(F_k)}{B^T\inv{C}}{\SinvT{B}C} &=&
\DXD{\SinvT{A}}{A^T}{nS},\\
\label{eqn_NW3_2}
\XDX{\T{A}{B}(C^2 F_k^T C^{-2})^T}{\SinvT{A}}{A^T} &=& 
 \DXD{B^T \inv{C}}{\SinvT{B}C}{nS}.
\end{eqnarray}
We conclude that $\widehat{F_k} \in \Nom{W}$.

By Equations (\ref{eqn_NW3_1}) and (\ref{eqn_NW3_2}), we have
\begin{equation*}
\widehat{F_k}
\begin{pmatrix}
\Evector{A^T}{i}{\SinvT{A}}{j} \\ \Evector{\SinvT{B}C}{i}{B^T \inv{C}}{j}\\
\end{pmatrix}
= 
n (C F_k^T \inv{C})_{i,j}
\begin{pmatrix}
\Evector{A^T}{i}{\SinvT{A}}{j} \\ \Evector{\SinvT{B}C}{i}{B^T \inv{C}}{j}\\
\end{pmatrix}.
\end{equation*}
Thus the image of $\widehat{F_k}$ under $\Theta_W$ is
\begin{equation*}
n
\begin{pmatrix}
CF_k^T \inv{C} & -CF_k^T \inv{C} \\
-CF_k^T \inv{C} & CF_k^T \inv{C} \\
\end{pmatrix} \in \Nom{W^T}.
\end{equation*}

Consider
\begin{equation*}
\widehat{F_k} \schur \widehat{F_l} = 
\begin{pmatrix}
0 & \T{A}{B}(F_k) \schur \T{A}{B}(F_l)\\
\left( \T{A}{B}(C^2 F_k^T C^{-2}) \schur \T{A}{B}( C^2 F_l^T C^{-2}) \right)^T & 0\\
\end{pmatrix},
\end{equation*}
for any $\range{l,k}{0}{r-1}$. 
By Lemma~\ref{lem_T_AB_Swap}, we have
\begin{equation*}
\widehat{F_k} \schur \widehat{F_l} = 
\begin{pmatrix}
0 & \T{A}{B}(F_k F_l)\\
\T{A}{B}\left((C^2 F_k^T C^{-2})(C^2 F_l^T C^{-2}\right)^T & 0\\
\end{pmatrix}.
\end{equation*}
Since $F_k$ and $F_l$ are the principal idempotents of $\N{A}{B}$, we
have $F_k F_l = \delta_{kl} F_k$ and 
$(C^2 F_k^T C^{-2})(C^2 F_l^T C^{-2}) = \delta_{kl}C^2 F_k^T C^{-2}$.
So $\widehat{F_k} \schur \widehat{F_l} = \delta_{kl} \widehat{F_k}$,
that is, $\seq{\widehat{F}}{0}{r-1}$  are the Schur idempotents of
$\Nom{W}$.
Moreover, for any Schur idempotent $M$ in $\Nom{W}$,
the matrix $\frac{1}{2n}\Theta_W(M)$ is a principal idempotent of $\Nom{W^T}$.
So the result follows.  
\eop

We conclude from Lemmas~\ref{lem_NW2} and \ref{lem_NW3} that the matrices
in (\ref{eqn_NW3}) form $r$ Schur idempotents which span the subspace of
$\Nom{W}$ consisting matrices with $\nbyn$ zero diagonal blocks.
Further, we see in the next corollary that 
$A$ and $B$ are encoded in $\Nom{W}$ and $\Nom{W^T}$. 

\begin{corollary}
\label{cor_NW23}
We have
\begin{eqnarray*}
\begin{pmatrix}
0&B\\B^T&0\\
\end{pmatrix}
&\in & \Nom{W},\\
\begin{pmatrix}
\inv{C} A C & -\inv{C} A C \\
-\inv{C} A C & \inv{C} A C \\
\end{pmatrix} 
&\in &\Nom{W^T}.
\end{eqnarray*}
\end{corollary}
{\sl Proof.}
Since $A \in \N{A}{B}$ and $A=C^2 A^T C^{-2}$,
we know that $\T{A}{B}(A)=B$ and 
$\T{A}{B}(C^2 A^T C^{-2})=\T{A}{B}(A)=B$.
Hence $B \in \dualN{A}{B}$ and it follows from Lemma~\ref{lem_NW3} that
\begin{equation*}
\begin{pmatrix}
0&B\\B^T&0\\
\end{pmatrix}
\in  \Nom{W}
\end{equation*}
Moreover its image under $\Theta_W$ equals
\begin{equation*}
n
\begin{pmatrix}
CA^T \inv{C} & -CA^T \inv{C} \\
-CA^T \inv{C} & CA^T \inv{C} \\
\end{pmatrix} 
=n
\begin{pmatrix}
\inv{C} A C & -\inv{C} A C \\
-\inv{C} A C & \inv{C} A C \\
\end{pmatrix}
\end{equation*} 
belongs to $\Nom{W^T}$.
\eop

\begin{lemma}
\label{lem_NW4}
All matrices in $\Nom{W}$ with $\nbyn$ zero off-diagonal blocks have the form
\begin{equation*}
\begin{pmatrix}
\Theta_A(M) & 0\\
0 & \Theta_B(M)^T\\
\end{pmatrix},
\end{equation*}
for some $M$ in $\Nom{A}$.
\end{lemma}
{\sl Proof.}
Suppose
\begin{equation*}
\begin{pmatrix}
C_1&0\\ 0&C_2\\
\end{pmatrix}
\in \Nom{W}.
\end{equation*}
By examining the eigenvectors
$\Evector{W}{i}{\Sinv{W}}{j}$, we conclude that
$C_1 \in \Nom{A^T}$ and $C_2 \in \Nom{\SinvT{B}C}$.
By Corollary~\ref{cor_NW23}, 
$
\begin{pmatrix}
0&B\\B^T&0\\
\end{pmatrix}
\in \Nom{W}$.
Since $\Nom{W}$ is commutative, we have
\begin{equation*}
\begin{pmatrix}
0&B\\B^T&0\\
\end{pmatrix}
\begin{pmatrix}
C_1&0\\ 0&C_2\\
\end{pmatrix}
= 
\begin{pmatrix}
C_1&0\\ 0&C_2\\
\end{pmatrix}
\begin{pmatrix}
0&B\\B^T&0\\
\end{pmatrix},
\end{equation*}
which leads to $BC_2=C_1B$ and $C_2 =\inv{B} C_1 B$.
Since $C_1 \in \Nom{A^T}$, 
there exists $M \in \Nom{A}$ such that $\Theta_A(M)=C_1$.
By Theorem~\ref{thm_Inv_JP}, 
we have 
\begin{equation*}
C_2= \inv{B} \Theta_A(M) B = \Theta_B(M)^T.  
\end{equation*}
\eop

Since $A$ and $B$ are type II, applying Theorem~\ref{thm_Dim}~(a) to
$A \in \N{A}{B}$ give $A \schur \Nom{A} = \N{A}{B}$.
Therefore the dimensions of $\N{A}{B}$ amd $\Nom{A}$ are both $r$.
We see below that the Schur idempotents of $\Nom{W}$ with zero diagonal
blocks can be expressed in terms of matrices in $\Nom{A}$.
\begin{lemma}
\label{lem_NW5}
Let $\{\seq{E}{0}{r-1}\}$ be the basis of the principal idempotents of $\Nom{A}$.
Then for $\range{k}{0}{r-1}$, the matrix
\begin{equation*}
\tilde{E_k} = 
\begin{pmatrix}
\Theta_{A}(E_k) & 0\\ 0 & \Theta_B(E_k)^T\\
\end{pmatrix}
\end{equation*}
is a Schur idempotent of $\Nom{W}$,
and
\begin{equation*}
\frac{1}{2}
\begin{pmatrix}
E_k^T & E_k^T \\ E_k^T & E_k^T \\
\end{pmatrix}
\end{equation*}
is a principal idempotent of $\Nom{W^T}$.
\end{lemma}
{\sl Proof.}
Since $\Theta_{A^T}(\Theta_A(E_k)) = nE_k^T$ and 
\begin{equation*}
\Theta_{\SinvT{B}}(\Theta_B(E_k)^T) =
\Theta_{\SinvT{B}}(\Theta_{\Sinv{B}}(E_k)) = nE_k^T,
\end{equation*}
we have
\begin{equation*}
\begin{pmatrix}
\Theta_{A}(E_k) & 0\\ 0 & \Theta_B(E_k)^T\\
\end{pmatrix}
\begin{pmatrix}
\pm \Evector{A^T}{i}{\SinvT{A}}{j} \\
\frac{C_{i,i}}{C_{j,j}} \ \Evector{\SinvT{B}}{i}{B^T}{j}\\
\end{pmatrix}
=n (E_k)_{j,i}
\begin{pmatrix}
\pm \Evector{A^T}{i}{\SinvT{A}}{j} \\
\frac{C_{i,i}}{C_{j,j}} \ \Evector{\SinvT{B}}{i}{B^T}{j}\\
\end{pmatrix}.
\end{equation*}
Hence $\tilde{E_k} \in \Nom{W}$
and 
\begin{equation*}
\Theta_W(\tilde{E_k})=
n
\begin{pmatrix}
E_k^T & E_k^T \\ E_k^T & E_k^T \\
\end{pmatrix}.
\end{equation*} 

Consider
\begin{equation*}
\tilde{E_k} \schur \tilde{E_l} = 
\begin{pmatrix}
\Theta_{A}(E_k) \schur \Theta_{A}(E_l) & 0\\
0 & \left(\Theta_B(E_k) \schur \Theta_B(E_l) \right)^T\\
\end{pmatrix},
\end{equation*}
for any $\range{k,l}{0}{r-1}$.
By Lemma~\ref{lem_T_AB_Swap}, we have
\begin{equation*}
\tilde{E_k} \schur \tilde{E_l} = 
\begin{pmatrix}
\Theta_{A}(E_kE_l) & 0\\
0 & \Theta_B(E_kE_l)^T\\
\end{pmatrix}.
\end{equation*}
Since $E_k$ and $E_l$ are the principal idempotents of $\Nom{A}$, 
we have $E_k E_l = \delta_{kl} E_k$.
So $\tilde{E_k} \schur \tilde{E_l} = \delta_{kl} \tilde{E_k}$.
That is, $\seq{\tilde{E}}{0}{r-1}$ are the Schur idempotents of
$\Nom{W}$.
Moreover, for any Schur idempotent $M$ in $\Nom{W}$, the matrix
$\frac{1}{2n}\Theta_W(M)$ is a principal idempotent of $\Nom{W^T}$.
So the result follows.
\eop

Combining the Lemmas~\ref{lem_NW3} and \ref{lem_NW5},
we find the basis of Schur idempotents for $\Nom{W}$
and the basis of principal idempotents for $\Nom{W^T}$.

\begin{theorem}
\label{thm_NW}
Suppose $\dim(\Nom{A})=\dim(\N{A}{B})=r$.
Let $\{\seq{E}{0}{r-1}\}$ be the basis of the principal idempotents of $\Nom{A}$.
Let $\{\seq{F}{0}{r-1}\}$ be the basis of the principal idempotents of $\N{A}{B}$.
Then the set
\begin{equation*}
\left\{
\begin{pmatrix}
\Theta_A(E_i) & 0\\ 0 & \Theta_B(E_i)^T\\
\end{pmatrix}, \quad
\begin{pmatrix}
0 & \T{A}{B}(F_j)\\
\T{A}{B}(C^2 F_j^T C^{-2})^T & 0\\
\end{pmatrix}
: \range{i,j}{0}{r-1} \right\}
\end{equation*}
is the basis of Schur idempotents for $\Nom{W}$.
Further, the set
\begin{equation*}
\left\{
\frac{1}{2}
\begin{pmatrix}
E_i^T & E_i^T \\ E_i^T & E_i^T \\
\end{pmatrix}, \quad
\frac{1}{2}
\begin{pmatrix}
C F_j^T \inv{C} & -C F_j^T \inv{C} \\
-C F_j^T \inv{C} & C F_j^T \inv{C} \\
\end{pmatrix}
: \range{i,j}{0}{r-1} \right\}
\end{equation*}
is the basis of principal idempotents for $\Nom{W^T}$.
Hence 
\begin{equation*}
\dim \Nom{W} = \dim \Nom{W^T} = 2 \dim \Nom{A}.
\end{equation*}
\eop
\end{theorem} 

In general $\Nom{W}$ is not equal to $\Nom{W^T}$.
We now examine the situation where these two algebras coincide.

\begin{lemma}
\label{lem_NW_Fsd}
If $\Nom{W}=\Nom{W^T}$, then there exists non-zero scalar $\alpha$ such
that $\alpha B$ is a spin model.
\end{lemma}
{\sl Proof.}
By Corollary~\ref{cor_NW23}, we have
\begin{equation*}
\widehat{B} = 
\begin{pmatrix}
0 & B \\ B^T & 0\\
\end{pmatrix}
\in \Nom{W^T}.
\end{equation*}
Now 
\begin{equation*}
W^T = 
\begin{pmatrix}
A & C \Sinv{B} \\ 
-A & C \Sinv{B} \\ 
\end{pmatrix},
\end{equation*}
and 
\begin{equation*}
\begin{pmatrix}
\Evector{C\Sinv{B}}{i}{\inv{C}B}{j}\\
\Evector{C\Sinv{B}}{i}{\inv{C}B}{j}\\
\end{pmatrix}
=
\begin{pmatrix}
\Evector{\Sinv{B}}{i}{B}{j}\\
\Evector{\Sinv{B}}{i}{B}{j}\\
\end{pmatrix}
\end{equation*}
is an eigenvector of $\widehat{B}$.
That is, 
$B \ (\Evector{\Sinv{B}}{i}{B}{j}) = \beta \ (\Evector{\Sinv{B}}{i}{B}{j})$,
for some $\beta \in \com$
This is equivalent to $B$ belongs to $\Nom{\Sinv{B}}$,
which is identical to $\Nom{B}$.  
By Theorem~\ref{thm_W_in_NW}, there exists non-zero scalar $\alpha$
such that $\alpha B$ is a spin model.  
\eop
This lemma tells us that if $(A,B)$ is not gauge equivalent to 
$(\alpha^{-1}d^{-1} \Sinv{B},\alpha B)$ for $\alpha B$ a spin model,
then we do have a formally dual pair of Bose-Mesner algebras.

\section{Nomura's Extension}
\label{Ext}

In this section, we build two $4n \times 4n$ symmetric spin models, $V$
and $V'$, from an invertible Jones pair.  We show that they share the
same formally self-dual Nomura algebra.  We obtain an explicit form of the
matrices belonging to this algebra.

In Section~5 of \cite{N_Twist}, Nomura constructed an $4n \times 4n$
symmetric spin model from a four-weight spin model 
whose matrices are symmetric.  
Our construction is a generalization of Nomura's, which is motivated by
the fact that Jones pairs are equivalent to four-weight spin models.
Moreover, our construction delivers the fourth
Bose-Mesner algebra associated with a four-weight spin model.  
It becomes evident in Section~\ref{JP_MI} that this Bose-Mesner algebra
is the ticket to our strategy of finding four-weight spin models.

By Lemma~\ref{lem_Odd_Sym} and the discussion prior to it, given any
invertible Jones pair $(A,B)$, we can construct an invertible Jones
pair $(A',B)$ with $A'$ symmetric which is odd-gauge equivalent to
$(A,B)$.
From this new invertible Jones pair $(A',B)$, 
we construct the $4n \times 4n$ symmetric spin models $V$ and $V'$.
So in the following results, we can focus on only the invertible Jones pairs 
with their first matrix symmetric.

Suppose $(A,B)$ is an invertible one-sided Jones pair with $A$ symmetric.
We define
\begin{equation*}
V = 
\begin{pmatrix}
dA & -dA & \Sinv{B} & \Sinv{B} \\ 
-dA & dA & \Sinv{B} & \Sinv{B} \\ 
\SinvT{B} & \SinvT{B} & dA & -dA \\
\SinvT{B} & \SinvT{B} & -dA & dA \\
\end{pmatrix}.
\end{equation*}
It follows from the construction of $V$ together with 
$A$ and $B$ being type II that $V$ is a type-II matrix.
For easier reading, we separate the columns and the rows into four
groups $\1$, $\2$, $\3$, and $\4$.
For instance, the $(\1,\2)$-block of $V$ equals $-dA$ while its
$(\2,\4)$-block equals $\Sinv{B}$.

For $\alpha, \beta \in \{ \1,\2,\3,\4\}$ and $\range{i,j}{1}{n}$,
we use $\y{\alpha}{\beta}$ to denote the eigenvectors we get from $V$.
For example, 
\begin{equation*}
\y{\3}{\1} = \Evector{V}{2n+i}{\Sinv{V}}{j}.
\end{equation*}

In the following, We show that $V$ is a spin model.
Our plan is to find the explicit block structure of the matrices in
$\Nom{V}$.  After that, we will prove that $V\in\Nom{V}$ and
$\Theta_V(V)=(2d)^{-1}\Sinv{V}$, which imply that $V$ is a spin model with
loop variable $2d$.

In the following, we use $J_k$ and $I_k$ denote the $k \times k$ matrix of all ones 
and the $k \times k$ identity matrix, respectively.
\begin{lemma}
\label{lem_NV1}
Let $(A,B)$ be an invertible Jones pair with $A$ symmetric.
We have
\begin{equation*}
(I_2 \otimes J_{2n}) = 
\begin{pmatrix}
J_n & J_n & 0 & 0\\
J_n & J_n & 0 & 0\\
0 & 0 & J_n & J_n\\
0 & 0 & J_n & J_n\\
\end{pmatrix}
\in \Nom{V},
\end{equation*}
and 
\begin{equation*}
\begin{pmatrix}
I_n & I_n & 0 & 0\\
I_n & I_n & 0 & 0\\
0 & 0 & I_n & I_n\\
0 & 0 & I_n & I_n\\
\end{pmatrix}
\in \Nom{V^T}.
\end{equation*}
\end{lemma}
{\sl Proof.}
Consider
\begin{equation*}
\y{\1}{\1}=\Evone.
\end{equation*}
Since $\Theta_A(J_n)=\Theta_{\SinvT{B}}(J_n) = nI_n$, we have
$(I_2 \otimes J_{2n})\ \y{\1}{\1}=2n \delta_{ij}\ \y{\1}{\1}$
and the $(\1,\1)$-block of $\Theta_V(I_2 \otimes J_{2n})$ equals 
$2n I_n$.
Similarly, when $\alpha,\beta \in \{ \1,\2\}$, we have
\begin{equation*}
(I_2 \otimes J_{2n})\ \y{\alpha}{\beta} = 2n \delta_{ij}\ \y{\alpha}{\beta},
\end{equation*}
for all $\range{i,j}{1}{n}$, 
and the same holds for $\alpha, \beta \in \{\3,\4\}$.

Consider
\begin{equation*}
\y{\1}{\3}=\Evfive.
\end{equation*}
We have $(I_2 \otimes J_{2n})\ \y{\1}{\3} = 0$ and so the
$(\1,\3)$-block of $\Theta_V(I_2 \otimes J_{2n})$ is the $\nbyn$ zero
matrix.
The same holds for 
$\y{\1}{\4}$, $\y{\2}{\3}$, $\y{\2}{\4}$, $\y{\3}{\1}$, $\y{\3}{\2}$,
$\y{\4}{\1}$ and $\y{\4}{\2}$.
Hence $I_2 \otimes J_{2n} \in \Nom{V}$ and 
\begin{equation*}
\Theta_V(I_2 \otimes J_{2n}) = 2n
\begin{pmatrix}
I_n & I_n & 0 & 0\\
I_n & I_n & 0 & 0\\
0 & 0 & I_n & I_n\\
0 & 0 & I_n & I_n\\
\end{pmatrix}
\end{equation*}
and the lemma holds.  
\eop

One consequence is that any Schur idempotent of $\Nom{V}$ has one of the
following two forms:
\begin{equation*}
\begin{pmatrix}
M_1 & N_1 & 0 & 0\\
P_1 & Q_1 & 0 & 0\\
0 & 0 & M_2 & N_2\\
0 & 0 & P_2 & Q_2\\
\end{pmatrix},
\quad
\begin{pmatrix}
0 & 0 & M_1 & N_1\\
0 & 0 & P_1 & Q_1\\
M_2 & N_2 & 0 & 0\\
P_2 & Q_2 & 0 & 0\\
\end{pmatrix}.
\end{equation*}

We need the next two lemmas to anatomize the Schur idempotents with
zero off-diagonal blocks.
\begin{lemma}
\label{lem_NV2_1}
Let $(A,B)$ be invertible Jones pair with $A$ symmetric.  
If $M \in \Nom{A}$, then
\begin{eqnarray*}
\Theta_A(M) &=& \Theta_{\SinvT{B}}(\inv{B}MB), \\
\Theta_{\Sinv{B}}(M) &=& \Theta_A(\inv{B}MB).
\end{eqnarray*}
\end{lemma}
{\sl Proof.}
Let $S = \Theta_A(M)$.
Then 
\begin{equation*}
\XDX{M}{\Sinv{A}}{A}=\DXD{\Sinv{A}}{A}{S}.
\end{equation*}
Multiplying each sides by $\D_A$ gives
\begin{equation*}
\XDX{M}{\Sinv{A}}{A}\D_A = \DXD{\Sinv{A}}{A}{S}\D_A.
\end{equation*}
By Lemma~\ref{lem_XDX_A_B},
we have $\DXD{\Sinv{A}}{A}{A}=\XDX{B}{B^T}{\inv{B}}$.
Replacing $\DXD{\Sinv{A}}{A}{A}$ by $\XDX{B}{B^T}{\inv{B}}$ on each side
yields
\begin{equation*}
X_{M} \XDX{B}{B^T}{\inv{B}} = \XDX{B}{B^T}{\inv{B}} \D_S.
\end{equation*}
The matrix $B$ is type II, therefore $\inv{B}=n^{-1}\SinvT{B}$ and 
the above is equivalent to
\begin{equation*}
\XDX{\inv{B}MB}{B^T}{\SinvT{B}}= \DXD{B^T}{\SinvT{B}}{S}.
\end{equation*}
So the first equation of the lemma follows.

Let $S'=\Theta_{\Sinv{B}}(M)$ and
\begin{equation*}
\XDX{M}{B}{\Sinv{B}} = \DXD{B}{\Sinv{B}}{S'}.
\end{equation*}
Multiplying both sides by $\D_A$,
\begin{equation*}
\XDX{M}{B}{\Sinv{B}}\D_A = \DXD{B}{\Sinv{B}}{S'}\D_A.
\end{equation*}
By Corollary~\ref{cor_T_AB_XDX}~(\ref{eqn_T_AB_XDX_d}) on $\T{A}{B^T}(A)=B^T$, 
we get $\DXD{B}{\Sinv{B}}{A}=\XDX{B}{\inv{A}}{A^T}$ and 
\begin{equation*}
X_{M} \XDX{B}{\inv{A}}{A^T} = \XDX{B}{\inv{A}}{A^T}\D_{S'}.
\end{equation*}
Since $A$ is symmetric, we have $n\inv{A} = \SinvT{A} = \Sinv{A}$ and
\begin{equation*}
\XDX{\inv{B}MB}{\Sinv{A}}{A}=\DXD{\Sinv{A}}{A}{S'}.
\end{equation*}
So the second equation of the lemma holds.  
\eop

\begin{lemma}
\label{lem_NV2_2}
Let $(A,B)$ be an invertible Jones pair with $A$ symmetric.
If $M$ lies in $\N{A}{B}$ and $N$ lies in $\N{A}{B^T}$,
then the following are equivalent:
\begin{enumerate}
\item
\label{lem_NV2_2_a}
$\T{A}{B}(M) = \T{\SinvT{B}}{\Sinv{A}}(N)$,
\item
\label{lem_NV2_2_b}
$\Theta_A(M \schur A) = \Theta_{B^T}(N \schur A)^T$,
\item
\label{lem_NV2_2_c}
$\T{\Sinv{B}}{\Sinv{A}}(M) = \T{A}{B^T}(N)$.
\end{enumerate}
\end{lemma}
{\sl Proof.}
First note that by Theorem~\ref{thm_Dim}~(a), we have
$\N{A}{B^T}=\N{A}{B}=A\schur \Nom{A}$.  Since $A$ is symmetric,
$\N{A}{B}$ is  closed under transpose.
Therefore $M^T, N^T$ also belong to $\N{A}{B}=\N{A}{B^T}$.

Now the right-hand side of (\ref{lem_NV2_2_a}) equals
$\T{\Sinv{A}}{\SinvT{B}}(N)^T$.
After applying Lemma~\ref{lem_Sinv_AB}, it becomes
$\T{A}{B^T}(N^T)^T$.  
Multiplying each side of (\ref{lem_NV2_2_a}) by $n^{-1} B^T$ gives
\begin{equation*}
\frac{1}{n} \T{A}{B}(M)B^T = \frac{1}{n} \T{A}{B^T}(N^T)^T B^T,
\end{equation*}
which is equivalent to
\begin{equation*}
\frac{1}{n} \T{A}{B}(M)\ \T{A}{B}(A)^T = 
\frac{1}{n} \left( \T{A}{B^T}(A)^T\ \T{A}{B^T}(N^T) \right)^T.
\end{equation*}
Applying the first part of Theorem~\ref{thm_TA_TB} to $(A,B)$ and 
its second part to $(A,B^T)$,
we have
\begin{equation*}
\Theta_A(M \schur A) =
\Theta_{B^T}(A\schur N)^T.
\end{equation*}
So (\ref{lem_NV2_2_a}) is equivalent to (\ref{lem_NV2_2_b}).

By Theorem~\ref{thm_TA_TB}, (\ref{lem_NV2_2_b}) can be rewritten as
\begin{equation*}
\frac{1}{n}\T{A}{B}(A)\ \T{A}{B}(M^T)^T =
\frac{1}{n}\left( \T{A}{B^T}(N)^T\ \T{A}{B^T}(A)\right)^T,
\end{equation*}
which is equal to
\begin{equation*}
B \T{A}{B}(M^T)^T = B \T{A}{B^T}(N).
\end{equation*}
So
\begin{equation*}
\T{B}{A}(M^T) = \T{A}{B^T}(N).
\end{equation*}
By Lemma~\ref{lem_Sinv_AB},
\begin{equation*}
\T{\Sinv{B}}{\Sinv{A}}(M) = \T{A}{B^T}(N).
\end{equation*}
Hence we have shown the equivalence of (\ref{lem_NV2_2_b}) and
(\ref{lem_NV2_2_c}).  
\eop

Now we are ready to determine the structure of the matrices in $\Nom{V}$
with zero off-diagonal blocks.
\begin{lemma}
\label{lem_NV2}
Let $(A,B)$ be an invertible Jones pair with $A$ symmetric.
The set of matrices
\begin{equation*}
\begin{pmatrix}
F+R & F-R & 0 & 0\\
F-R & F+R & 0 & 0\\
0 & 0 & \inv{B}FB+R_1 & \inv{B}FB-R_1\\
0 & 0 & \inv{B}FB-R_1 & \inv{B}FB+R_1\\
\end{pmatrix}
\end{equation*}
satisfying
\begin{eqnarray*}
F &\in & \Nom{A},\\ 
R &\in & \N{A}{B}, \\
\Theta_{B^T}(A \schur R_1)^T &=& \Theta_A(A \schur R)
\end{eqnarray*}
equals the subspace of $\Nom{V}$ consisting matrices with $2n \times 2n$
zero off-diagonal blocks.
\end{lemma}
{\sl Proof.}
Suppose
\begin{equation*}
Z=
\begin{pmatrix}
M_1 & N_1 & 0 & 0\\
P_1 & Q_1 & 0 & 0\\
0 & 0 & M_2 & N_2\\
0 & 0 & P_2 & Q_2\\
\end{pmatrix}
\in \Nom{V}.
\end{equation*}
Since 
\begin{equation*}
\y{\1}{\1}=\y{\2}{\2}=\Evone,
\end{equation*}
$Z$ having $\y{\1}{\1}$ and $\y{\2}{\2}$ as eigenvectors for all
$\range{i,j}{1}{n}$ implies
\begin{equation}
\label{eqn_NV2_1}
\Theta_A(M_1+N_1) = \Theta_A(P_1+Q_1)
=\Theta_{\SinvT{B}}(M_2+N_2)=\Theta_{\SinvT{B}}(P_2+Q_2).
\end{equation}
From the first and the third equalities, we know that $M_1+N_1=P_1+Q_1$ and 
$M_2+N_2=P_2+Q_2$.

By the first equation of Lemma~\ref{lem_NV2_1}, the second equality in 
(\ref{eqn_NV2_1}) holds if and only if
$M_2+N_2 = \inv{B}(M_1+N_1)B$, which is true if we let $F=\frac{1}{2}(M_1+N_1)$.
Therefore both the $(\1,\1)$- and the $(\2,\2)$-blocks
of $\Theta_V(Z)$ equal $\Theta_A(M_1+N_1)$.
It is an easy consequence that for $\range{i,j}{1}{n}$,
\begin{equation*}
\y{\1}{\2}=\y{\2}{\1}=\Evtwo
\end{equation*}
are also eigenvectors of $Z$, 
and the $(\1,\2)$- and the $(\2,\1)$-blocks of $\Theta_V(Z)$ also
equal $\Theta_A(M_1+N_1)$.

Since $M_2+N_2 = \inv{B} (M_1+N_1) B$, 
the second equation of Lemma~\ref{lem_NV2_1} with $M=M_1+N_1$
implies  the following
\begin{eqnarray*}
\XDX{(M_1+N_1)}{B}{\Sinv{B}} &=& \DXD{B}{\Sinv{B}}{S'}, \\
\XDX{(M_2+N_2)}{\Sinv{A}}{A} &=& \DXD{\Sinv{A}}{A}{S'},
\end{eqnarray*}
where $S'=\Theta_{\Sinv{B}}(M_1+N_1)=\Theta_A(\inv{B} (M_1+N_1)B)$.
Hence for all $\range{i,j}{1}{n}$,
\begin{equation*}
\y{\3}{\3}=\y{\4}{\4}=\Evthree
\end{equation*}
and 
\begin{equation*}
\y{\3}{\4}=\y{\4}{\3}=\Evfour
\end{equation*}
are eigenvectors of $Z$.
We see that both the $(\3,\3)$- and the $(\4,\4)$-blocks of
$\Theta_V(Z)$ are
\begin{equation*}
\Theta_{\Sinv{B}}(M_1+N_1) =
\Theta_B(M_1+N_1)^T=\inv{B}\Theta_A(M_1+N_1)B,
\end{equation*}
with the last equality implied by Theorem~\ref{thm_Inv_JP}.

Consider
\begin{equation*}
\y{\1}{\3} = -\y{\2}{\4} = \Evfive,
\end{equation*}
and 
\begin{equation*}
\y{\1}{\4} = -\y{\2}{\3} = \Evsix.
\end{equation*}
They are eigenvectors of $Z$ if and only if the following hold:
\begin{eqnarray*}
M_1-N_1=-P_1+Q_1 &\in &\N{A}{B},\\
M_2-N_2=-P_2+Q_2 &\in &\N{A}{B^T}=\N{A}{B},\\
\T{A}{B}(M_1-N_1) &=& \T{\SinvT{B}}{\Sinv{A}}(M_2-N_2).
\end{eqnarray*}
If we let $R=\frac{1}{2} (M_1-N_1)$ and $R_1 =\frac{1}{2} (M_2-N_2)$ then
by Lemma~\ref{lem_NV2_2}, the third equation above is equivalent to
\begin{equation}
\label{eqn_NV2_2}
\Theta_A(A \schur R)= \Theta_{B^T}(A\schur R_1)^T.
\end{equation}
So the $(\1,\3)$-, $(\2,\4)$-, $(\1,\4)$- and $(\2,\3)$-blocks of
$\Theta_V(Z)$ equals $\T{A}{B}(M_1-N_1)=2\T{A}{B}(R)$.

Now, consider
\begin{equation*}
\y{\3}{\1} = -\y{\4}{\2} = \Evseven
\end{equation*}
and 
\begin{equation*}
\y{\4}{\1} = -\y{\3}{\2} = \Eveight.
\end{equation*}
They are eigenvectors of $Z$ if and only if 
\begin{equation*}
\T{\Sinv{B}}{\Sinv{A}}(M_1-N_1) = \T{A}{B^T}(M_2-N_2).
\end{equation*}
This equation holds because it is equivalent to
Equation~(\ref{eqn_NV2_2}) by Lemma~\ref{lem_NV2_2}.
So these vectors are indeed eigenvectors of $Z$ and the
$(\3,\1)$-, $(\4,\2)$-, $(\4,\1)$- and $(\3,\2)$-blocks of
$\Theta_V(Z)$ equal
\begin{equation*}
\T{\Sinv{B}}{\Sinv{A}}(M_1-N_1)=\T{\Sinv{A}}{\Sinv{B}}(2R)^T=2\T{A}{B}(R^T)^T.
\end{equation*}

If we let $F = \frac{1}{2} (M_1+N_1)$, $R=\frac{1}{2}(M_1-N_1)$
and $R_1 = \frac{1}{2}(M_2-N_2)$, then
the result follows and
\begin{equation*}
\Theta_V(Z) = 
2
\begin{pmatrix}
\Theta_A(F) & \Theta_A(F) & \T{A}{B}(R) & \T{A}{B}(R) \\
\Theta_A(F) & \Theta_A(F) & \T{A}{B}(R) & \T{A}{B}(R) \\
\T{A}{B}(R^T)^T & \T{A}{B}(R^T)^T & \inv{B}\Theta_A(F) B & \inv{B}\Theta_A(F) B\\
\T{A}{B}(R^T)^T & \T{A}{B}(R^T)^T & \inv{B}\Theta_A(F) B & \inv{B}\Theta_A(F) B\\
\end{pmatrix}.
\end{equation*} 
\eop

\begin{corollary}
\label{cor_NV2}
Let $(A,B)$ be an invertible Jones pair with $A$ symmetric.
Suppose
$F \in \Nom{A}$ and 
$R \in \N{A}{B}$
then
\begin{equation*}
\begin{pmatrix}
\Theta_A(F) & \Theta_A(F) & \T{A}{B}(R) & \T{A}{B}(R) \\
\Theta_A(F) & \Theta_A(F) & \T{A}{B}(R) & \T{A}{B}(R) \\
\T{A}{B}(R^T)^T & \T{A}{B}(R^T)^T & \inv{B}\Theta_A(F) B & \inv{B}\Theta_A(F) B\\
\T{A}{B}(R^T)^T & \T{A}{B}(R^T)^T & \inv{B}\Theta_A(F) B & \inv{B}\Theta_A(F) B\\
\end{pmatrix}
\in \Nom{V^T}.
\end{equation*}
\eop
\end{corollary} 

If $\Nom{A}$ has dimension $r$, then Lemma~\ref{lem_NV2} says that
the subspace of $\Nom{V}$ spanned by the matrices with zero off-diagonal
blocks has dimension $2r$.
Unfortunately, we are not able to determine the dimension of 
the subspace of $\Nom{V}$ spanned by the matrices with zero diagonal
blocks.  Lemma~\ref{lem_NV3} describes the structure of this subspace.
We need the following two lemmas to prove Lemma~\ref{lem_NV3}.

\begin{lemma}
\label{lem_NV3_1}
Let $A$ and $B$ be $\nbyn$ type-II matrices.  Assume further that $A$ is
symmetric.
If $M \in \dualN{A}{B}$ then
$\T{A}{B}(R^T) = nM$ if and only if
\begin{equation*}
\XDX{M}{B^T}{\SinvT{B}} = \DXD{\Sinv{A}}{A}{R}.
\end{equation*}
\end{lemma}
{\sl Proof.}
Since $\SinvT{B}=n\inv{B}$, the equality above becomes
\begin{equation*}
n \DXD{A}{M}{B^T}=\XDX{A}{R}{B}.
\end{equation*}
Applying the Exchange Lemma yields
\begin{equation*}
n \DXD{M}{A}{B}=\XDX{A}{B}{R}.
\end{equation*}
Taking the transpose of each side, we get
\begin{equation*}
n \DXD{B}{A}{M}=\XDX{R^T}{B}{A}
\end{equation*}
and $\T{A}{B}(R^T) = nM$.  
\eop

\begin{lemma}
\label{lem_NV3_2}
Let $A$ be a symmetric matrix.
If $(A,B)$ is an invertible one-sided Jones pair
and $R\in\N{A}{B}$, then 
\begin{eqnarray*}
\T{A}{B}(R)^T &=& \inv{B}\ \T{A}{B}(R^T) B^T\\
&=& B^T\ \T{A}{B}(R^T) \inv{B}.
\end{eqnarray*}
\end{lemma}
{\sl Proof.}
Since $A$ is symmetric, $\N{A}{B}=A \schur \Nom{A}$ is closed under 
transpose and $R^T$ belongs to $\N{A}{B}$.
Applying Theorem~\ref{thm_TA_TB}, we get
\begin{equation*}
\Theta_A(R^T \schur A) = \frac{1}{n}\T{A}{B}(R^T) B^T
\end{equation*}
and
\begin{equation*}
\Theta_A(A \schur R^T) = \frac{1}{n}B\ \T{A}{B}(R)^T.
\end{equation*}
Hence the first equality follows.

Similarly, applying the same theorem, we get
\begin{equation*}
\Theta_A(R^T \schur \inv{A}) = \frac{1}{n}\T{A}{B}(R^T) \SinvT{B}
\end{equation*}
and
\begin{equation*}
\Theta_A(\inv{A} \schur R^T) = \frac{1}{n}\Sinv{B}\ \T{A}{B}(R)^T.
\end{equation*}
So the second equality holds.  
\eop

Now we are ready to examine the matrices of $\Nom{V}$ with zero diagonal
blocks.
\begin{lemma}
\label{lem_NV3}
Suppose $A$ is a symmetic matrix and $(A,B)$ is an invertible Jones
pair.
Let ${\cal S}$ be the set of matrices
\begin{equation*}
\begin{pmatrix}
0 & 0 &  G+H & G-H \\
0 & 0 &  G-H & G+H \\
\inv{B}GB^T+H_1 & \inv{B}GB^T-H_1 & 0 & 0\\
\inv{B}GB^T-H_1 & \inv{B}GB^T+H_1 & 0 & 0\\
\end{pmatrix}
\end{equation*}
satisfying
\begin{equation*}
G \in  \dualN{A}{B}
\end{equation*}
and there exist $\nbyn$ matrices $S$ and
$S_1$  such that
\begin{eqnarray*}
\XDX{\inv{A}}{\Sinv{B}}{H} \D_{\Sinv{A}} X_{\inv{B}} &=&  \D_S =
\XDX{B}{A}{H_1} \D_B X_A\\
\XDX{B^T}{A}{H} \D_{B^T} X_A &=&
\D_{S_1} =
\XDX{\inv{A}}{\SinvT{B}}{H_1} \D_{\Sinv{A}} X_{\invT{B}}.
\end{eqnarray*}
Then ${\cal S}$ equals the subspace of $\Nom{V}$ consisting matrices with 
$2n \times 2n$ zero diagonal blocks.
\end{lemma}
{\sl Proof.}
Suppose
\begin{equation*}
Z =
\begin{pmatrix}
0 & 0 & M_1 & N_1\\
0 & 0 & P_1 & Q_1\\
M_2 & N_2 & 0 & 0\\
P_2 & Q_2 & 0 & 0\\
\end{pmatrix}.
\end{equation*}
There exists an $\nbyn$ matrix $R$ satisfying the following two equations 
if and only if $\y{\1}{\1}$, $\y{\2}{\2}$, $\y{\1}{\2}$ and $\y{\2}{\1}$
are eigenvectors of $Z$.
\begin{eqnarray}
\label{eqn_NV3_1}
\XDX{P_1+Q_1}{B^T}{\SinvT{B}} =
\XDX{M_1+N_1}{B^T}{\SinvT{B}}&=&\DXD{\Sinv{A}}{A}{R},\\
\label{eqn_NV3_2}
\XDX{P_2+Q_2}{\Sinv{A}}{A} =
\XDX{M_2+N_2}{\Sinv{A}}{A} &=& \DXD{B^T}{\SinvT{B}}{R}.
\end{eqnarray}
By Lemma~\ref{lem_NV3_1}, Equation~(\ref{eqn_NV3_1}) implies 
$P_1+Q_1=M_1+N_1=n^{-1}\T{A}{B}(R^T)$.
Applying Corollary~\ref{cor_T_AB_XDX}~(\ref{eqn_T_AB_XDX_d}),
Equation~(\ref{eqn_NV3_2}) is equivalent to 
\begin{equation*}
\T{A}{B}(R) = n (M_2+N_2)^T = n(P_2+Q_2)^T.
\end{equation*}
By Lemma~\ref{lem_NV3_2}, we conclude that 
\begin{eqnarray*}
M_2+N_2 &=& \inv{B} (M_1+N_1) B^T\\
&=& B^T (M_1+N_1) \inv{B}.
\end{eqnarray*}
Note that the $(\1,\1)$- and the $(\2,\2)$-blocks of $\Theta_V(Z)$ equal
to $R$, while its $(\1,\2)$- and $(\2,\1)$-blocks equal $-R$.

We know that $M_1+N_1=\T{A}{B}(n^{-1}R^T)$ lies in $\dualN{A}{B}$ and
$M_2+N_2 = \inv{B} (M_1+N_1) B^T$.  
Now by Corollary~\ref{cor_N_AB_ABT}, 
$M_2+N_2$ also belongs to $\dualN{A}{B^T}$.
Let $R_1$ be the matrix such that $M_2+N_2 = \T{A}{B^T}(n^{-1}R_1^T)$.
By Lemma~\ref{lem_NV3_1}, this equation is equivalent to 
\begin{equation}
\label{eqn_NV3_3}
\XDX{M_2+N_2}{B}{\Sinv{B}}=\DXD{\Sinv{A}}{A}{R_1}.
\end{equation}
Applying the second equation of Lemma~\ref{lem_NV3_2} to $(A,B^T)$, we get
\begin{eqnarray*}
\T{A}{B^T}(R_1)^T &=& B\ \T{A}{B^T}(R_1^T) \invT{B}\\
&=& n B (M_2+N_2) \invT{B} \\
&=& n (M_1+N_1).
\end{eqnarray*}
Applying Corollary~\ref{cor_T_AB_XDX}~(\ref{eqn_T_AB_XDX_d}),
$\T{A}{B^T}(R_1) = n (M_1+N_1)^T$ is equivalent to
\begin{equation}
\label{eqn_NV3_4}
\XDX{M_1+N_1}{\Sinv{A}}{A} = \DXD{B}{\Sinv{B}}{R_1}.
\end{equation}

Equations~(\ref{eqn_NV3_3}) and (\ref{eqn_NV3_4}) imply that
$\y{\3}{\3}$, $\y{\4}{\4}$, $\y{\3}{\4}$ and $\y{\4}{\3}$ are 
eigenvectors of $Z$.
Moreover, the $(\3,\3)$- and the $(\4,\4)$- blocks of $\Theta_V(Z)$ 
equal $R_1$, while its $(\3,\4)$- and $(\4,\3)$-blocks equal $-R_1$.

Now 
$\y{\1}{\3}$, $\y{\2}{\4}$, $\y{\1}{\4}$ and $\y{\2}{\3}$ 
are eigenvectors of $Z$
if and only if there exists an $\nbyn$ matrix $S$ such that
\begin{eqnarray*}
d^{-1} \XDX{M_1-N_1}{\Sinv{A}}{\SinvT{B}} &=& d\DXD{B}{A}{S}\\
d\XDX{M_2-N_2}{B}{A} &=& d^{-1} \DXD{\Sinv{A}}{\SinvT{B}}{S}
\end{eqnarray*}
which is equivalent to the first equation in the lemma
if we let $H=M_1-N_1$ and $H_1=M_2-N_2$.
The $(\1,\3)$- and the $(\2,\4)$ blocks of $\Theta_Z(V)$ are $S$,
while its $(\1,\4)$- and $(\2,\3)$-blocks are $-S$.

Similarly, there exists an $\nbyn$ matrix $S_1$ such that
\begin{eqnarray*}
d \XDX{M_1-N_1}{B^T}{A} &=& d^{-1}\DXD{\Sinv{A}}{\Sinv{B}}{S_1}\\
d^{-1}\XDX{M_2-N_2}{\Sinv{A}}{\Sinv{B}} &=& d \DXD{B^T}{A}{S_1}
\end{eqnarray*}
if and only if 
$\y{\3}{\1}$, $\y{\4}{\2}$, $\y{\3}{\2}$ and $\y{\4}{\1}$ 
are eigenvectors of $Z$.
The above equations are equivalent to the second equation in the lemma.
Hence the $(\3,\1)$- and the $(\4,\2)$ blocks of $\Theta_Z(V)$ are $S_1$,
while its $(\3,\2)$- and $(\4,\1)$-blocks are $-S_1$.

The lemma follows after we let $G=M_1+N_1$, $H=M_1-N_1$ and
$H_1=M_2-N_2$.
In addition, it is worth noting that
\begin{equation*}
\Theta_V(Z) = 
\begin{pmatrix}
R & -R & S & -S \\
-R & R & -S & S \\
S_1 & -S_1 & R_1 & -R_1\\
-S_1 & S_1 & -R_1 & R_1\\
\end{pmatrix}.
\end{equation*} 
\eop

The lemma above describes the matrices in $\Nom{V}$ with zero diagonal
blocks.  The dimension of this subspace of $\Nom{V}$ is greater than or
equal to the dimension of $\dualN{A}{B}$.  The matrix $H_1$ is
determined by $H$.  However, we do not understand the conditions on $H$
listed in the lemma.  We only know that $\Nom{V}$ has dimension at least
three times the dimension of $\Nom{A}$.  

\begin{theorem}
\label{thm_NV}
Suppose $(A,B)$ is an invertible Jones pair with $A$ symmetric.
Let ${\cal S}$ be the set of matrices of the following form
\begin{equation*}
\begin{pmatrix}
F+R & F-R & G+H & G-H\\
F-R & F+R & G-H & G+H\\
\inv{B}GB^T+H_1 & \inv{B}GB^T-H_1 & 
\inv{B}FB+R_1 & \inv{B}FB-R_1\\
\inv{B}GB^T-H_1 & \inv{B}GB^T+H_1 & 
\inv{B}FB-R_1 & \inv{B}FB+R_1\\
\end{pmatrix}
\end{equation*}
where $F \in \Nom{A}$, $R \in \N{A}{B}$, $G \in \dualN{A}{B}$,
$H$ and $H_1$ satisfy the conditions in Lemma~\ref{lem_NV3}
and $R_1$ such that
\begin{equation*}
\Theta_{B^T}(A \schur R_1)^T = \Theta_A(A \schur R).
\end{equation*}
Then the Nomura algebra $\Nom{V}$ equals ${\cal S}$.
\eop
\end{theorem}

\begin{theorem}
\label{thm_V}
Let $A$ and $B$ be type-II matrices.
If $A$ is symmetric, then
the $4n \times 4n$ matrix 
\begin{equation*}
V = 
\begin{pmatrix}
dA & -dA & \Sinv{B} & \Sinv{B} \\ 
-dA & dA & \Sinv{B} & \Sinv{B} \\ 
\SinvT{B} & \SinvT{B} & dA & -dA \\
\SinvT{B} & \SinvT{B} & -dA & dA \\
\end{pmatrix}
\end{equation*}
is a symmetric spin model with loop variable $2d$
if and only if $(A,B)$ is an invertible Jones pair.
\end{theorem}
{\sl Proof.}
Suppose $(A,B)$ is an invertible Jones pair.
If we let $F=H=H_1=0$, $R_1=R=dA$ and $G=\Sinv{B}$ in Theorem~\ref{thm_NV},
then we have $V \in \Nom{V}$.

Let 
\begin{equation*}
\breve{A} = 
\begin{pmatrix}
dA & -dA & 0 & 0 \\
-dA & dA & 0 & 0 \\
0 & 0 & dA & -dA\\
0 & 0 & -dA & dA\\
\end{pmatrix}
\end{equation*}
and 
\begin{equation*}
\ddot{B}=
\begin{pmatrix}
0 & 0 & \Sinv{B} & \Sinv{B}\\
0 & 0 & \Sinv{B} & \Sinv{B}\\
\SinvT{B} & \SinvT{B} & 0 & 0 \\
\SinvT{B} & \SinvT{B} & 0 & 0 \\
\end{pmatrix}.
\end{equation*}
From the proof of Lemma~\ref{lem_NV2}, we know that
\begin{equation*}
\Theta_V (\breve{A}) 
=
2d
\begin{pmatrix}
0 & 0 & B & B\\
0 & 0 & B & B\\
B^T & B^T & 0 & 0 \\
B^T & B^T & 0 & 0 \\
\end{pmatrix}.
\end{equation*}
Moreover, from the proof of Lemma~\ref{lem_NV3}, we have
\begin{equation*}
\Theta_V (\ddot{B})
=
2d
\begin{pmatrix}
d^{-1}\Sinv{A} & -d^{-1}\Sinv{A} & 0 & 0 \\
-d^{-1}\Sinv{A} & d^{-1}\Sinv{A} & 0 & 0 \\
0 & 0 & d^{-1}\Sinv{A} & -d^{-1}\Sinv{A}\\
0 & 0 & -d^{-1}\Sinv{A} & d^{-1}\Sinv{A}\\
\end{pmatrix}.
\end{equation*}
So we have $\Theta_V(V) = 2d \Sinv{V}$ and
$(V, 2d\Sinv{V})$ is therefore an invertible one-sided Jones pair.
It follows from Lemmas~\ref{lem_W_JP} and \ref{lem_JP_Spin} that $V$ is in
fact a spin model with loop variable $2d$.

Conversely, if $(V,2d\Sinv{V})$ is an invertible Jones pair, then
\begin{equation*}
V \y{\1}{\3} = 2d \Sinv{V} \y{\1}{\3}
\end{equation*}
implies $A \ (\Evector{A}{i}{B}{j}) = B_{i,j} \ (\Evector{A}{i}{B}{j})$.
So $(A,B)$ is an invertible one-sided Jones pair.
Furthermore,
\begin{equation*}
V \y{\3}{\1} = 2d \Sinv{V} \y{\3}{\1}
\end{equation*}
implies $A \ (\Evector{A}{i}{B^T}{j}) = B_{j,i} \ (\Evector{A}{i}{B^T}{j})$.
So $(A,B^T)$ is also an invertible one-sided Jones pair.
Hence $(A,B)$ is an invertible Jones pair.  
\eop

\begin{theorem}
\label{thm_V'}
Let $(A,B)$ be an invertible Jones pair with $A$ symmetric.
Then the matrix
\begin{equation*}
V' = 
\begin{pmatrix}
dA & -dA & -\Sinv{B} & -\Sinv{B}\\
-dA & dA & -\Sinv{B} & -\Sinv{B}\\
-\SinvT{B} & -\SinvT{B} & dA & -dA\\
-\SinvT{B} & -\SinvT{B} & -dA & dA\\
\end{pmatrix}
\end{equation*}
is a symmetric spin model with loop variable $-2d$.
Moreover, $\Theta_{V'} = \Theta_V$.
\end{theorem}
{\sl Proof.}
Let 
\begin{equation*}
D=
\begin{pmatrix}
I_{2n} & 0 \\ 0 & -I_{2n}\\
\end{pmatrix}
\end{equation*}
Note that $V' = DVD$.
By Lemma~\ref{lem_Scaling}, we conclude that $\Nom{V'}=\Nom{V}$.
By setting $F=H=H_1=0$, $R_1=R=dA$ and $G=-\Sinv{B}$ in
Theorem~\ref{thm_NV},
we have $V' \in \Nom{V'}$.  
Moreover, by the proof of the previous theorem, we have
\begin{eqnarray*}
\Theta_V(V') &=& 
2d 
\begin{pmatrix}
0 & 0 & B & B\\
0 & 0 & B & B\\
B^T & B^T & 0 & 0 \\
B^T & B^T & 0 & 0 \\
\end{pmatrix}
- 2d
\begin{pmatrix}
d^{-1}\Sinv{A} & -d^{-1}\Sinv{A} & 0 & 0 \\
-d^{-1}\Sinv{A} & d^{-1}\Sinv{A} & 0 & 0 \\
0 & 0 & d^{-1}\Sinv{A} & -d^{-1}\Sinv{A}\\
0 & 0 & -d^{-1}\Sinv{A} & d^{-1}\Sinv{A}\\
\end{pmatrix}
\\
&=& -2d \Sinv{V'}.
\end{eqnarray*}

Since $D^2=I_{4n}$,
\begin{eqnarray*}
\Evector{V'}{h}{\Sinv{(V')}}{k} &=&
\Evector{(DVD)}{h}{D\Sinv{V}D}{k}\\
&=& \Evector{(VD)}{h}{\Sinv{V}D}{k}\\
&=&\pm \Evector{(V)}{h}{\Sinv{V}}{k},
\end{eqnarray*}
for all $\range{h,k}{1}{4n}$.
So for any $M\in \Nom{V}$, we have 
$\Theta_{V'}(M)_{h,k}= \Theta_{V}(M)_{h,k}$.
As a result, the two duality maps $\Theta_V$ and $\Theta_{V'}$ are
identical.
Further, $V'$ is a spin model because
\begin{equation*}
\Theta_{V'}(V')=\Theta_V(V')=-2d\Sinv{(V')}.
\end{equation*}
\eop
\section{The Modular Invariance Equation}
\label{JP_MI}

In \cite{Ban1}, Bannai asked for a strategy to find four-weight spin
models.  In this section, we answer her question.
In \cite{BBJ}, Bannai, Bannai and Jaeger used the modular invariance
equation to design a method that exhaustively searches for spin models.
Their design relies on the fact that a spin model always belongs to some
formally self-dual Bose-Mesner algebra.  Since this is usually not the
case for the matrices in a four-weight spin model, their method does not
apply directly.

Given an $\nbyn$ invertible Jones pair $(A,B)$, 
there exist two $4n \times 4n$ spin models, $V$ and $V'$, 
with $A$ and $\Sinv{B}$ as their submatrices.  
So we can apply Bannai, Bannai and Jaeger's method to formally
self-dual Bose-Mesner algebras on $4n$ elements, hoping to find a
symmetric spin model with the same structure as $V$.  From $V$, we
retrieve $A$ and $B$.  The algorithm we outline at the end of 
this section is the only known strategy to find four-weight spin models.

Suppose that $(A,B)$ is an invertible Jones pair and that $\Nom{V}$ is
the Nomura algebra constructed from this pair.
We want to see how we can recover $(A,B)$ from this algebra.
Let $\A = \{ \seq{A}{0}{m}\}$ be the basis of Schur idempotents of 
$\Nom{V}$.  Let ${\cal E}=\{\seq{E}{0}{m}\}$ be the basis of 
principal idempotents of $\Nom{V}$ such that 
$\Theta_V(E_i) = A_i$ and $\Theta_V(A_i)=nE_i^T$.
Let $P$ be the matrix of eigenvalues with respect to this ordering.
So if $T$ is the matrix representing the transpose map with respect to $\A$,
then $P^2=nT$.  Note that the ordering of the matrices in ${\cal E}$ and 
the matrix of eigenvalues depend on the duality map $\Theta_V$.

From Lemma~\ref{lem_NV1}, we can always define
two sets of indices $\I, \J \subset \{0,\ldots,m\}$
satisfying
\begin{equation}
\label{eqn_I_J}
\sum_{i\in \I} A_i^T =
\begin{pmatrix}
J_n & J_n & 0 & 0\\
J_n & J_n & 0 & 0\\
0 & 0 & J_n & J_n\\
0 & 0 & J_n & J_n\\
\end{pmatrix},
\quad
\sum_{j\in \J} A_j^T =
\begin{pmatrix}
0 & 0 & J_n & J_n\\
0 & 0 & J_n & J_n\\
J_n & J_n & 0 & 0\\
J_n & J_n & 0 & 0\\
\end{pmatrix}.
\end{equation}
So $\I \cup \J = \{0,\ldots,m\}$.

Suppose $V=\sum_{k=0}^m t_k A_k^T$.
Then $V'= \sum_{i\in \I} t_i A_i^T - \sum_{j\in\J} t_j A_j^T$.
Define two $(m+1) \times (m+1)$ diagonal matrices $D_{\I}$ and $D_{\J}$ as
\begin{equation*}
(D_{\I})_{i,i} = 
\begin{cases}
t_i &\text{if }  i\in \I,\\
0 &\text{if }  i \in \J
\end{cases}
\end{equation*}
and
\begin{equation*}
(D_{\J})_{j,j} = 
\begin{cases}
t_j &\text{if }  j\in \J,\\
0 &\text{if }  j \in \I.
\end{cases}
\end{equation*}
By Theorem~\ref{thm_Mod_Inv}, 
since $V$ is a spin model with loop variable $2d$, we get
\begin{equation*}
\left( P(D_{\I}+D_{\J}) \right)^3 = 8d^3 t_0\ I_{m+1}.
\end{equation*}
Similarly, since $V'$ is a spin model with loop variable $-2d$ and 
$\Theta_{V'}=\Theta_V$, we have
\begin{equation*}
\left( P(D_{\I}-D_{\J}) \right)^3 = -8d^3 t_0\ I_{m+1}.
\end{equation*}

Define $D_{\I}^{-}$ to be the diagonal matrix with its $ii$-entry
equals to $t_i^{-1}$ if $i\in \I$, and zero otherwise.
Define $D_{\J}^{-}$ similarly.
Then the above two equations can be rewritten as
\begin{eqnarray*}
(D_{\I}+D_{\J}) P (D_{\I}+D_{\J}) &=&  
8d^3 t_0\ \inv{P} ({D_{\I}^{-}} + {D_{\J}^{-}}) \inv{P},\\
(D_{\I}-D_{\J}) P (D_{\I}-D_{\J}) &=&  
-8d^3 t_0\ \inv{P} ({D_{\I}^{-}} - {D_{\J}^{-}}) \inv{P},
\end{eqnarray*}
which are equivalent to
\begin{eqnarray}
\label{eqn_MI_1}
D_{\I} P D_{\I} + D_{\J} P D_{\J} &=& 8d^3 t_0\ \inv{P} {D_{\J}^{-}}P,\\
\label{eqn_MI_2}
D_{\I} P D_{\J} + D_{\J} P D_{\I} &=& 8d^3 t_0\ \inv{P} {D_{\I}^{-}}P.
\end{eqnarray}

So Equations~(\ref{eqn_MI_1}) and (\ref{eqn_MI_2}) are necessary conditions
on $\Nom{V}$ for the existence of the invertible Jones pair $(A,B)$
where $A$ is symmetric.
By Lemma~\ref{lem_NV1}, the matrix
$I_2 \otimes J_{2n} \in \Nom{V}$.
We say that $\Nom{V}$ is the Bose-Mesner algebra of an \SL{imprimitive}\ 
association scheme.
Note that it is possible for a formally self-dual Bose-Mesner to have
more than one duality map, each of which has a corresponding 
matrix of eigenvalues.

Here is the strategy for constructing invertible Jones pairs:
\begin{enumerate}
\item
Given a formally self-dual Bose-Mesner algebra on $4n$ elements that contains
$I_2 \otimes J_{2n}$, define $\I$ and $\J$ as in Equation~(\ref{eqn_I_J}).
\item
Enumerate its duality maps.
\item
For each duality map, form the matrix of eigenvalues $P$.
\item
Solve Equations~(\ref{eqn_MI_1}) and (\ref{eqn_MI_2}) for the $t_i$'s.
\item
For each solution, compute the matrix
$V =\sum_{i\in \I} t_i A_i^T + \sum_{j\in\J} t_j A_j^T$.
\item
Check if it has the same structure as the matrix described in Theorem~\ref{thm_V}.
If so, retrieve $A$ and $B$.
\item
Verify if $(A,B)$ is an invertible Jones pair.
\end{enumerate}

If the formally self-dual Bose-Mesner algebra we start with is not the
$\Nom{V}$ constructed from some invertible Jones pair $(A,B)$, then
we have three possibile outcomes at each iteration:
\begin{description}
\item[Step~(d):]
We have no solution at this step.
\item[Step~(f):]
We have a solution at Step~(d).
But the $4n \times 4n$ matrix $V$ constructed does not have the 
desired structure.
\item[Step~(g):]
The $4n \times 4n$ matrix $V$ has the desired structure, but the pair
$(A,B)$ retrieved is not an invertible Jones pair.
\end{description}

Note that any invertible Jones pair $(A,B)$ with $A$ symmetric could be
found  by this method.  
It follows from Lemma~\ref{lem_Odd_Sym} that
this method does provide all invertible Jones pairs up to odd-gauge
equivalence. 
Although this method is admittedly quite inefficient, it is all we have.

\section{Quotients}
\label{JP_Quotient}

Up to this point, we have four Bose-Mesner algebras constructed from 
an invertible Jones pair $(A,B)$ where $A$ is symmetric.
In this section and the next one, 
we examine the relations among these algebras.
We summarize them in a picture at the end of the next section.
They are interesting because the constructions in
Sections~\ref{N_W} and \ref{Ext} are done independently.  
In particular, we have not been able to spot any connections or 
common features between these two constructions.  
However their resulting Bose-Mesner algebras are 
closely related to each other.

When $A$ is symmetric, we have $A^T=I^{-2} A I^2$.
So we let $C=I$ and define
\begin{equation*}
W = 
\begin{pmatrix} 
A^T & -A^T \\ \SinvT{B} & \SinvT{B} \\
\end{pmatrix}.
\end{equation*}

Suppose $\BM$ is the Bose-Mesner algebra of an association scheme
on $n$ elements.
Let $\pi = (\seq{C}{1}{r})$ be a partition of $\{1,\ldots,n\}$.
Define the characteristic matrix $S$ of $\pi$ to be the $n \times r$ matrix
with
\begin{equation*}
S_{i,k} =
\begin{cases}
1 & \text{if } i \in C_k,\\
0 & \text{otherwise}.
\end{cases}
\end{equation*}
We say $\pi$ is \SL{equitable relative to $\BM$}\  if and only if
for all $M \in \BM$, there exists matrix $B_M$ satisfying
\begin{equation*}
MS=SB_M.
\end{equation*}
We call the $r \times r$ matrix $B_M$ the \SL{quotient of $M$
with respect to $\pi$}.
We name the set $\{ B_M : M\in \BM\}$ the \SL{quotient of $\BM$ with respect to $\pi$}.

\begin{theorem}
\label{thm_Quotient_NA}
For $\range{i}{1}{n}$, let $C_i = \{i ,n+i \}$
and let $\pi=(\seq{C}{1}{n})$.
If $(A,B)$ is an invertible Jones pair with $A$ symmetric, then
$\Nom{A}$ is the quotient of $\Nom{W^T}$ with respect to $\pi$.
\end{theorem}
{\sl Proof.}
It follows from Theorem~\ref{thm_NW} that any matrix in 
$\Nom{W^T}$ can be expressed as
\begin{eqnarray*}
M &=& 
\begin{pmatrix}
F & F \\F & F \\
\end{pmatrix}
+
\begin{pmatrix}
R^T & -R^T \\-R^T & R^T\\
\end{pmatrix}\\
&=&
\begin{pmatrix}
F + R^T & F- R^T \\F- R^T & F + R^T\\
\end{pmatrix}
\end{eqnarray*}
for some $F \in \Nom{A}$ and $R \in \N{A}{B}$.
The characteristic matrix of $\pi$ is 
\begin{equation*}
S = 
\begin{pmatrix} I_n \\ I_n \end{pmatrix},
\end{equation*}
and 
\begin{equation*}
MS = 
\begin{pmatrix} 2F \\ 2F \end{pmatrix}
= S (2F).
\end{equation*}
So the quotient of $\Nom{W^T}$ with respect to $\pi$ equals $\Nom{A}$.  
\eop
Although we are assuming $A$ is symmetric throughout this section,
the above proof does not require this condition.

\begin{theorem}
\label{thm_Quotient_NW}
For $i=1,\ldots,n,2n+1,\ldots,3n$, let $C_i = \{i ,n+i \}$
and let $\pi=(\seq{C}{1}{n},\seq{C}{2n+1}{3n})$.
If $(A,B)$ is an invertible Jones pair with $A$ symmetric, then
$\Nom{W}$ is the quotient of $\Nom{V}$ with respect to $\pi$.
\end{theorem}
{\sl Proof.}
The characteristic matrix of $\pi$ is
\begin{equation*}
S = 
\begin{pmatrix}
I_n & 0 \\ I_n & 0\\ 0& I_n \\0& I_n\\
\end{pmatrix}.
\end{equation*}
Let $N$ be a matrix in $\Nom{V}$.
By Theorem~\ref{thm_NV}, we have
\begin{equation*}
NS = 2
\begin{pmatrix}
F & G \\ 
F & G \\
\inv{B}G B^T & \inv{B} F B\\
\inv{B}G B^T & \inv{B} F B\\
\end{pmatrix}
= S (2
\begin{pmatrix}
F & G \\ 
\inv{B}G B^T & \inv{B} F B\\
\end{pmatrix}),
\end{equation*}
for some $F\in \Nom{A}$ and $G \in \dualN{A}{B}$.

Now $F$ belongs to $\Nom{A}$ which equals $\Nom{A^T}$, 
so there exists matrix $M$ in $\Nom{A}$ such that $F=\Theta_A(M)$.
By Theorem~\ref{thm_Inv_JP}, $\inv{B} F B = \Theta_B(M)^T$.
Since $G \in \dualN{A}{B}$, there exists $M_1 \in \N{A}{B}$ 
such that $G = \T{A}{B}(M_1)$.
By Lemma~\ref{lem_NV3_2}, we know $\inv{B}G B^T = \T{A}{B}(M_1^T)^T$.

So the quotient of $\Nom{V}$ with respect to $\pi$ is
\begin{equation*}
\left\{
\begin{pmatrix}
\Theta_A(M) & \T{A}{B}(M_1) \\ \T{A}{B}(M_1^T)^T & \Theta_B(M)^T\\
\end{pmatrix}
: M \in\Nom{A}, M_1\in\N{A}{B}
\right\},
\end{equation*}
which equals $\Nom{W}$ according to Theorem~\ref{thm_NW}.  
\eop

\section{Induced Schemes}
\label{JP_Induced_Schemes}

We now complete the picture of the relations of the association schemes 
obtained from an invertible Jones pair.

Suppose $\A = \{\seq{A}{0}{d}\}$ is an association scheme on $n$ elements.
Let $Y$ be a non-empty subset of $\{1,\ldots,n\}$.
For any $\nbyn$ matrix $M$, we use $M_Y$ to denote the 
$|Y| \times |Y|$ matrix obtained from the rows and the columns of $M$ indexed
by the elements in $Y$.
We define $\A_Y=\{(A_0)_Y,\ldots,(A_d)_Y\}$.
If $\A_Y$ is also an association scheme, we say it is 
an \SL{induced scheme}\  of $\A$.

\begin{theorem}
\label{thm_Induced_Scheme_NA}
If $(A,B)$ is an invertible Jones pair, then $\Nom{A}$
is an induced scheme of $\Nom{W}$.
\end{theorem}
{\sl Proof.}
Let $Y=\{1,\ldots,n\}$.  
The result follows immediately from Theorem~\ref{thm_NW}.  
\eop

\begin{theorem}
\label{thm_Induced_Scheme_NWT}
Let $A$ be an $\nbyn$ symmetric matrix.
If $(A,B)$ is an invertible Jones pair, then
$\Nom{W^T}$ is an induced scheme of $\Nom{V}$.
\end{theorem}
{\sl Proof.}
Let $Y=\{1,\ldots,2n\}$.
Define $(\Nom{V})_Y$ to be the set $\{M_Y : M \in \Nom{V}\}$.
By Theorem~\ref{thm_V}, we have
\begin{equation*}
(\Nom{V})_Y = 
\left\{
\begin{pmatrix}
F+R & F-R \\ F-R & F+R\\
\end{pmatrix}
: F \in \Nom{A}, R \in \N{A}{B} \right\}.
\end{equation*}
By Theorem~\ref{thm_NW}, $\Nom{W^T}$ is spanned by
\begin{equation*}
\begin{pmatrix}
F & F \\ F & F \\
\end{pmatrix},
\quad
\begin{pmatrix}
R^T & -R^T \\ -R^T & R^T \\ 
\end{pmatrix},
\end{equation*}
for all $F \in \Nom{A}$ and $R \in \N{A}{B}$.
Since $A$ is symmetric, both $\Nom{A}$ and $\N{A}{B}$ are closed under
the transpose.  
Therefore $R^T\in \N{A}{B}$ and $\Nom{W^T} = (\Nom{V})_Y$.  
\eop

The following diagram summarizes the relations among the four Bose-Mesner
algebras obtained from an invertible Jones pair $(A,B)$ where $A$ is
symmetric.

\setlength{\unitlength}{1mm}
\begin{center}
\begin{picture}(85,65)(0,0)
\put(40,10){\line(1,1){20}}
\put(40,10){\line(-1,1){20}}
\put(20,30){\line(1,1){20}}
\put(60,30){\line(-1,1){20}}
\put(40,52){\makebox(0,0)[b]{$\Nom{V}$}}
\put(40,5){\makebox(0,0)[b]{$\Nom{A}$}}
\put(11,30){\makebox(0,0)[l]{$\Nom{W}$}}
\put(61,30){\makebox(0,0)[l]{$\Nom{W^T}$}}
\put(52,40){\makebox(0,0)[l]{induced scheme}}
\put(8,40){\makebox(0,0)[l]{quotient}}
\put(54,20){\makebox(0,0)[l]{quotient}}
\put(0,20){\makebox(0,0)[l]{induced scheme}}
\end{picture}
\end{center}

\section{Subschemes}
\label{JP_Subschemes}

Now we give a stronger version of the two theorems in the previous
section.

Godsil noticed that if we let $U=\{2n+1,\ldots,4n\}$, then the set of $2n
\times 2n$ matrices of $\Nom{V}$ induced by $U$ also equals $\Nom{W^T}$.
From this observation, he spotted a subscheme of $\Nom{V}$ 
which is a union of two copies of $\Nom{W^T}$. 
This fact is stronger than Theorem~\ref{thm_Induced_Scheme_NWT}.
It turns out that similar situation happens to $\Nom{W}$, which we describe
below.

Suppose $\BM$ is the Bose-Mesner algebra of some association scheme
$\A$.  If $\BM' \subset \BM$ is also a Bose-Mesner algebra, we call its
corresponding association scheme the \SL{subscheme}\  of $\A$.

\begin{lemma}
\label{lem_Subscheme_NV}
Let $U=\{2n+1,\ldots,4n\}$.  Then the set
\begin{equation*}
S=\{ M_U : M \in \Nom{V}\}
\end{equation*}
equals $\Nom{W^T}$.
\end{lemma}
{\sl Proof.}
From Lemma~\ref{lem_NV2}, we see that $S$ is the span of
\begin{equation*}
\left\{ 
\begin{pmatrix}
\inv{B} F B & \inv{B} F B \\
\inv{B} F B & \inv{B} F B \\
\end{pmatrix}, \quad
\begin{pmatrix}
R_1 & -R_1 \\ -R_1 & R_1\\
\end{pmatrix}:
F \in \Nom{A}, A \schur R_1 \in \Nom{B^T} 
\right\}.
\end{equation*}
Since $\Nom{A}=\Nom{A^T}=\Nom{B^T}=\Nom{B}$,
we have $\inv{B}FB \in \Nom{A}$ by Theorem~\ref{thm_Dim}~(d). 
Moreover, part~(a) of the same theorem implies that $R_1 \in \N{A}{B}$.
Consequently, by Theorem~\ref{thm_NW}, we conclude that $S$ and
$\Nom{W^T}$ are equal.
\eop

Now we define $\widehat{S}$ to be the set of $4n \times 4n$ matrices
having the form
\begin{equation*}
\begin{pmatrix}
F+R & F-R & 0 & 0\\
F-R & F+R & 0 & 0\\
0 & 0 & \inv{B}FB+R_1 & \inv{B}FB-R_1\\
0 & 0 & \inv{B}FB-R_1 & \inv{B}FB+R_1\\
\end{pmatrix},
\end{equation*}
where $F \in \Nom{A}$, $R \in \N{A}{B}$ and
$\Theta_{B^T}(A \schur R_1)^T=\Theta_A(A \schur R)$.
We also define
\begin{equation*}
\widehat{J} = 
\begin{pmatrix}
0 & 0 & J_n & J_n\\
0 & 0 & J_n & J_n\\
J_n & J_n & 0 & 0\\
J_n & J_n & 0 & 0\\
\end{pmatrix}.
\end{equation*}

\begin{theorem}
\label{thm_Subscheme_NV}
The space $\BM'$ spanned by the matrices in $\widehat{S} \cup \{\widehat{J}\}$
is a Bose-Mesner algebra contained in $\Nom{V}$.
\end{theorem}
{\sl Proof.}
It is straightforward to show that this space is a commutative algebra
containing $I_{4n}$ and $J_{4n}$, 
and it is also closed under the Schur product and transpose.
\eop

Note that if $\Nom{A}$ has dimension $r$, 
then $\BM'$ has dimension $2r + 1$.

Similarly, if we let $U=\{n+1,\ldots,2n\}$, then $\Nom{A}$ equals the set of matrices
$\{M_U : M\in \Nom{W}\}$.
Moreover, if $\seq{E}{0}{r-1}$ is the basis of principal idempotents of
$\Nom{A}$, then the following set is a subscheme of $N_W$ with
$r$ classes:
\begin{equation*}
\left\{ 
\begin{pmatrix}
\Theta_A(E_i) & 0 \\ 0 & \Theta_B(E_i)^T \\
\end{pmatrix} 
: \range{i}{0}{r-1}
\right\}
\bigcup
\left\{
\begin{pmatrix}
0 & J_n\\J_n & 0\\
\end{pmatrix}
\right\}.
\end{equation*}

\section{Dimension Two}
\label{JP_Dim2}

In this section, we look at the simplest kind of invertible Jones pairs.

We define the \SL{dimension of a Jones pair}\  $(A,B)$ to be the dimension of
$\N{A}{B}$, and we define the \SL{degree of $(A,B)$}\  to be the number of
distinct entries of $B$.
If $B$ is invertible, then the degree of $(A,B)$ is at least two.
Moreover, the number of distinct entries of $B$ equals the number of
eigenvalues of $A$, and $A \in \N{A}{B}$.  
Therefore the dimension of $(A,B)$ is greater than
or equal to its degree.
In the following, we consider invertible Jones pairs of dimension two.

\begin{lemma}
\label{lem_Dim2_B}
If $(A,B)$ is an $\nbyn$ invertible Jones pair of dimension two, then
$B = a(J-N)+bN$, where $N$ is the incidence matrix of some symmetric
design.
\end{lemma}
{\sl Proof.}
Let $B = a(J-N)+bN$ for some $01$-matrix $N$ and some non-zero scalars
$a$ and $b$.
By Corollary~\ref{cor_NW23}, 
\begin{equation*}
\begin{pmatrix}
0 & B \\ B^T &0 \\
\end{pmatrix}
\in \Nom{W},
\end{equation*}
which implies
\begin{equation*}
\widehat{N}=
\begin{pmatrix}
0 & N \\ N^T & 0 \\
\end{pmatrix}
\in \Nom{W}.
\end{equation*}
Since $\Nom{W}$ is closed under matrix multiplication, we have
\begin{equation*}
\widehat{N}^2=
\begin{pmatrix}
NN^T & 0 \\ 0 & N^T N\\
\end{pmatrix}
\in \Nom{W}.
\end{equation*}
By Theorem~\ref{thm_NW}, we know that $N N^T$ belongs to $\Nom{A}$.
Since $\Nom{A} = \dim \N{A}{B}$ have dimension two,
$\Nom{A}$ equals to the span of $\{I,J\}$.
Therefore we can write 
\begin{equation*}
N N^T = \lambda (J-I) + kI,
\end{equation*}
for some integers $\lambda$ and $k$.
As a result, $N$ is the incidence matrix of a symmetric $(n,k,\lambda)$-design.
\eop
This lemma is a weaker version of Bannai and Sawano's result \cite{BS1}.
They showed that $N$ is the incidence matrix of a symmetric design whose
derived design with respect to any block is a quasi-symmetric design.

Later, Godsil showed that if $A$ is an $\nbyn$ symmetric matrix and 
$(A,B)$ is an invertible Jones pair of dimension two, then
$A$ comes from a regular two-graph.
A \SL{regular two-graph}\  is a symmetric matrix $M$ that satisfies
the following conditions:
\begin{enumerate}[(a)]
\item
all diagonal entries of $M$ are zero and all off-diagonal entries of $M$
are equal to $\pm 1$, and
\item
$M$ has quadratic minimal polynomial.
\end{enumerate}
Now we give Godsil's argument.
\begin{lemma}
\label{lem_Dim2_A}
If $A$ is an $\nbyn$ symmetric matrix and 
$(A,B)$ is an invertible Jones pair of dimension two, then
$A = cI+dM$, where $M$ is a regular two-graph.
\end{lemma}
{\sl Proof.}
By Theorem~\ref{thm_Dim}~(a), we have $A \schur A \in \Nom{A}$.
Since $\Nom{A}$ equals to the span of $\{I,J\}$, there exist
non-zero scalars $a$ and $b$ such that $A\schur A=aI+b(J-I)$.
Hence we can write
\begin{equation*}
A = c I + d M,
\end{equation*}
for some symmetric matrix $M$ such that $M \schur I = 0 $ and $M\schur M = J-I$.
So the off-diagonal entries of $M$ equal to $\pm 1$.

Moreover, since $\N{A}{B}$ has dimension two, the minimal polynomial
of $A$ is quadratic.  As a result, $M$ also has quadratic minimal polynomial,
and so it is a regular two-graph.
\eop

When $(A,B)$ is an $\nbyn$ invertible Jones pair of dimension two, then
$B$ comes from a symmetric design on $n$ points.  
From such design, we can construct a bipartite distance regular
graph on $2n$ vertices with diameter three.  
Theorem~\ref{thm_NW} tells us that $\Nom{W}$ is 
the Bose-Mesner algebra of this graph.
Furthermore, if $A$ is symmetric, then $A$ comes from a regular
two-graph $M$.  By the same theorem, we see that $\Nom{W^T}$ is the
Bose-Mesner algebra of the association scheme associated to $M$.
For more information on two-graphs, please see Seidel's survey in \cite{Seidel}.

\section{Unfinished Business}
\label{JP_Future}

We believe Jones pairs are the natural way to view problems on
spin models and four-weight spin models.  
We hope to extend their existing theory using Jones pairs.
In a narrower scope, we feel that our results have great potential for
extensions.
Now we propose several research problems.

Firstly, from Theorem~\ref{thm_NV}, we know that the dimension of $\Nom{V}$ is
bounded between $3r$ and $3r+n$ where $r$ is the dimension of $\Nom{A}$.
We conjecture that
\begin{equation*}
\dim \Nom{V} = 4r.
\end{equation*}
If we can determine the dimension of $\Nom{V}$, then we can limit our
search for four-weight spin models 
to the formally self-dual Bose-Mesner algebras having the right dimension.
In order to prove this, we need to better understand the conditions in
Lemma~\ref{lem_NV3} and the interactions
among $\Nom{A}$, $\N{A}{B}$ and $\dualN{A}{B}$.

Secondly, given an invertible Jones pair $(A,B)$, we get a link invariant.
But we also get a link invariant from the $4n \times 4n$ symmetric spin
models constructed from $(A,B)$.  It is natural to ask for the relations
between these two link invariants.

Further, we want to investigate more thoroughly the two constructions 
described in this chapter.
In particular, we would like to understand why the resulting Nomura algebras
relate to each other in such interesting fashion.  
In addition, we want to examine other known constructions for possible
extensions.
Nomura's construction of the symmetric and non-symmetric
Hadamard spin models, described in Section~\ref{Had_Spin},
is a good candidate.

In Section~\ref{JP_Dim2}, we have studied invertible Jones pairs of
dimension two.
The next open case is the invertible Jones pairs of dimension three.
On the other hand,
Jaeger characterized the three-dimensional Bose-Mesner algebras
that contain spin models.  
The characterization of symmetric Bose-Mesner algebras of dimension four that
contain spin models still remains open.
This case corresponds to a special class
of invertible Jones pairs of dimension four, so
any results on these Jones pairs may help solving this open case.

Lastly, we have designed Jones pairs in an attempt to answer Jones'
question in \cite{VJ_Knot}.  
However we have not been able to answer his question yet.
One obstacle is that there are very few 
Jones pairs where the matrices are not type II.  
Here we ask for new examples of non-invertible Jones pairs that are not
tensor products of existing Jones pairs.

\bibliographystyle{acm}
\bibliography{thesis}
\printindex
\end{document}